\documentclass[11pt]{amsart}
 \usepackage{graphicx} 
 \usepackage{amscd}

 \def\a{\alpha}
 \def\da{{\dot\alpha}}
 \def\be{\beta}
 \def\dbe{{\dot\beta}}
 \def\de{\delta}
 
 \def\e{\varepsilon}
 \def\deta{{\dot{\eta}}}
 \def\etanm{{\eta_{n,m}}}
 \def\detanm{{\dot{\eta}_{n,m}}}
 \def\teta{{\tilde{\eta}}}
 \def\ga{\gamma}
 \def\dga{{\dot{\gamma}}}
 \def\hga{{\widehat{\gamma}}}
 \def\tga{{\tilde{\gamma}}}
 \def\Ga{\Gamma}
 
 \def\vr{\varphi}
 \def\la{\lambda}
 
 \def\La{\Lambda}
 \def\si{\sigma}
 \def\Si{\Sigma}
 \def\om{\omega}
 \def\Om{\Omega}

 \def\dxi{{\dot{\xi}}}
 \def\xinm{{\xi_{n,m}}}
 \def\dxinm{{\dot{\xi}_{n,m}}}
 \def\dzeta{{\dot{\zeta}}}

 \def\re{{\mathbb R}}
 \def\na{{\mathbb N}}

 \def\ov{\overline}
 \def\Z{{\mathbb Z}}

 \def\cA{{\mathcal A}}
 \def\cAk{{{\mathcal A}_k}}

 \def\tB{{\widetilde{B}}}
 
 \def\cC{{\mathcal C}}
 \def\co{{\mathbb C}}
 \def\D{{\mathbb D}}
 \def\dnm{{d_{n, m}}}
 \def\E{{\mathbb E}}
 \def\Ek{{E^{-1}\{k\}}}
 \def\cF{{\mathcal F}}
 \def\H{{\mathbb H}}
 \def\cH{{\mathcal H}}
 \def\Hk{{H^{-1}\{k\}}}
 \def\Hx{{{\mathbb H}_x}}
 \def\cK{{\mathcal K}}
 \def\uk{{\underline{k}}}
 \def\L{{\mathbb L}}
 \def\cF{{\mathcal F}}
 
 \def\cL{{\mathcal L}}
 \def\fL{{\mathfrak L}}
 \def\hL{{\widehat{L}}}
 \def\tL{{\widetilde{L}}}

 \def\hM{{\hat{M}}}
 \def\tM{{\widetilde{M}}}
 \def\cN{{\mathcal N}}
 
 \def\P{{\mathbb P}}
 \def\cS{{\mathcal S}}
 \def\T{{\mathbb T}}
 \def\bT{{\mathbb T}}
 \def\cT{{\mathcal T}}
 \def\cU{{\mathcal U}}
 \def\tU{{\widetilde{U}}}
 
 \def\V{{\mathbb V}}
 \def\cW{{\mathcal W}}

 \def\tqo{{\tilde{q}_0}}
 \def\tq1{{\tilde{q}_1}}
 
 \def\dwn{{\dot{w}_n}}
 \def\tw{{\tilde{w}}}
 \def\dx{{\dot{x}}}
 
 \def\tx{{\tilde{x}}}
 \def\txn{{\tilde{x}_n}}
 \def\dxn{{{\dot x}_n}}

 \def\yn{{y_n}}
 \def\ty{{\tilde{y}}}
 \def\tyn{{\tilde{y}_n}}
 \def\dyn{{{\dot {y}}_n}}
 \def\dy{{\dot{y}}}
 \def\dz{{\dot{z}}}

 \def\tz{{\tilde{z}}}
 \def\tzn{{\tilde{z}_n}}

 \def \lv{\left\vert}
 \def \rv{\right\vert}
 \def \lV{\left\Vert}
 \def \rV{\right\Vert}
 \def \ov{\overline}

 \DeclareMathOperator*{\tsum}{{\textstyle \sum}} 
 \DeclareMathOperator{\supp}{supp}
 
 \DeclareMathOperator{\diam}{diam}
 \DeclareMathOperator{\length}{length}
 
 \DeclareMathOperator{\interior}{int}
 \DeclareMathOperator{\ima}{im} 

  \renewcommand{\proofname}{{\bf Proof:}}

 \theoremstyle{plain}

 \newtheorem{MainThm}{Theorem}
 \newtheorem{MainCor}[MainThm]{Corollary}
 \newtheorem{MainProp}[MainThm]{Proposition}

 \newtheorem{NThm}{\bf Theorem}
 \newtheorem{NProp}[NThm]{\bf Proposition}

  \swapnumbers

 \newtheorem{Thm}{Theorem}[section]
 
 \newtheorem{Lemma}[Thm]{\bf Lemma}
 
 \newtheorem{Corollary}[Thm]{\bf Corollary}
 \newtheorem{Theorem}[Thm]{\bf Theorem}
 \newtheorem{Proposition}[Thm]{\bf Proposition}
  
 \newtheorem{eexample}[Thm]{\bf Example}

 \theoremstyle{definition}

 \theoremstyle{remark}
 \newtheorem{Remark}[Thm]{\bf Remark}

 \newtheoremstyle{Cl}
  {5pt}
  {3pt}
  {\sl}
  {}
  {\it}
  {:}
  {.5em}
  {}

 \theoremstyle{Cl}
 \newtheorem*{Claim}{Claim}

 \def\begincproof{
                  \renewcommand{\proofname}{\it Proof:}
                  \begin{proof}
                 }

 \def\endcproof{
                \renewcommand{\qedsymbol}{$\diamondsuit$}
                \end{proof} 
                \renewcommand{\qedsymbol}{\openbox}
                \renewcommand{\proofname}{\bf Proof:}
               }

 \renewcommand{\proofname}{{\bf Proof:}}

 \title[The Palais-Smale condition]
 {The Palais-Smale condition on contact type energy levels 
 for convex lagrangian systems} 
 \author[G. Contreras]{Gonzalo Contreras}

\address{CIMAT \\
          A.P. 402, 3600 \\
          Guanajuato. Gto. \\
          M\'exico}
\email{gonzalo@cimat.mx}

\thanks{Partially supported by Conacyt, Mexico, grant  36496-E}

\begin{document}

   \begin{abstract} We prove that for a uniformly convex 
                    Lagrangian system $L$ on a compact manifold $M$,
		    almost all energy levels contain a periodic
		    orbit. We also prove that  below  Ma\~n\'e's 
		    critical value of the lift of the Lagrangian 
		    to the universal cover, $c_u(L)$, almost all energy levels
		    have conjugate points.

		    We prove that if the
		    energy level $[E=k]$ is of contact type and $M\ne
		    \T^2$ then the free time action functional  
		    of $L+k$ satisfies the Palais-Smale condition. 
		    
   \end{abstract}

 \maketitle


 \section{Introduction}

 In this paper we continue the study of the Morse theory 
 of the free time action functional for convex lagrangian    
 systems that we begun in~\cite{CIPP2}.			     
 This time we try to include the case of low energy levels, 
  where very little is known.				     
 The main problem with the free time action functional
 is that it may fail to satisfy the Palais-Smale condition,
 usually required for variational methods. 
 Here we prove that if an energy level is of contact type
 and the configuration space $M\ne \T^2$ is not the 2-torus,
 then it satisfies the Palais-Smale condition.
 We also prove that when an energy level projects onto the
 whole configuration space $M$, the set of closed loops
 has a mountain pass geometry.  An adaptation of an argument by
 Struwe to the mountain pass geometry shows the existence of
 convergent Palais-Smale sequences for almost all energy levels.
 This implies that for almost all energy levels which project 
 onto $M$ the Euler-Lagrange flow has a periodic orbit,
 has closed orbit loops starting at any $x\in M$, and has
 conjugate points if the energy is below Ma\~n\'e's 
 critical value of the universal cover.
 The same holds for an energy level which satisfies
 the Palais-Smale condition, and hence in particular
  for contact type energy levels. 
 
 In~\cite{CIPP2} we proved that high energy levels have a 
 periodic orbit. Very low energy levels which do not project
 onto $M$ are displaceable, and then,
  by results of Frauenfelder and Schlenk ~\cite{FraSch},
  \cite{Schlenk}, they have finite Hofer-Zehnder capacity.
 Combining these results we get that almost all energy
 levels have a periodic orbit. Our class of Lagrangian
 systems include  exact magnetic flows on compact manifolds.

 \subsection{Critical energy values}\quad

 Let $M$ be a closed Riemannian  manifold
 with $\dim M\ge 2$. 
 Let $\pi:TM\to M$ be the projection.
 A lagrangian on $M$ is a $C^\infty$ function $L:TM\to\re$.
 We shall assume that
  $L$ is (uniformly) {\it convex:} 
   \quad  there is $a>0$ such that
   $$
   w^*\cdot\tfrac{\partial^2L}{\partial v\;\partial v}\Big\vert_{(x,v)}
    \cdot w  > a\;\lv w\rv_x^2
  \quad\text{ for all }x\in M,\; v,w\in T_xM.
   $$
  This uniform convexity and the compactness of $M$ imply 
  (see e.g. Lemma~\ref{Lo} below) that  
 $L$ is superlinear:
     $$
     \lim_{|v|_x\to+\infty}\frac{L(x,v)}{|v|_x}=+\infty
     \quad\text{ uniformly on }TM.
     $$ 
 
 Since $M$ is compact and $L$ is autonomous,
 the {\it Euler-Lagrange equation}
 \begin{equation}\label{EL}
  \frac{d}{dt}\frac{\partial L}{\partial v}(x,\dx) 
  =\frac{\partial L}{\partial x}(x,\dx)
  \tag{E-L}
 \end{equation} 
 defines a complete flow $\vr_t$  on $TM$ called the {\it Euler-Lagrange flow}
 of $L$.
 The {\it energy function} $E:TM\to\re$,
 $$
 E(x,v):=\frac{\partial L}{\partial v}(x,v)\cdot v -L(x,v),
 $$
 is invariant under the Euler-Lagrange flow.
 
 The {\it action} of an absolutely continuous curve $\ga\in C^{ac}([a,b],M)$
 is defined by
 $$
 A_L(\ga)=\int_a^b L(\ga(s),\dga(s))\;ds.
 $$
 As noticed by Dias Carneiro~\cite{Carneiro} and Ma\~n\'e~\cite{Ma7},
 critical points for the action of $L+k$ among  curves 
 with free time interval are solutions of the Euler-Lagrange 
 equation which have energy $E\equiv k$.
 The most direct way to obtain critical points is to look for 
 minima. It turns out that if $k$ is low enough there are no minima
 because then the action of $L+k$ is not bounded from  below. The
 exact threshold is given by {\it Ma\~n\'e's critical value}:
 $$
 c(L)=\min\{\,k\in\re\,|\,A_{L+k}(\ga)\ge 0
 \text{ for all closed curves $\ga$ on $M$}\}.
 $$
 The action functional $A_{L+k}$ is bounded from below
 on the space of curves with fixed endpoints and on the
 space of closed curves if and only if $k\ge c(L)$.
   It is also known~\cite{Ma7},~\cite{CDI} that
 $$
 c(L)\ge e_0(L):=\min\{\,k\in\re\,|\,\pi(E^{-1}\{k\})=M\,\}.
 $$

 If $p:N\to M$ is a covering map and $L_1=L\circ dp$ is
 the lift of the lagrangian, it is easy to check that 
 $c(L_1)\le c(L)$.
 Thus we have that
 $$
 e_0(L) \le c_u(L) \le c_0(L) \le c(L),
 $$
 where $c_u$ and $c_a=c_0$ are the critical values of the
 lifts of $L$ to the universal cover  and the abelian 
 cover. The number $c_0(L)$ is also called the {\it strict critical value}
 and has the following characterization~\cite{PP4}:
 \begin{align}
 c_0(L)&= -\min\Big\{\,\int L\;d\mu\,\Big|\, \mu 
 \text{ is a  $\vr_t$-invariant probability with homology }
 \rho(\mu)=0\,\Big\}
 \notag
 \\
 &=\min\big\{\,c(L-\om)\;\big\vert\;[\om]\in H^1(M,\re)\,\},
 \label{c0}
  \end{align}
 where the homology $\rho(\mu)\in H_1(M,\re)\approx H^1(M,\re)^*$ of an
 invariant measure with compact support $\mu$
 is defined by
 $$
 \langle [\om],\rho(\mu)\rangle 
 = \int_{TM} \om_x(v)\;d\mu(x,v)
 $$
 for any closed 1-form $\om$ on $M$.
 Here $[\om]\in H^1(M,\re)$ is the
 cohomology class of $\om$.

 Given a covering map $p:N\to M$
 let $L_1=L\circ dp$ be the lift of the
 Lagrangian $L$ to $TN$ and $c_1=c(L_1)$ 
 its critical value.
 The {\it Peierls barrier} 
 $h_{c_1}:N\times N\to \re$ is defined by
 \begin{align*}
   h_{c_1}(q_0,q_1):&=\liminf_{T\to+\infty}\Phi_{c_1}(q_0,q_1;T),
   \\
  \Phi_{c_1}(q_0,q_1;T):&=\inf\,\{\,A_{L_1+c_1}(\ga)\,|\,
 \ga\in C^{ac}([0,T],N), \;\ga(0)=q_0,\;\ga(T)=q_1\;\}.
 \end{align*}

%
 \bigskip

 \bigskip

 \bigskip

 \subsection{The Palais-Smale condition}\quad
 
 We describe now our setting for the Morse theory of the 
 free time action functional. 
 Let  $\cH^1(M)$ be the set of 
 absolutely continuous curves  $x:[0,1]\to M$ such that
 $$
 \int_0^1\lv\dx(s)\rv_{x(s)}^2\;ds <\infty.
 $$
 Then $\cH^1(M)$ is a Hilbert manifold and its tangent space
 at $x$ consists of weakly differentiable vector fields 
 along $x$ whose  covariant derivative is bounded in $\cL^2$.
 We shall use the Hilbert manifold $\cH^1(M)\times \re^+$
 with the Riemannian metric
  \begin{equation}\label{1metric}
 \big\langle\,(\xi,\a),(\eta,\be)\,\big\rangle_{(x,T)}
 =\a\,\be+f(T)\,\big\langle\xi(0),\eta(0)\big\rangle_{x(s)}
 + g(T)\,\int_0^1\left\langle\tfrac{D}{ds}\xi(s),
 \tfrac{D}{ds}\eta(s)\right\rangle_{x(s)}\;ds,
 \end{equation}
 where $\tfrac{D}{ds}$ is the covariant derivative along $x(s)$
 and  $f$, $g:\re^+\to\re^+$ are smooth positive functions such that
 $\max\{f,\,g\}\le 2$,
 $$
 f(T)=\begin{cases} 
 T^2 & \text{ if } T\le 1,
 \\
   1  & \text{ if } T\ge 10.
 \end{cases} 
 \qquad\text{ and }\qquad
 g(T)=\begin{cases} 
 T^2 & \text{ if } T\le 1,
 \\
 \tfrac1T \;{e^{-4 T^2}} & \text{ if } T\ge 10.
 \end{cases} 
 $$ 
 We shall discuss this choice of metric in more detail later on.
 Observe that this metric is locally equivalent to the metric
 obtained when $f\equiv g\equiv 1$. In particular, the set of 
 differentiable functions on $\cH^1(M)\times \re^+$
 is the same for this metric and for the one with
  $f\equiv g\equiv 1$.
 
 Given $k\in\re$ define the {\it free time action functional} 
 $\cA_k:\cH^1(M)\times\re^+\to\re$ by
 $$
 \cA_k(x,T)=
 \int_0^1 \left[L\big(x(s),\tfrac{\dx(s)}T\big)+k\right]T\;ds.
 $$
 Observe that if $y(t):=x(t/T)$ then
 $$
  \cA_k(x,T)=A_{L+k}(y).
 $$

 We say that $L$ is {\it Riemannian at infinity} if there exists $R>0$ 
 such that $L(x,v)=\tfrac 12\,|v|_x^2$ if $|v|_x>R$.
 In~\cite[prop. 18]{CIPP2} it is proven that given a uniformly convex 
 lagrangian $L$ and $k\in\re$, there exists a convex lagrangian $L_0$
 such that $L=L_0$ on $[E\le k+1]$ and $L_0$ is Riemannian at infinity.
 In~\cite[Lemma 19]{CIPP2} it is proven that if $L=L_0$ on 
 $[E\le c(L)+1]$ then $c(L)=c(L_0)$. Thus if our objective is to
 find solutions of the Euler-Lagrange equation with prescribed energy,
 we can assume that $L$ is Riemannian at infinity.

 Given $q_0,\;q_1\in M$ let $\Om_M(q_0,q_1)$ be the set of curves 
 $(x,T)\in\cH^1(M)\times \re^+$ with endpoints $x(0)=q_0$ and $x(1)=q_1$.
 Also, let  $\La_M$ be the set of closed curves in $\cH^1(M)\times\re^+$.
 The sets  $\Om_M(q_0,q_1)$ and $\La_M$ are Hilbert submanifolds of 
 $\cH^1(M)\times\re^+$. A connected component of $\Om_M(q_0,q_1)$ (resp.
 $\La_M$) consists of closed  curves in the same homotopy class
 with fixed endpoints (resp. in the same free homotopy class).

 A theorem of Smale~\cite{Sm} implies that $\cA_k$ is $C^2$
 on $\cH^1(M)\times\re^+$ for the metric with $f\equiv g\equiv 1$,
 and hence also for the metric~\eqref{1metric}.
 We show in Lemma~\ref{solution} that a critical point
 of $\cA_k$ restricted to $\Om_M(q_0,q_1)$ or to $\La_M$
 is a solution of the Euler-Lagrange equation with energy
 $E\equiv k$. 
 
 We say that $\cA_k$ satisfies the {\it Palais-Smale condition
 on $\Om_M(q_0,q_1)$} [resp. on $\La_M$]
 or that the {\it energy level $k$ satisfies the Palais-Smale
 condition on $\Om_M(q_0,q_1)$} [resp. on $\La_M$]
 if every sequence $(x_n,T_n)$ in {\it the same connected component}
 of $\Om_M(q_0,q_1)$ [resp. on $\La_M$] such that
 $\lv \cA_k(x_n,T_n)\rv$ is bounded and 
 $\lim_n\lV d_{(x_n, T_n)}\cA_k\rV_{(x_n,T_n)}=0$
 has a convergent subsequence.

 We shall prove 

 \medskip

 \noindent
 \begin{MainThm}\label{NPS}\quad

  If $L$ is Riemannian at infinity 
 and $\cA_k$ does not satisfy the Palais-Smale condition
 on $\Om_M(q_0,q_1)$,  or on $\La_M$,
 then there exists a Borel probability measure $\mu$, invariant
 under the Euler-Lagrange flow, supported in a connected
 component of the energy level $E\equiv k$, 
 which has homology $\rho(\mu)=0$
 and whose $(L+k)$-action is zero:
 $$
 A_{L+k}(\mu)=\int \big[L+k\big]\;d\mu = 0.
 $$
 \end{MainThm}

 In appendix~\ref{Reeb} we give an example in which
 the measure obtained in Theorem~\ref{NPS} can not be ergodic.
 In~\cite[th.~C]{CIPP2} we found counterexamples
 to the Palais-Smale condition at $k=c(L)$, but in~\cite{CIPP2}
 we didn't require the Palais-Smale sequences to be in the same
 connected component of the space of curves. 
 Combining the arguments in~\cite{CIPP2} with those
 of Theorem~\ref{NPS} we get the following Corollary~\ref{SPS}.
 The novelty is that it allows curves with trivial homotopy class.

  \begin{MainCor}\label{SPS}\quad
    If $L$ is Riemannian at infinity then
   $\cAk$ satisfies the
  Palais-Smale condition for all $k>c_u(L)$.
   On $\Om_M(q_0,q_1)$,
   $\cA_{c_u}$ satisfies the Palais-Smale condition
  if and only if the Peierls barrier on the universal cover
  is $h_{c_u}\equiv+\infty$.
 \end{MainCor}

 Another example is the lagrangian
 $\L:T\D\to\re$ on the hyperbolic disc $\D \subset\co$, 
 where $\L(x,v)=\tfrac 12\,|v|_x^2+\eta_x(v)$, $|\cdot|_x$ is the
 hyperbolic metric and $\eta$ is a
 1-form on $\D$ whose differential $d\eta$ is the hyperbolic area form.
 In this case the Peierls barrier at $k=c(\L)=c_u(\L)$ is finite
 (cf.~\cite[ex.~6.2]{wkam}) 
 and $\cA_{c_u}$ does not satisfy the Palais-Smale condition.
 If $M$ is a compact surface with constant 
 curvature $K\equiv -1$, the Euler-Lagrange flow of  $\L$
 projects to a (non-exact) magnetic flow on $TM$. 
  At the energy level $k=c_u(\L)$ the projection of
  the Euler-Lagrange flow  of $\L$ 
   is the horocycle flow\footnote{When $k>c_u(L)$ the
   flow is Anosov on $d\pi(\E^{-1}\{k\})$ and for
   $k<c_u$ the energy level $d\pi(\E^{-1}\{k\})$
   is foliated by contractible periodic orbits.}
   for $M$ which has no closed orbits.

 The idea of the proof of Theorem~\ref{NPS} is the following.
 Let $(x_n,T_n)$ be a Palais-Smale sequence in the
 same connected component of $\La_M$ or $\Om_M(q_0,q_1)$.
 We first prove in Proposition~\ref{Tfinito}, 
 similar to~\cite{CIPP2}, that
 if the times $T_n$ are bounded away from 
 $0$ and $+\infty$ then there is a convergent subsequence.
 In Corollary~\ref{Om0} we prove that if $q_0\ne q_1$
 and $(x_n,T_n)\in\Om_M(q_0,q_1)$ then $T_n$ is bounded 
 away from zero. In Proposition~\ref{La0} we prove that if
 $\liminf_n T_n=0$ and $(x_n,T_n)\in\La_M$ or
 $(x_n,T_n)\in\Om_M(q_0,q_0)$  then 
 $x_n$ converges to a singularity $(q_0,0)\in TM$ of the
 Euler-Lagrange flow with zero action $L(q_0,0)+k=0$
 and energy $k$. In this case the measure $\mu$ is the
 Dirac probability supported at the point $(q_0,0)$.

 The most delicate case is when $\lim_n T_n=+\infty$.
 Since the gradient of $-\cA_k$ at $(x_n,T_n)$ 
 converges to zero one expects that the curves 
 $y_n(sT_n):=x_n(s)$ are approximate solutions of the 
 Euler-Lagrange equation with average energy $k$. 
 If $\mu_n$ is the
 probability measure on $TM$ defined by
 $$
 \int_{TM} f\; d\mu_n= \int_0^{1}
 f(x_n,\tfrac{\dxn}{T_n})\, ds
 =\frac1{T_n}\int_0^{T_n}f\big(y_n,\dyn)\; dt,
 \qquad \Big[ y_n(sT_n):=x_n(s)\Big],
 \label{mun}
 $$
 we prove that $\mu_n$
 converges to an invariant probability for the Euler-Lagrange flow
 with support in the energy level $k$.
 Since the $L+k$ action of the curves $y_n$ is bounded
 and $\lim_n T_n= +\infty$
 their average action converges to zero.
 Since their homotopy class is fixed, and $\lim_n T_n=+\infty$,
 their average homology class tends to zero.

 We use the functions $f$ and $g$ in the definition of the
 metric~\eqref{1metric} to deal with the cases 
 $\lim_nT_n=0$ and $\lim_n T_n=+\infty$.
 In order to justify their choice observe that
 by suitably expanding the metric near the endpoints
 any bounded function on the open interval $]-1,1[$ can be 
 made not to satisfy the  Palais-Smale condition.
 For example let $\psi(x)=x^2$ on $|x|<1$.
 Let $h:\re\to]-1,1[$ be a diffeomorphism.
 Then $\psi\circ h$ does not satisfy the Palais-Smale condition
 because $\lim_{x\to\pm\infty}d_x(\psi\circ h)=0$.
 So, if one is going to obtain any conclusion from the
 fact that the Palais-Smale condition does not hold,
 one needs to  use an appropriate metric.
 Since our metric is locally equivalent to the usual one
 with $f= g\equiv 1$, the critical points are still
 solutions of the Euler-Lagrange equation and also the
 change of metric does not prevent finding Palais-Smale
 sequences by, say, a minimax argument.
 
 \bigskip

 \subsection{\bf The mountain pass geometry.}\quad
 \medskip

  We show that for low energy levels $e_0(L)<k<c_u(L)$,
  the action functional $\cA_k$ exhibits a mountain pass geometry
  on the space of loops $\Om_M(q_0,q_0)$ and closed curves $\La_M$.
  This result is suggested by Ta\u{\i}manov in~\cite[p. 362]{tai1}
  for a different action functional for magnetic flows 
  saying that ``one-point curves form the
  manifold of local minima of the functional $\ell$''.
  S. Bolotin (cf.\cite[p. 362]{tai1})
  observed that the results of the papers~\cite{nov8}, \cite{nov7},
  \cite{nov9}, \cite{NovSh},
    \cite{NovTai}, \cite{tai10}
    may not be valid because the Palais-Smale condition could fail.
  The approach in this paper recovers the (a.e.)--validity of some
  of those results.
  
  Let $k<c_u(L)$. 
  By the definition of $c_u(L)$,
   there are a closed curve $(x_1,T_1)\in\La_M$ and
  for any $q_0\in \pi(E^{-1}\{k\})$ a loop
  $(x_2,T_2)\in\Om_M(q_0,q_0)$,
  both with trivial homotopy class and negative $(L+k)$-action.

  
  \begin{MainProp}
     \label{CritVal}\quad
    \begin{enumerate}
    \item\label{CVloop}
 Let $q_0\in M$ and $k>E(q_0,0)$.
 Then there exists $c>0$ such that
 if $\Ga:[0,1]\to \Om_M(q_0,q_0)$ is a continuous
  path joining a constant loop 
  $\Ga(0)=q_0:[0,T]\to \{q_0\}\subset M$ (with any $T>0$)
  to any closed loop $\Ga(1)\in\Om_M(q_0,q_0)$ 
  with negative $(L+k)$-action, $A_{L+k}(\Ga(1))<0$, then
 $$
 \sup_{s\in[0,1]} A_{L+k}(\Ga(s)) > c >0.
 $$
 
   \item\label{CVper}
 Let $k>e_0(L)$.  Then there exists $c>0$ such that
 if $\Ga:[0,1]\to \La_M$ is a continuous path joining any constant curve
 $\Ga(0)=q_0:[0,T]\to \{q_0\}\subset M$ to any closed curve
 $\Ga(1)$ with negative $(L+k)$-action, $A_{L+k}(\Ga(1))<0$, then
 $$
 \sup_{s\in[0,1]} A_{L+k}(\Ga(s)) > c >0.
 $$
 \end{enumerate}
 \end{MainProp}

 Standard critical point theory gives
 contractible periodic orbits on
 any energy level $e_0(L)<k<c_u(L)$
 where the Palais-Smale condition holds.
 Since the failure of the Palais-Smale 
 condition can only be due to one direction
 of non-compactness, namely the time parameter
 $T$ on $\cH^1(M)\times\re^+$,
 an argument originally due to Struwe in~\cite{Struwe1}
 (see also Struwe~\cite{StruweB}, Jeanjean~\cite{Jeanjean}
 and Jeanjean, Toland~\cite{JT}) can be applied
 to the mountain pass geometry of Proposition B 
 to overcome the Palais-Smale condition for almost every $k$.

 Previous results on higher energy
 levels (cf.~\cite{Ma7}, \cite{CDI}, \cite{CIPP2})
 give that $E^{-1}\{k\}$ has a periodic orbit for every
 $k>c_u(L)$. When the energy level does not project
 onto the whole configuration space $M$ (i.e. $k<e_0(L)$)
 we show that the displacement energy of 
 $[E\le k]$ is finite. Then by results of 
 U. Frauenfelder and F. Schlenk~\cite{FraSch}, \cite{Schlenk}, 
 the $\pi_1$-sensitive Hofer-Zehnder capacity of $[E\le k]$ is finite
 and so standard arguments (cf. ~\cite{HZ}) show that almost any
 energy level $E^{-1}\{k\}$, $k<e_0(L)$ has a
 contractible periodic orbit. We summarize this in the
 following:

   \begin{MainThm}\label{LC}\quad
 \renewcommand{\theenumi}{\alph{enumi}}
 \begin{enumerate}
 \item\label{LC.c} There is a total Lebesgue measure set 
 $A\subset \re$ such that for all $k\in A$ either
 the energy level $E^{-1}\{k\}$ is empty or it contains
 a periodic orbit.
 
 Moreover, \begin{itemize}
  \item The set $A$ contains $]c_u(L),+\infty[$. 
  \item If $k<c_u(L)$ and $k\in A$ this periodic orbit is
  contractible.
  \item  If $e_0(L)<k<c_u(L)$ and $k\in A$  it has positive
 $(L+k)$-action.
 \end{itemize}
 
 \item \label{LC.l}
 For any $q_0\in M$, there is a total Lebesgue measure
 subset $]c_u(L),+\infty[\subset B\subset ]E(q_0,0),+\infty[$
 such that for all $k\in B$ there is a solution of the Euler-Lagrange
 equation in $\Om_M(q_0,q_0)$ with energy $k$.

 \item \label{LC.ps}
 The above items hold for a specific $k\in]e_0(L),c_u(L)[$ (resp.
 $k\in]E(q_0,0),c_u(L)[$) if the energy level $k$ satisfies the
 Palais-Smale condition.
 
 \end{enumerate}
 \end{MainThm}
 
  As an example in appendix~\ref{non-magnetic} we prove that
 a lagrangian with no magnetic term has a closed
 orbit on every energy level.

 \bigskip

 Two points $\theta_0,\,\theta_1\in TM$, are said to be 
 {\it conjugate} \label{conjugate} 
 if there is $\tau\in\re$ such that 
 $\theta_1=\vr_\tau(\theta_0)$ and 
 $\V(\theta_1)\cap d_{\theta_0}\vr_\tau(\V(\theta_0))\ne\{0\}$,
 where $\V\subset T(TM)$ is the vertical sub-bundle
 $\V(\theta)=\ker d_\theta\pi$.
 R. Ma\~n\'e asked whether if $k<c_0(L)$ there is always an orbit
 with energy $k$ and conjugate points. G. Paternain and M. Paternain
 in~\cite{PP4} showed examples of magnetic  flows with Anosov energy 
 levels without conjugate points with energy $k\in]c_u(L),c_0(L)[$.
 At $k=c_u(L)$ these examples do not have conjugate points.
 The question remains open for $k<c_u(L)$.

 In~\cite[p. 663]{CIPP2} we gave an example of an orbit  segment
 without conjugate points which is not a local minimizer of
 the free time action functional.
 In Proposition~\ref{SLM} we prove
 that in an energy level without conjugate points
 every orbit segment is a strict local minimizer of the action
 functional. Since a mountain pass critical point can not
 be a strict local minimizer we get

 \begin{MainThm}\label{CP}
 
 Let $e_m(L)=\inf_{(x,v)\in TM}E(x,v)$.

 There is an open subset with total Lebesgue measure
 $A\subset [e_m(L),c_u(L)[$
 such that if $k\in A$ then there is an orbit with energy $k$
 and conjugate points.
 
 If $e_m(L)<k<c_u(L)$ and $\cA_k$ satisfies the Palais-Smale
 condition, then the energy level $k$ has conjugate points.
 \end{MainThm}
 
 In~\cite[prop. 8]{CIPP} and in~\cite[cor.~1.13]{CI} 
 we proved that if $k$ is a regular value of the energy
 function $E$ and $k<e_0(L)$ then $E^{-1}\{k\}$ has conjugate points.

 We don't know if the following holds:
 
 {\it Question:} Is it true that for the universal cover $\tilde{M}$,
 $$
 c_u(L)=\inf\{ k\in\re\,|\, \forall x,y\in\tilde{M}
 \;\exists \text{ orbit }\ga\in\Om_M(x,y), \; E(\ga,\dga)=k\;\}\;?
 $$
 
 An {\it exact magnetic flow} is the lagrangian flow of
 $$
 L(x,v)=\tfrac 12\, |v|_x^2 - \eta_x(v),
 $$ where
 $|\cdot|_x$ is the Riemannian metric of $M$ and
 $\eta_x$ is a non-closed 1-form on $M$.
 Thus for exact magnetic flows we get periodic orbits 
 for almost all energy levels and in particular
 for contact type energy levels, as seen below.

  \bigskip

 \subsection{Contact type energy levels.}\quad
 
 We now concentrate on a property that ensures the
 Palais-Smale condition.
 Let $H:T^*M\to\re$ be the hamiltonian associated to $L$:
 \begin{equation}\label{H}
 H(x,p)=\max_{v\in T_xM}\big[ p(v)-L(x,v)\big],
 \end{equation}
 and let $\om=dp\wedge dx$ be the canonical symplectic form on
 $T^*M$.
 The {\it hamiltonian vector field} $X$ on $T^*M$ is defined by 
 $i_X \om=-dH$.
 The induced {\it hamiltonian flow} is conjugate to the 
 lagrangian flow of $L$ by the Legendre transform 
 $\cL:TM\to T^*M$, $\cL(x,v)=L_v(x,v)$.
 The energy function satisfies $E= H\circ \cL$, so that
 energy levels for $L$ are sent to level sets of $H$.

 An energy level $\Si=H^{-1}\{k\}$
 is said to be of {\it contact type}
 if there exists a 1-form $\la$ on $\Si$
 such that $d\la=\om\vert_{T\Sigma}$ and $\la(X)\ne 0$.
 We call such a form $\la$ a {\it contact-type} form 
 for $\Si$.

 \begin{MainProp}\label{kPS}
 If $[H=k]$ is of contact type, $\dim M\ge 2$  and
 \begin{itemize}
 \item $M\ne \T^2$ or
 \item $M=\T^2$ and $k\notin[e_0,c_0]$,
 \end{itemize}
  then $\cA_k$ satisfies the Palais-Smale condition.
 \end{MainProp}

  \bigskip

   In Section~\ref{sAF} we introduce the space
   of curves with free time interval and the action functional
   and compare various metrics on the space of curves.
   In Section~\ref{sPSC} we prove Theorem~\ref{NPS}.
   In Section~\ref{sSPS} we prove Corollary~\ref{SPS}.
   In Section~\ref{sMPG} we prove Proposition~\ref{CritVal}
    on the mountain pass geometry. In Section~\ref{sSRMT} we prove 
   some results in Morse theory that we need
   and the relative completeness of the gradient flow of
    the action functional.
   In Section~\ref{sPSMPG} we give the argument to overcome
   the Palais-Smale condition in a mountain pass geometry for the
   action functional. In Section~\ref{sLC} we prove
   Theorems~\ref{LC} and~\ref{CP} and in Section~\ref{CTO}
   we prove Proposition~\ref{kPS}.
   In Appendix~\ref{Reeb} we give an example in which
   the measure of Theorem~A can not be ergodic.
   In Appendix~\ref{nctp} we show energy levels of non-contact type.
   In Appendix~\ref{non-magnetic} we prove that non-magnetic
   lagrangians have periodic orbits on every energy level.   

   \bigskip
   
   The author wishes to thank Patrick Bernard that suggested the 
   possibility of using Struwe's argument in our situation.

\section{The action functional and the space of curves.}
 \label{sAF}
 Given a Riemannian metric on $M$, by Nash's Theorem there exists 
 an isometric embedding of $M$ into some $\re^N$.
 Let 
 $$
 \cH^1:=\cH^1(\re^N):=\Big\{\,\xi:[0,1]\to\re^N\text{ absolutely
 continuous }\;\Big|\;
 \int_0^1|\dxi(s)|^2\;ds<+\infty\;\Big\}
 $$ 
 be endowed with the metric
 \begin{equation*}\label{NH1}
 \langle\xi,\eta\rangle_{\cH^1}:=\langle\xi(0),\eta(0)\rangle
 +\int_0^1\langle\dxi(s),\deta(s)\rangle\;ds.
 \end{equation*}
 The corresponding norm is given by
 \begin{equation*}
 \lV\xi\rV_{\cH^1}^2:=\lv\xi(0)\rv^2+\int_0^1|\dxi(s)|^2\;ds.
 \end{equation*}
 On $\cH^1\times\re^+$ we shall use the Riemannian metric
 \begin{equation}\label{metric}
 \big\langle\,(\xi,\a),(\eta,\be)\,\big\rangle_{(x,T)}
 =\a\,\be+f(T)\,\langle\xi(0),\eta(0)\rangle
 + g(T)\,\int_0^1\langle\dxi(s),\deta(s)\rangle\;ds,
 \end{equation}
 where $f$ and $g$ are smooth positive functions such that
  $\max\{f,\,g\}\le 2$,
 $$
 f(T)=\begin{cases} 
 T^2 & \text{ if } T\le 1,
 \\
   1  & \text{ if } T\ge 10.
 \end{cases} 
 \qquad\text{ and }\qquad
 g(T)=\begin{cases} 
 T^2 & \text{ if } T\le 1,
 \\
 \tfrac1T \;{e^{-4 T^2}} & \text{ if } T\ge 10.
 \end{cases} 
 $$

 Let $\zeta(t):=\xi(t/T)$, $0\le t\le T$. Then
 $\dxi(t/T)=T\cdot\dzeta(t)$ and
 $$
 \int_0^1|\dxi(s)|^2\;ds=T\cdot\int_0^T|\dzeta(t)|^2\;dt.
 $$
 In the variables $(\zeta,\a)$ the Riemannian metric above
 is written as
 \begin{equation}\label{NormForm}
 \begin{aligned}
 \lV(\xi,\a)\rV_{(x,T)}^2
 &=\a^2+T^2\,|\zeta(0)|^2 +{T^3}\int_0^T|\dzeta|^2\;dt
 &\text{ if } T\le 1,
 \\
  \lV(\xi,\a)\rV_{(x,T)}^2
 &=\a^2+|\zeta(0)|^2 + {e^{-4T^2}}\int_0^T|\dzeta|^2\;dt
 &\text{ if } T\ge 10.
 \end{aligned}
 \end{equation}
 This metric is locally equivalent to the metric of the product
 Hilbert space $\cH^1\times\re$.
 
 Given $q_0,q_1\in M$, let
 \begin{align*}
 \Om(q_0,q_1):&=\big\{\,(x,T)\in\cH^1\times\re^+\;\big|\;
 x(0)=q_0,\;\;x(1)=q_1\,\big\},
 \\
 \La :&=\big\{\,(x,T)\in\cH^1\times\re^+\;\big|\;
 x(0)=x(1)\,\big\}.
 \end{align*}
 Their tangent spaces at $(x,T)$ are given by
 \begin{align*}
 T_{(x,T)}\Om(q_0,q_1)
 &=\big\{(\xi,\a)\in\cH^1\times\re\;\big|\;
 \xi(0)=\xi(1)=0\,\big\},
 \\
 T_{(x,T)}\La
 &=\big\{(\xi,\a)\in\cH^1\times\re\;\big|\;
 \xi(0)=\xi(1)\,\big\}.
 \end{align*}
 Endow $\Om(q_0,q_1)$ and $\La$ with the Riemannian metric~\eqref{metric}.

 Let
 \begin{gather*}
 \cH^1(M):=\big\{\; x\in\cH^1(\re^N)\;\vert\; x([0,1])\subset M\;\big\},
 \\
 \Om_M(q_0,q_1):= \Om(q_0,q_1)\cap\cH^1(M)\times\re^+,
 \\
 \La_M:=\La\cap \cH^1(M)\times\re^+.
 \end{gather*}
 Then $\cH^1(M)\times\re^+$, $\Om_M(q_0,q_1)$ and $\La_M$
 are Hilbert   submanifolds of 
 $\cH^1(\re^N)\times\re^+$, $\Om(q_0,q_1)$ and $\La$  
 respectively.
 A connected component of $\La_M$ is given by closed curves 
 in the same free homotopy class.
  A connected component of $\Om_M(q_0,q_1)$ is given by the curves $(x,T)$
 in $\Om_M(q_0,q_1)$ which have a given homotopy class with fixed
 endpoints. 

 On $\cH^1(M)\times\re^+$ we shall use the intrinsic Riemannian metric
 defined by
 \begin{equation}\label{metricM}
  \big\langle\,(\xi,\a),(\eta,\be)\,\big\rangle_{(x,T)}
 :=\a\,\be+f(T)\,\langle\xi(0),\eta(0)\rangle_{x(0)}
 + g(T)\,\int_0^1\left\langle\tfrac{D}{ds}\xi(s),
 \tfrac{D}{ds}\eta(s)\right\rangle_{x(s)}\;ds
 \end{equation} 
 where $\langle\cdot,\cdot\rangle_x$ is the Riemannian metric on
 $M$ and $\tfrac{D}{ds}$ are covariant derivatives.
 Since $M$ is isometrically embedded into $\re^N$,
 the covariant derivative
 $\tfrac{D}{ds}\xi(s) = \P\big(\dxi(s)\big)$
 is the orthogonal projection $\P:T_x\re^N\to T_xM$
 of the derivative $\dxi(s)$ taken in $\re^N$.
 Thus the norm in $\cH^1(M)\times\re^+$ is
 smaller than the induced norm from $\cH^1(\re^N)\times\re^+$. 
 Formulas analogous  to~\eqref{NormForm} hold for
 the norm on $\cH^1(M)\times\re^+$.
  We also want to compare the metric on $\cH^1(M)\times\re^+$
  with the metric induced by $\cH^1(\re^m)\times\re^+$ on a
 local chart $\re^m\supset U\hookrightarrow M$. 
 In Lemma~\ref{equivM} below we shall prove that the three norms
 are locally equivalent.

 \bigskip

 Given $k\in\re$, define the action functionals 
 $\cA_k:\Om_M(q_0,q_1)\to\re$ and 
 $\cA_k:\La_M\to\re$ of $L+k$ by
 $$
 \cA_k(x,T)=\int_0^1T\,\left[L\big(x(s),\tfrac{\dx(s)}T\big)+k\,\right]
 \;ds.
 $$
 Writing $y(t):=x\big(\tfrac tT\big)$, $0\le t\le T$, we have that 
 $$
 \cA_k(x,T)=\int_0^T\left[L(y,\dy)+k\,\right]\;dt =:A_{L+k}(y).
 $$

     We say that a  lagrangian $L$ is {\it quadratic at infinity }
 if there is $R>0$, a 1-form $\theta_x$ on $M$ and $a,\psi\in
 C^\infty(M,\re)$, $a>0$, such that
 $L(x,v)=\tfrac 12\,a(x)\, |v|_x^2+\theta_x(v)+\psi(x)$
 for all $|v|_x\ge R$,
 where $|v|_x$ is the Riemannian norm of $v$ in $TM$.

 We say that $L$ is {\it Riemannian at infinity} 
 if there exists $R>0$ such
 that $L(x,v)=\tfrac 12 |v|^2_x$ for all $|v|_x>R$.
 Since we are assuming that $M$ is isometrically embedded in $\re^N$,
 this is equivalent to $L(x,v)=\tfrac 12\,|v|^2$ for $|v|>R$, where
 $|v|$ is the euclidean norm of $v$ and $(x,v)\in
 TM\subset\re^N\times\re^N$.
 In a coordinate chart such a lagrangian is given as
  $L(x,v)=\tfrac 12\; v^* G(x) v$,
 when $|v|$ is large enough, where $G(x)$ is the matrix
 of the Riemannian metric in the chart.
 Then  in coordinate charts $L$ is  quadratic at infinity.
 
 It follows from a result of Smale~\cite{Sm} that if $L$ is Riemannian
 at infinity then the action functional $\cA_k$ is $C^2$ on
 $\cH^1(M)\times\re^+$ with the metrics with $f\equiv g\equiv 1$.
  Since the Riemannian  metrics~\eqref{metric}, \eqref{metricM}
 are locally equivalent to the metrics with $f\equiv g\equiv 1$,
 then $\cA_k$ is $C^2$ on $\Om_M(q_0,q_1)$
 and on $\La_M$, with respect to all three Riemannian metrics.

  The derivative of $\cA_k$ is given by
  \begin{equation}\label{deriv}
  \begin{aligned}
   d_{(x,T)}\cA_k(\xi,\a)
  &=\int_0^1T\,\left[\,L_x\big(x,\tfrac{\dx}T\big)\,\xi
  +L_v\left(x,\tfrac{\dx}T\right)\,\tfrac{\dxi}{T}\,\right]\;ds
 +\a\int_0^1\left[\,k-E\left(x,\tfrac{\dx}T\right)\,\right]\;ds
 \\
 &=\int_0^T\left[\,L_x(y,\dy)\,\zeta+L_v(y,\dy)\,\dzeta\,\right]\;dt
 +\frac{\a}T\int_0^T\big[\,k-E(y,\dy)\,\big]\;dt,
  \end{aligned}
\end{equation}
 where $y(t)=x\big(\tfrac tT\big)$, $\zeta(t)=\xi\big(\tfrac tT\big)$,
 for $0\le t\le T$ and $E:TM\to\re$, 
 $$
 E(x,v)=v\,L_v(x,v)-L(x,v)
 $$ 
 is the energy function.
 The formulas~\eqref{deriv} can be interpreted either in
 local charts with usual derivatives or in covariant derivatives.
 In the former case, 
 \begin{align*}
 L_x\,\xi =\langle \nabla_x L,\xi\rangle_{x(s)}
 \quad \text{ and }\quad
 L_v\,\dxi =\langle \nabla_v L,\tfrac{D}{ds}\xi\rangle_{x(s)},
 \end{align*}
 where  $\nabla_x L$ and $\nabla_v L$ are the projections of
 the gradient of $L$ to the splitting $T_{(x,\dx)}TM=H\oplus V$
 and $\tfrac{D}{ds}\xi$ is the covariant derivative of $\xi$.
 The splitting $T_\theta TM=H(\theta)\oplus V(\theta)$ is described
 on page~\pageref{splitting}.

  \bigskip

     Fix $C_1>0$ we say that $f:U\subset \re^m\to M$ is
    a {\it bounded chart} if $f$ is an embedding
    such that the pull-back $f^* g$ of the 
    Riemannian metric has matrix $G(x)$  such that
    $G$ and $G^{-1}$ have $C^1$ norm 
    bounded by $C_1>0$. Fix a finite atlas of bounded
    charts $\cU=\{U_i\}_{i\in I}$ such that each $U_i\subset\re^m$
    is a convex set.

    We fix some constants used repeatedly. 
    Observe that the property
    of being quadratic at infinity is invariant
    under transformations by bounded charts.
    The following constants are 
    taken to hold in any bounded chart of our finite atlas, 
    i.e. in equations
    ~\eqref{a_0}--\eqref{b1} below the same constants
    are assumed to hold when the norm $|\cdot|_x$ 
    is interpreted as either the riemannian metric on $M$,
    the euclidean norm on $\re^N\supset M$
    or the euclidean norm on any bounded chart $U_i\subset\re^m$
    of our atlas.
    The norms $|L_x|$, $|L_{xv}|$  are interpreted as the euclidean
    norm in any bounded chart $U_i\in\cU$.

  Since $L$ is convex, 
  \begin{equation}\label{a_0}
  a_0:=\inf_{(x,v)\in TM} \frac{v\cdot L_{vv}(x,v)\cdot v}{|v|_x^2}>0.
  \end{equation}
   Since
  $L$ is quadratic at infinity there are $a_1,\;a_2>0$ such that
  \begin{equation}\label{a1a2}
  L(x,v)\ge a_1\,|v|^2_x-a_2,\qquad \text{ for all }(x,v)\in TM.
  \end{equation}
     \begin{gather}
   A_0:=\sup_{(x,v)} \lV L_{vv}(x,v)\rV_x<+\infty,
   \label{A0}
   \\
   b_1:=\sup_{x\in M}\lV d_x\psi\rV_x, \qquad
   \text{ where }\quad\psi(x):=L(x,0),
  \label{b1}
  \\
  \intertext{  Since  $L$ is quadratic at infinity,}
  b_2:=\sup_{(x,v)\in TM}\frac{\lv L_x(x,v)\rv}{1+|v|_x^2}<+\infty,
  \label{b2}
  \\
  b_3:=\sup_{(x,v)\in TM}\frac{\lv L_{xv}(x,v)\rv}{1+|v|_x}<+\infty.
  \label{b3}
  \end{gather}

  \bigskip

 \begin{Lemma}\label{solution}
 If $(x,T)\in\Om_M(q_0,q_1)$ {\rm [resp. $(x,T)\in\La_M$]}
 is a critical point of the action functional $\cA_k:\Om_M(q_0,q_1)\to\re$
 {\rm [resp. $\cA_k:\La_M\to\re$]} then the curve
 $y:[0,T]\to M$, $y(t):=x(t/T)$
 is a differentiable solution of the Euler-Lagrange equation
 $\tfrac{d\,}{dt}L_v(y,\dy)=L_x(y,\dy)$
 with $y(0)=q_0$, $y(T)=q_1$, {\rm [resp. $(y,\dy)$ is a closed orbit of
 the Euler-Lagrange flow]} with energy $E(y,\dy)\equiv k$.
 \end{Lemma}

  \begin{proof}
 Cover the image $y([0,T])$ by  images $\tU_i\subset M$ 
 of charts $U_i\in\cU$, $U_i\subset\re^m$.
 It is enough to prove that $y$ is a solution of
 the Euler-Lagrange equation on each intersection
 $y([0,T])\cap \tU_i$. Assume for a while that
 $y([0,T])\subset U_i\subset\re^m$.
 Using the same notation as in~\eqref{deriv}, we have that
 \begin{align}
 d_{(x,T)}\cA_k(\xi,0)
 &=\int_0^T \big[\,L_x(y,\dy)\;\zeta
 +L_v(y,\dy)\;\dzeta\,]\;dt 
 \label{derx}
 \\
 &=\L_x\cdot\zeta\big\vert_0^T
 +\int_0^T\big[\,L_v(y,\dy)-\L_x(t)\,]\cdot\dzeta\;dt
 \notag\\
 &=0, 
 \notag
 \end{align}
 where 
 $
 \L_x(t):=\int_0^t L_x\big(y(s),\dy(s)\big)\;ds.
 $
 Since $L_x(y,\dy)\le b_2(1+|\dy|_x^2)$ and $x\in\cH^1(M)$,
 we have that $L_x(y,\dy)\in\cL^1([0,T],\re^m)$ and that
 $\L_x$ is continuous in view of Lebesgue's theorem. 

 Since  for both $\La_M$ and $\Om_M(q_0,q_1)$ we can choose 
 $\zeta(0)=\zeta(T)=0$, 
 $$
 \int_0^T\big[\,L_v(y,\dy)-\L_x(t)\,]\cdot\dzeta\;dt
 =0
 $$
 for all $\dzeta\in\cL^2([0,T],\re^m)$ with $\int_0^T\dzeta\;dt=0$.
 
 This implies that $L_v(y,\dy)-\L_x$ is constant a.e. in $[0,T]$.
 Since $\L_x(t)$ is continuous,  it is bounded on $[0,T]$.
 Since $L$ is superlinear and $\L_x$ is bounded,
 $\dy$ is bounded by a constant almost everywhere. Since $L$ is convex,
 $v\mapsto L_v(y,v)$ is a continuous bijection. Hence we can uniquely
 extend $\dy$ to $[0,T]$ so that $L_v(y,\dy)-\L_x$ is constant on
 all $t\in[0,T]$.
  Since $\L_x(t)$ is continuous, $L_v(y,\dy)$ is also continuous
  and hence $\dy(t)$ is continuous.
  
 We have that 
 \begin{equation}\label{Lv}
 L_v\big(y(t),\dy(t)\big)=A+\int_0^tL_x(y,\dy)\;dt
 \end{equation}
 for some constant $A\in\re^m$. Since $\dy(t)$ is continuous, 
 the right hand side of~\eqref{Lv} is differentiable and
 $$
 \tfrac{d\;}{dt}\,L_v(y,\dy)=L_x(y,\dy).
 $$
 \noindent
  Hence $y(t)$ is a differentiable solution of the Euler-Lagrange
  equation. The theory of ordinary differential equations implies that
  $y$ is $C^r$ if $L$ is $C^{r+2}$.
  
  Since $y(t)$ is a solution of the Euler-Lagrange equation,
  its energy $E(y(t),\dy(t))$ is constant. Since
  $$
  \left.\frac{\partial \cA_k}{\partial T}\rv_{(x,T)}
  =\frac 1T \;\int_0^T \!\!\big[\,k-E(y,\dy)\,\big]\;dt =0,
  $$
  $E(y,\dy)\equiv k$. This completes the case of $\Om_M(q_0,q_1)$.
  
  For the case of $\La_M$ it remains to prove that $\dy(0)=\dy(T)$.
  Choose a chart $U_i\subset \cU$ whose image contains $y(0)=y(T)$
  and restrict ourselves to vector fields $\zeta$ over $y$
  with support in the connected component of $y([0,T])\cap \tU_i$
  containing $y(0)$.
   Since we already know that $(y,\dy)$ is a differentiable solution of the
  Euler-Lagrange equation, integrating by parts in~\eqref{derx} we
  have that
     \begin{align*}
     d_{(x,T)}\cA_k(\xi,0)     &=L_v\;\zeta\big\vert_0^T 
    +\int_0^T\left(L_x-\tfrac {d\,}{dt}L_v\right)\,\zeta\;dt
    \\
       &=\big[\,L_v\big(y(T),\dy(T)\big)-L_v\big(y(0),\dy(0)\big)\,\big]
    \cdot\zeta(0) + 0
   \end{align*}
   whenever $\zeta(0)=\zeta(T)\in\re^m$. 
   Then $L_v\big(y(T),\dy(T)\big)=L_v\big(y(0),\dy(0)\big)$.
   Since $y(T)=y(0)$ and  $v\mapsto L_v(y(0),v)$ is injective, 
   $\dy(T)=\dy(0)$.
    
    \end{proof}

\bigskip

 The following lemma shows that the intrinsic Riemannian metric 
 $\lV\,\cdot\,\rV^{\cH^1(M)\times\re^+}_{(x,T)}$ given by~\eqref{metricM}
 on $\cH^1(M)\times\re^+$
 and the induced metric from $\cH^1(\re^N)\times\re^+$ are locally
 equivalent. Also the metric on $\cH^1(M)\times\re^+$ and the
 metric on $\cH^1(\re^m)\times\re^+$, $m=\dim M$ on a bounded coordinate
 chart are locally equivalent.
 
 \begin{Lemma}\label{equivM}\quad
 \begin{enumerate}
 \item\label{equivRN} 
      Given $A_1, T_1>10$, there exists $B=B(A_1,T_1,k,\{f,g\})>0$ such that
 
      if $(x,T)\in\cH^1(M)\times\re^+$, $\lv\cA_k(x,T)\rv<A_1$ 
      and $T<T_1$, then 
      
      for all $(\xi,\a)\in
      T_{(x,T)}\Om_M(q_0,q_1)\cup T_{(x,T)}\La_M$,
      \begin{equation}\label{equiv1}
      \tfrac 1B \:
      \lV (\xi,\a)\rV_{(x,T)}^{\cH^1(\re^N)\times\re^+}
      \le
      \lV (\xi,\a)\rV_{(x,T)}^{\cH^1(M)\times\re^+}
      \le
      \lV (\xi,\a)\rV_{(x,T)}^{\cH^1(\re^N)\times\re^+}.
      \end{equation}
      
  \item\label{equivCH} 
       Let $\psi:U\subset\re^m\to M$ be an immersion such that
      the pull-back $\psi^*g_M(v,w)=v^*\,G(x)\,w$ of the Riemannian
      metric on $M$ has matrix $G(x)$ which is bounded in the
      $C^1$-norm:
      $$
      \max\left\{\lV G\rV_{C^1(U,\re^{m\times m})},
      \lV G^{-1}\rV_{C^1(U,\re^{m\times m})}\right\}<C_1.
      $$
      
      \noindent
      For all $A_1,T_1>10$ there exist $B=B(A_1,T_1,C_1,k,\{f,g\})>0$ such that 
         
      if $(x,T)\in\cH^1(\re^m)\times\re^+$,
      $\lv\cA_k(\psi\circ x,T)\rv<A_1$ and $T<T_1$, then
      
      for all $(d_x\psi\circ\xi,\a)\in
      T_{(x,T)}\Om_M(q_0,q_1)\cup T_{(x,T)}\La_M$,
      $$
      \tfrac 1B \:
      \lV (\xi,\a)\rV_{(x,T)}^{\cH^1(\re^m)\times\re^+}
      \le
      \lV (d_x\psi\circ\xi,\a)\rV_{(x,T)}^{\cH^1(M)\times\re^+}
      \le
       B\;\lV (\xi,\a)\rV_{(x,T)}^{\cH^1(\re^m)\times\re^+}.
      $$      
 \end{enumerate}
 \end{Lemma}

 \noindent{\bf Proof:}\quad 
 (1). Let $y(sT):=x(s)$. Then
 \begin{gather}
 A_1>A_{L+k}(y) \ge \int_0^T \big(a_1\,|\dy|_y^2-a_2+k\big)\;dt
 = a_1\int_0^T|\dy|_y^2\;dt - (a_2-k)\,T,
 \notag\\
 \int_0^T|\dy(t)|^2_y\;dt \le \frac{A_1+(a_2-k)\,T}{a_1},
 \label{yL2}
 \\
 \int_0^1\lv\dx(s)\rv_x^2\;ds = T\int_0^T\lv\dy(t)\rv_y^2\;dt
 \le T\left[\frac{A_1+(a_2-k)\,T}{a_1}\right].
 \label{xL2}
 \end{gather}
 Let $\bT(s,r):T_{x(r)}M\to T_{x(s)}M$ be the parallel 
 transport along $x(s)$. Then
 $$
 \xi(s)=\bT(s,0)\cdot \xi(0)
 +\int_0^s\bT(s,r)\cdot\tfrac{D}{dr}\xi(r)\;dr.
 $$
 Thus
 \begin{align}
 \lv\xi(s)\rv_{x(s)} 
 &\le \lv\xi(0)\rv_{x(0)}
 +\lV\tfrac{D}{ds}\xi\rV_{\cL^1([0,1])}
 \le \lv\xi(0)\rv_{x(0)}
 +\lV\tfrac{D}{ds}\xi\rV_{\cL^2([0,1])}
 \notag\\
 &\le \max\left\{\tfrac{1}{\sqrt{f(T)}},\tfrac 1{\sqrt{g(T)}}\right\}
 \; \sqrt{2}\;\lV(\xi,0)\rV_{(x,T)}^{\cH^1(M)\times\re^+}
 \notag\\
 &\le \sqrt{2}\;\max\big\{T^{-1},\a^{-\frac12},T_1^{\frac 12} {e^{2T_1^2}}\big\}
 \; \lV(\xi,0)\rV_{(x,T)}^{\cH^1(M)\times\re^+}
 \label{xiC0}
 \end{align}
 where
 $$
 \a:=\min_{1\le t\le 10}\big\{f(t),\;g(t)\big\}.
 $$
 
  Observe that  for $\xi, \eta\in T_xM$ 
  the first and second terms in~\eqref{metric},
  \eqref{metricM} are equal, so we only have to bound the
  $\cL^2$ norm of the derivatives  $\dxi$ and $\tfrac{D}{ds}\xi$.
  
  Let $U$ be a small tubular neighbourhood of $M$ in $\re^N$ and 
  let $F:U\to M$ be the orthogonal projection onto $M$.
  Since $\xi(s)\in T_{x(s)}M$, we have that
  $d_{x(s)}F\cdot\xi(s)=\xi(s)$.
  Differentiating this equation with respect to $s$ we get that
  $$
  d^2_{x(s)}F\big(\dx(s),\xi(s)\big)
  +d_{x(s)}F\cdot\dxi(s)=\dxi(s).
  $$
  The second term is the projection of $\dxi$ to $T_xM$:
  $$
  d_{x(s)}F\cdot\dxi(s)=\P\cdot\dxi(s)=\tfrac{D}{ds}\xi(s);
  $$
  and the first term is the projection of $\dxi$ to the orthogonal
  complement $T_xM^\perp$ of $T_xM$:
  $$
  \dxi^\perp(s):= d^2_{x(s)}F\big(\dx(s),\xi(s)\big).
  $$
  Let $c_1:=\sup_{x\in M}\lV d^2_xF\rV^2$, then using~\eqref{xL2} 
  and~\eqref{xiC0} we have that
  \begin{align}
  g(T)&\int_0^1 \vert \dxi^\perp(s)\vert^2\;ds
  \le g(T)\int_0^1 c_1 \,|\dx(s)|^2\,|\xi(s)|^2\;ds
  \le c_1\;g(T)\;\big\Vert\xi\big\Vert_\infty^2\int_0^1|\dx(s)|^2\;ds
  \notag\\
  &\le c_1\,
   \min\{2\, T^2,2\}\;
   2\,\max\big\{T^{-2},\a^{-1},T_1\,e^{4T_1^2}\big\}
 \;\left[ \lV(\xi,0)\rV_{(x,T)}^{\cH^1(M)\times\re^+}\right]^2
  T\left[\frac{A_1+(a_2-k)\,T}{a_1}\right]
  \notag\\
  &\le c_1\,B_1\;
  \left[\lV(\xi,0)\rV_{(x,T)}^{\cH^1(M)\times\re^+}\right]^2,
 \label{L2prod}
 \end{align}
 where
 $$
  B_1=B_1(A_1,a_1,a_2,T_1,k,\{f,g\}):= 
         4\,T_1\,\left[\frac{A_1+(a_2+|k|)\,T_1}{a_1}\right]
    \max\big\{\a^{-1},T_1\,e^{4T_1^2}\big\}.
 $$
 Observe that the bound $c_1\,B_1$ above holds for all $0<T<T_1$.
  
 Since $\dxi=\P\;\dxi+\dxi^\perp$, we have that
 \begin{align*}
 g(T)\;\big\Vert\dxi\big\Vert_{\cL^2}^2
 &\le g(T)\,\left[\,\big\Vert\P\;\dxi\big\Vert_{\cL^2}^2
 +\big\Vert\dxi^\perp\big\Vert_{\cL^2}^2\,\right]
 \\
 &\le g(T)\,\big\Vert\tfrac{D}{ds}\xi\big\Vert_{\cL^2}^2
 +c_1\, B_1\;\left[\lV(\xi,0)\rV_{(x,T)}^{\cH^1(M)\times\re^+}\right]^2.
 \end{align*}
Then, for all $0<T<T_1$,
 $$
 \left[\lV(\xi,\a)\rV_{(x,T)}^{\cH^1(\re^N)\times\re^+}\right]^2
 \le \big[1+c_1\,B_1\big]\;
 \left[\lV(\xi,\a)\rV_{(x,T)}^{\cH^1(M)\times\re^+}\right]^2.
 $$

 Since 
 $\lv\tfrac{D}{ds}\xi\rv=\big\vert\P\;\dxi\big\vert
 \le\big\vert\dxi\big\vert$, 
 then $\lV\tfrac{D}{ds}\xi\rV_{\cL^2}\le\big\Vert\dxi\big\Vert_{\cL^2}$. 
 This implies the second inequality in~\eqref{equiv1}.
 
 \medskip

 (2). We have that
  $$
 \tfrac{D}{ds}\xi
  =\dxi(s)+\tsum_{ijk}\,\Ga^k_{ij}(x(s))\;\dx_i(s)\;\xi_j(s)\;e_k,
  $$
 where the $\Ga^k_{ij}(x)$ are the Christoffel symbols for the 
 Riemannian metric of $M$ in the coordinate chart $\psi^{-1}$ and 
 $e_k$ is the $k$-th vector of the canonical basis of $\re^m$.
 Our hypothesis on $\psi$ implies that 
 $c_2=c_2(C_1):=m^2\,\sup_{ijk,x\in U }
 \big\vert \Ga^k_{ij}(x)\big\vert^2$ is finite.
 Then
 $$
 \lv\tfrac{D}{ds}\xi\rv\le\big\vert\dxi(s)\big\vert
 +\sqrt{c_2}\,\lV\xi\rV_\infty\,\lv\dx(s)\rv.
 $$
 Similar calculations as in~\eqref{xiC0} and~\eqref{L2prod}
 using $\dxi$ in~\eqref{xiC0} instead of the covariant derivative show that
 if $0<T<T_1$, then
 $$
 g(T)\,\lV\tfrac{D}{ds}\xi\rV_{\cL^2}^2
 \le 2\,g(T)\,\big\Vert\dxi\big\Vert_{\cL^2}^2
 + 2\,c_2\, B_2\,
 \left(\lV(\xi,0)\rV_{(x,T)}^{\cH^1(\re^m)\times\re^+}\right)^2,
 $$
 where  $B_2:=B_1(A_1,\ov{a}_1,\ov{a}_2,T_1,k,\{f,g\})$
 and $\ov{a}_1$, $\ov{a}_2$ are constants such that
 the inequality~\eqref{a1a2} holds in our coordinate system $\psi$
 for the euclidean metric in $U\subset\re^m$ instead of 
 the riemannian metric $|\cdot|_x$ on $M$.
 Then, if $0<T<T_1$, we have that
 $$
 \lV(\xi,\a)\rV_{(x,T)}^{\cH^1(M)\times\re^+}
 \le \sqrt{2}\;\big[\,1+c_2\,B_2\,\big]^{\frac 12}\;
 \lV(\xi,\a)\rV_{(x,T)}^{\cH^1(\re^m)\times\re^+}.
 $$

 Now write
 \begin{gather*} 
  \dxi(s) =\tfrac{D}{ds}\xi
 -\tsum_{ijk}\,\Ga^k_{ij}(x(s))\;\dx_i(s)\;\xi_j(s)\;e_k,
 \\
 \big\vert\dxi(s)\big\vert
 \le \lv\tfrac{D}{ds}\xi\rv
 +\sqrt{c_2}\,\lV\xi\rV_\infty\,\lv\dx(s)\rv.
 \end{gather*}
 The same calculations as in~\eqref{xiC0} and~\eqref{L2prod} give
  $$
  g(T)\,\big\Vert\dxi\big\Vert_{\cL^2}^2
 \le 2\, g(T)\,\lV\tfrac{D}{ds}\xi\rV_{\cL^2}^2
 + 2\,c_2\,
 B_2\,\left(\lV(\xi,0)\rV_{(x,T)}^{\cH^1(M)\times\re^+}\right)^2.
 $$ 
 And then
 $$
 \lV(\xi,\a)\rV_{(x,T)}^{\cH^1(\re^m)\times\re^+}
 \le \sqrt{2}\;\big[\,1+c_2\,B_2\,\big]^{\frac 12}\;
 \lV(\xi,\a)\rV_{(x,T)}^{\cH^1(M)\times\re^+}.
 \qquad\qed
  $$
 
 In the next lemma, we write $d$ for the distance
 $d_{\cH^1(M)\times\re^+}$ on $\Om_M(q_0,q_1)$ or on $\La_M$.
  \begin{Lemma}\label{nearby}\quad
 Given $T_0>0$ there exists $C=C(T_0)>0$ and $\e=\e(T_0)>0$ 
 such that if $T\in[\tfrac 1{T_0},T_0]$, 
 $(x,T),\,(y,S)\in\Om_M(q_0,q_1)\cup\La_M$
 and $d\big((x,T),(y,S)\big)<\e$
 then for the Hausdorff distance $d_H$ induced by the
 Riemannian metric, we have that
 $$
 d_H\big(x([0,1],y([0,1])\big)<
  C\;d\big((x,T),(y,S)\big).
 $$
    \end{Lemma}
 
  \begin{proof}
  Let
  $$
  A_3:=4\;\max\left\{\tfrac 1{f(t)},\;\tfrac1{g(t)}\;\Big|\;
  t\in\left[\tfrac 1{2 T_0}, 2 T_0\right]\;\right\},
  $$
  where $f(t)$ and $g(t)$ are as given in the definition of
  the Riemannian metric on $\cH^1(M)$.
  Let $0<\e_0<1$ be such that
  $$
  \tfrac 1{T_0}-2 \e_0>\tfrac 1{2T_0}
  \quad\text{ and }\quad
  T_0+2 \e_0< 2T_0.
  $$
  Let $0<\e=\e(T_0)<\e_0$ and write
   $\de:=d\big((x,T),(y,S)\big)<\e$.
     There is a curve $\Ga(\la)=(z_\la,T_\la)$, $\la\in[0,1]$,
    from $(x,T)$   to $(y,S)$ in  $\Om_M(q_0,q_1)$ or in $\La_M$
    such that
   \begin{gather*}
   \length(\Ga)=\int_0^1 \lV \tfrac {d}{d\la}\Ga(\la)\rV\;d\la
               < 2\de.
   \end{gather*}
   We can reparametrize $\Ga$ so that the norm of
   its tangent vector is constant:
   $$
   \lV\tfrac {d}{d\la}\Ga(\la)\rV^2
   =  \lv\frac{d\,T_\la}{d\la}\rv^2 +
   f(T_\la) \; \lv\frac{\partial
   z_\la(0)}{\partial\la}\rv_{z_\la(0)}^2
   + g(T_\la)\int_0^1 
    \lv\tfrac{D} {d\la}\dz_\la\rv^2_{z_\la(s)}\;ds
      <4\,\de^2.
   $$

   Since $|T_\la-T|\le d\big(\Ga(\la),(x,T)\big)\le 2\de<2\e_0$ then 
   $S=T_{\la=1},\,T_\la,\,T\in[\tfrac 1{2T_0},2T_0]$.
   Hence 
    \begin{align}
  \lv\frac{\partial
   z_\la(0)}{\partial\la}\rv_{z_\la(0)}^2
   &< A_3\;\de^2 \qquad\text{ for }\la\in[0,1],
   \label{pzl}\\
    \int_0^1 
    \lv\tfrac{D} {d\la}\dz_\la\rv_{z_\la(s)}^2\;ds
    &< A_3\;\de^2 \qquad\text{ for }\la\in[0,1].
    \notag
    \end{align}
    Let
    $$
     F(s):=\frac 12\int_0^1
    \lv\frac{\partial z_\la(s)}{\partial\la}\rv^2_{z_\la(s)}\;d\la.
    $$
    From~\eqref{pzl}, $|F(0)|<A_3\,\de^2$. We have that
    \begin{align*}
    F(s)-F(0)&=\int_0^s \tfrac{d}{ds} F(s)\;ds
    \\
    &=\int_0^s\int_0^1\left\langle\frac{D}{ds}
    \frac{\partial }{\partial\la} z_\la(s),
    \frac{\partial}{\partial\la}z_\la(s)\right\rangle_{z_\la(s)}\;d\la\;ds
    \\
    &=\int_0^s\int_0^1\left\langle\frac{D}{d\la}
     \dz_\la(s),
    \frac{\partial}{\partial\la}z_\la(s)\right\rangle_{z_\la(s)}
    d\la\;\;ds\;,
    \\
    \lv F(s)-F(0)\rv&\le
    \left[\int_0^1\int_0^1 \lv \frac{D}{d\la}
    \dz_\la(s)
    \rv^2\; ds\;d\la  \,   \right]^{\frac 12}
    \left[\int_0^s\int_0^1 
    \lv \frac{\partial
    z_\la(s)}{\partial\la}\rv^2_{z_\la(s)}
    \;d\la\;ds\,\right]^{\frac 12}
    \\
    &\le \sqrt{A_3\,\de^2}\; \left[ 2 \int_0^s F(s)\;ds\right]^{\frac 12}.
    \\
    F(s) &\le A_3\,\de^2 +\sqrt{2 A_3}\,\de
    \left[\int_0^s F(s)\;ds\right]^{\frac 12}.
    \end{align*}
   
    Write 
    $$
    u(s):=\left[\int_0^s F(s)\;ds\right]^{\frac 12}.
    $$
    Then
    $$
    \tfrac{d\;}{ds}\;u(s)^2
    \le A_3\,\de^2 +\de\,\sqrt{2 A_3}\;\,u(s).
    $$
    Let $t_0:=\sup\{\,t\in[0,1]\;|\;
    u(s)\le\sqrt{2A_3}\,\de\,(s+\tfrac 12),\;\forall s\in[0,t]\,\}$.
    Then
    \begin{align*}
    \tfrac{d}{ds} u(s)^2 
    &\le A_3\,\de^2+2\,A_3\,\de^2\,(s+\tfrac 12)
    \qquad\text{ if } s\in[0,t_0],
    \\
    &\le 2\, A_3\,\de^2\,(s+1)
    \qquad\qquad\text{ if } s\in[0,t_0].
    \\
    u(s)^2 &\le 2\, A_3\,\de^2\,(\tfrac 12 s^2+s),
    \\
    u(s) &\le \sqrt{2\, A_3}\,\de\,\sqrt{\tfrac 12 s^2+s}
    \qquad \text{ if } s\in[0,t_0].
    \end{align*}
    
    Since $\sqrt{\frac 12 s^2+s\;}<s+\tfrac 12$ for all $s>0$,
    we have that $t_0=1$. Hence, for all $s\in[0,1]$,
    \begin{align*}
    F(s) &\le A_3\,\de^2+\sqrt{2A_3}\,\de\,u(s)
    \\
    &\le A_3\,\de^2+2\,A_3\,\de^2\,\sqrt{\tfrac 32}
    \\
    &\le \tfrac 12\,C^2\,\de^2.
    \end{align*}

    We have that 
    \begin{align*}
    d_M(y(s),x(s))=d_M(z_1(s),z_0(s))
    \le \int_0^1 \lv\frac{\partial z_\la(s)}{\partial \la}\rv\; d\la
    \le \sqrt{2 \,F(s)} \le C\;\de.
    \end{align*}
    for all $s\in[0,1]$. This implies the lemma.
   \end{proof}

 \bigskip
 \bigskip

  \section{The Palais-Smale condition.}
  \label{sPSC}

 \bigskip
 
 In this section
  we are interested in the validity of the Palais-Smale condition
 for the action functional $\cA_k$ on a connected component $\Om_1$
 (resp. $\La_1$) of $\Om_M(q_0,q_1)$ (resp. $\La_M$).

 \bigskip

 \renewcommand{\theNThm}{\ref{NPS}}
 \begin{NThm}\quad

  If $L$ is Riemannian at infinity 
 and $\cA_k$ does not satisfy the Palais-Smale condition
 on $\Om_M(q_0,q_1)$,  or on $\La_M$,
 then there exists a Borel probability measure $\mu$, invariant
 under the Euler-Lagrange flow, supported in a connected component
 of the energy level $E\equiv k$, which has homology $\rho(\mu)=0$
 and whose $(L+k)$-action is zero:
 $$
 A_{L+k}(\mu)=\int \big[L+k\big]\;d\mu = 0.
 $$
 \end{NThm}

 \begin{proof}[\bf Proof of Theorem~\ref{NPS}:]\quad
 
 \smallskip

 Let $(x_n,T_n)$ be a sequence in a connected component $\Om_1$ of
 $\Om_M(q_0,q_1)$ (resp. $\La_1$ of $\La_M$) such that 
 $$
 \lv\cA_k(x_n,T_n)\rv<A_1 \qquad\text{and}\qquad
 \lV d_{(x_n,T_n)}\cA_k\rV<\tfrac 1n.
 $$
 Assume that $(x_n,T_n)$ does not have an accumulation point 
 in $\Om_M(q_0,q_1)$ (resp. $\La_M$).
 Then
 Proposition~\ref{Tfinito} implies that 
 either $\liminf_n T_n=0$ or $\limsup_n T_n=+\infty$.
 
 If $\limsup_n T_n=+\infty$, 
 Proposition~\ref{Tinf} implies the thesis of the theorem.
 So assume that $\liminf_n T_n=0$.
 
 If $\langle(x_n,T_n)\rangle\subset\Om_M(q_0,q_1)$ with $q_1\ne q_0$
 then Corollary~\ref{Om0} shows that
 $\liminf_n T_n>0$. This contradicts our assumption.
 Hence $q_0=q_1$. If either $\langle(x_n,T_n)\rangle\subset\Om_M(q_0,q_1)$
 with $q_0=q_1$ or $\langle(x_n,T_n)\rangle\subset\La_M$,
 then Proposition~\ref{La0} and Remark~\ref{La0b} imply the thesis
 of the theorem.
 
 \end{proof}

 \subsection{Preliminary lemmas}
 \quad
 \medskip

  \begin{Lemma}\label{Lo}
 If $L$ is convex and quadratic at infinity, then
 \begin{gather*}
 \tfrac 12\,a_0\,|v|_x^2+\theta_x(v)+\psi(x)\le
 L(x,v)\le \tfrac 12\,A_0\,|v|_x^2+\theta_x(v)+\psi(x),
 \\
 -\psi(x)+\tfrac 12\, a_0\,|v|_x^2\le
 E(x,v)\le -\psi(x)+\tfrac 12\, A_0\,|v|_x^2,
 \end{gather*}
 where
  $\theta_x(v):=L_v(x,0)\cdot v$ and
  $\psi(x) :=L(x,0)$.
 \end{Lemma}
 
 \begin{proof}
 Let $L_0(x,v):= L(x,v)-\theta_x(v)-\psi(x)$.
 Let $f(t):=L_0(x,tv)$. Then $f(0)=0$, $f'(0)=0$
 and $f''(t)=v\cdot \frac{\partial L}{\partial v^2}(x,tv)\cdot v$,
 so that $  a_0\;|v|_x^2\le f''(t)\le A_0\;|v|_x^2$.
 Hence
 \begin{align*}
 L_0(x,v) =\int_0^1\int_0^t f''(s)\;ds\,dt
  &\ge\int_0^1\int_0^t a_0\;|v|_x^2\;ds\,dt
 \ge \tfrac 12\,a_0\,|v|_x^2
 \\
 &\le \tfrac 12\,A _0\,|v|_x^2.
 \end{align*}
 
 Now let $g(t):=E(x,tv)= tv\cdot L_v(x,tv)-L(x,tv)$.
 Then $g(0)=-\psi(x)$ and 
 $g'(t)=t\;v\cdot L_{vv}(x,tv)\cdot v$,
 so that $ t\,a_0 |v|_x^2 \le g'(t)\le t\,A_0 |v|_x^2$.
 Therefore, 
 \begin{align*}
 E(x,v)=g(1)= g(0)+\int_0^1 g'(t)\;dt
 &\ge
 -\psi(x)+\tfrac 12\, a_0\, |v|_x^2
 \\
 &\le -\psi(x)+\tfrac 12\, A_0\, |v|_x^2\,.
 \end{align*}
 \end{proof}

 Let $\la>0$ be a Lebesgue number  for our finite atlas
 $\cU=\{U_i\}$ of bounded charts.

 \begin{Lemma}\label{III.8}
 
 Suppose that $L$ is quadratic at infinity. 
 
 If $U_i\in\cU$, $x,y\in U_i$ and $d_M(x,y)<\la$
 then in the chart $U_i$ we have:
 \begin{gather}
  \big[\,L_v (x,v) - L_v (y,w)\,\big]\cdot\zeta 
 +b_3\,\big(|v|+|w|+1\big)\,|\zeta|\,|x-y|
 +A_0\,|\zeta|\,|v-w| \ge 0 .
 \tag{i}\label{lv1}
 \\
 a_0\,|v-w|^2 \leq \big[\,L_v (x,v) - L_v (y,w)\,\big]\cdot(v-w) 
 +b_3\,\big(|v|+|w|+1\big)\,|v-w|\,|x-y|.
 \tag{ii}\label{lv2}
 \end{gather}
 \end{Lemma}
 
 \begin{proof} 
 Recall that the domains $U_i\subset\re^m$ are convex.
 We work in local coordinates as if $L$
 were  defined in $TU_i\subset \re ^{2n}$.
 \begin{multline*}
 \int_{0}^{1} \zeta\cdot L_{vv}\bigl(t\,(x,v) + (1-t)\,(y,w)\bigr)
     \cdot (v-w) \; dt =
  \\ 
 =\bigl(L_v (x,v) - L_v(y,w)\bigr)\cdot\zeta
 \\
 -\int_{0}^{1} 
      \zeta\cdot L_{xv}\bigl(t\,(x,v) + (1-t)\,(y,w)\bigr)
      \cdot (x-y) \; dt.
 \end{multline*}
 This implies~\eqref{lv1}. Using $\zeta=v-w$ one gets~\eqref{lv2}.
 
  \end{proof}


 \begin{Lemma}\label{regu}
 Let $\cC:=C^\infty([0,1],M)\times\re^+$.

 \noindent
 The subsets  $\cC\cap\Om_M(q_0,q_1)$ and  $\cC\cap\La_M$
 are dense in $\Om_M(q_0,q_1)$ and  $\La_M$ respectively.
 \end{Lemma}  

 \begin{proof}
 We prove the lemma for $\Om_M(q_0,q_1)$. The proof for $\La_M$
 is similar. 
  Let $\la>0$ be a Lebesgue number for our finite atlas $\cU$.

  Suppose first that $\length(x)<\la$. Then the image of $x$ lies
  inside of a domain of a chart $U_i\in\cU$
   and by lemma~\ref{equivM} we can assume
  that $M=U_i\subset\re^n$.
 Extend $x:[0,1]\to M$ to $\re$ by setting
 $x(t)=x(0)$ for $t<0$ and $x(t)=x(1)$ for $t>1$.
 Then the extension is also in the Sobolev space 
 $W^{1,2}_{\text{loc}}(\re,\re^m)$, see~\cite[\S 4.1]{EG}.
 Let $\eta\in C^\infty(\re,\re)$ be
 $$
 \eta(t):=\begin{cases} 
          C\,\exp(\tfrac 1{t^2-1}) &\text{if } |t|\le 1,\\
          0 &\text{if } |t|\ge 1,
         \end{cases}
 $$
 where the constant $C$ is chosen such that $\int_\re\eta\; dt=1$.
 For $\e>0$ let 
 $
 \eta_\e(t):=\tfrac 1\e\;\eta(\tfrac t\e).
 $
 Define 
 $$
 x^\e(t):=\int_\re \eta_\e(s-t)\;x(s)\; ds.
 $$
 Then~\cite[\S 4.2.1]{EG}, $x^\e\in C^\infty(\re,\re^m)$, 
 $x^\e\to x$ uniformly on compact subsets 
 and $\dx^\e\to\dx$ in $\cL^2([0,1],\re^m)$.
 Let 
 $$
 y^\e(s):= x^\e(s)+ (1-s)\big(x(0)-x^\e(0)\big)+s\big(x(1)-x^\e(1)\big).
 $$
 Then $(y^\e,T)\in\Om_M(q_0,q_1)\cap\cC$.
 Since $\lim\limits_{\e\to 0}x^\e(0)=x(0)$, 
 $\lim\limits_{\e\to 0}x^\e(1)=x(1)$
 and $\dx^\e\to\dx$ in $\cL^2([0,1],\re^m)$
 we have that
 $\lim\limits_{\e\to 0} (y^\e,T)=(x,T)$ in $\Om_M(q_0,q_1)$.
 
 Now assume that $\length(x)>\la$.  Let $0=s_0<s_1<\cdots<s_N=1$
  be such that $\length(x|_{[s_{i-1},s_{i}]})<\tfrac\la 8$.
  Let $x_i(t)=x\big(s_i+t\,(s_{i+1}-s_i)\big)$, $t\in[0,1]$.
  For each $i$ do the construction above and obtain a $C^\infty$ 
  curve $y_i$ with the same endpoints as $x_i$ and which is 
  near $x_i$ in $\cH^1(M)$.
  The curve $y=y_1*\cdots *y_{N}$, appropriately
  defined in $[0,1]$, is piecewise $C^\infty$ and is near 
  $x$ in $\cH^1(M)$.  
  Let $0=s_0<t_1<s_1<t_2<\cdots<t_N<s_N=1$
  be such that $\length(y|_{[t_j,t_{j+1}]})<\tfrac\la 4$.
  Then $y|_{[t_{j},t_{j+1}]}$ is in the domain of a chart
  $U_j\subset\re^m$ and it is $C^\infty$ in  neighbourhoods
  of $t_j$ and $t_{j+1}$. Let $z_j(s)=y\big(t_j+s\,(t_{j+1}-t_j)\big)$,
  $s\in[0,1]$. Let $c_j:[0,1]\to U_j$ be a $C^\infty$ curve
  such that $c_j=z_j$ in neighbourhoods of $0$ and $1$.
  Extend $(z_j-c_j)$ to $\re$ by setting $(z_j-c_j)(s)=0$ if
  $s\in\re\setminus[0,1]$. Then $(z_j-c_j)$ is $C^\infty$.
  Let $\eta_\e$ be as above and let
  \begin{align*}
  w^\e_j(t):=\int_\re\eta_\e(s-t)\cdot(z_j-c_j)(s)\;ds.
  \end{align*}
  Then $w^\e\in C^\infty(\re,\re^m)$ and $w^\e$ is near
  $(z_j-c_j)$ in $\cH^1(\re^m)$. Since $(z_j-c_j)\equiv 0$
  in neighbourhoods of $0$ and $1$ and $\supp(\eta_\e)\subset[-\e,\e]$,
  if $\e$ is small enough then
   $w^\e=0$ in neighbourhoods of $0$ and $1$. Let
  $$
  z_j^\e(t):=c_j(t)+w^\e(t), \qquad t\in[0,1].
  $$
  Then $z_j^\e$ is $C^\infty$, it is near $z_j$ in 
  $\cH^1(\re^m)$ and coincides with $z_j$ in neighbourhoods of
  $0$ and $1$. Let $z^\e:=y|_{[0,t_1]}*z_1*\cdots*z_{N-1}*y|_{[t_N,1]}$
  appropriately defined on $[0,1]$. Then $z^\e$ is $C^\infty$,
  and $(z^\e, T)$ is near $(x,T)$ in $\Om_M(q_0,q_1)$.

 \end{proof}

\medskip

%
%
%

 \begin{Lemma}\label{Jacobi}
 Suppose that the injectivity radius of $M$ is larger than $2$.
  There is $K>0$ 
 such that if $\ga:[0,1]\to M$ is a geodesic with
 $\lv\dga\rv\le 1$ and $J$ is a Jacobi field along $\ga$
 then
 \begin{align*}
 \max_{s\in[0,1]}\big\{\,|J(s)|,|J'(s)|,\,|J''(s)|\,\big\}\,
 \le K\,\Big[|J(0)|+|J(1)|\Big]
  \end{align*}
 \end{Lemma}
 \begin{proof}\quad
 
  We first prove that there is $K_0>0$ such that
  if $J$ is a Jacobi field along $\ga$
  and $J(0)=0$ then 
  $$
  |J(s)|\le K_0\, |J(1)|.
  $$

   Suppose that $K_0$ does not exist. Then for all $N\in\na^+$
   there is a geodesic $\ga_N:[0,1]\to M$ with $|\dga_N|\le 1$,
   a Jacobi field $J_N$ along $\ga_N$  with $J_N(0)=0$,
   and $s_N\in[0,1]$
   such that $\lv J_N(s_N)\rv> N\,\lv J_N(1)\rv$.
   Since $M$ is compact, taking a subsequence 
   of $\langle s_N\rangle$ we can assume that
   the limits    $s_0=\lim_N s_N\in[0,1]$,
    $(\ov{x},\ov{v})=\lim_N\left(\ga_{N}(0),{\dga_{N}(0)}\right)
          \in TM$ and
    $\lim_N\dfrac{J_N(s_N)}{\lv J_N(s_N)\rv}
          \in TM$
    exist.
   Since Jacobi fields are the projection of the derivative
   of the geodesic flow, which is $C^1$,  the map
   $u\mapsto \dfrac{J_N(u)}{\lv J_N(s_N)\rv}$
   converges to a Jacobi field $I(u)$
   along the geodesic $\de(u)$ with 
   $(\de(0),{\de'}(0))=(\ov{x},\ov{v})$
   such that $I(0)=0$ and  
   \begin{gather*}
   \lv I(1)\rv
   =\lim_N\frac{\lv J_N(1)\rv}{\lv J_N(s_N)\rv}
   \le \lim_N\frac{\lv J_N(1)\rv}{N\lv J_N(1)\rv}=0.
   \end{gather*}
   Also,
   $$
   \lv I(s_0)\rv
   =\lim_N\frac{\lv J_N(s_N)\rv}{\lv J_N(s_N)\rv}=1.
   $$
   Then  $I$ is a non-trivial Jacobi field   
   along a geodesic $\de$ of length $|\ov{v}|\le 1$, which
   is zero at the endpoints. Therefore the
    geodesic $\de$ has conjugate points. This
    contradicts\footnote{When $\ov{v}=0$ and $\de$ is a constant
    geodesic, the Jacobi equation along $\de$ is
    $J''=0$, which has no conjugate points.}
    the hypothesis that the injectivity radius of $M$
    is larger than 2.
    
    Now let $J$ be any  Jacobi field along $\ga$.
   Let $A(s)$, $B(s)$ be the Jacobi fields 
   along $\ga$
   satisfying $A(0)=0$, $A(1)=J(1)$
   and
   $B(0)=J(0)$, $B(1)=0$.
    Since  $\length(\ga)\le 1$
    and the injectivity radius of $M$ is larger than 2,   
    the geodesic $\ga$ has no conjugate points. 
    This implies that such Jacobi fields $A$ and $B$ exist.
       
    By the estimate above   $\lv A(s)\rv\le K_0\;\lv J(1)\rv$.
    Considering the Jacobi field $\tB(s):=B(-s)$ along
    the geodesic $\tga(s):=\ga(-s)$ we get that
    $\lv B(s)\rv\le K_0\;\lv J(0)\rv$.
    Since $J(s)=A(s)+B(s)$, we get that
   \begin{align*}
   |J(s)|\le K_0\,\Big[|J(0)|+|J(1)|\Big]
   \qquad\text{ for all }\quad s\in[0,1].
   \end{align*}

  Since $|\dga|\le 1$,
  from the Jacobi equation
  $J''+R(\dga,J)\,\dga=0$, we get that
  $$
  |J''(s)|\le b\, |J(s)|\le b\,K_0\,\left(|J(0)|+|J(1)|\right),
  $$
  for some $b=b(M)>0$. 
 
   We have that
  $$
  J(1)-T_1\cdot J(0)=\int_0^1 T_s\cdot J'(s)\;ds, 
  $$
  where $T_s:T_{\ga(0)}M\to T_{\ga(s)}M$ is the parallel
  transport along $\ga$. Then there is $s_0\in[0,1]$ such that
  $|J'(s_0)|\le |J(0)|+|J(1)|$. Therefore
  \begin{align*}
  |J'(s)|&\le |J'(s_0)|+\lv\int_{s_0}^s|J''(u)|\;du\rv
  \\
  &\le |J(0)|+|J(1)| + b\, K_0\,\left(|J(0)|+|J(1)|\right),
  \end{align*}
  Now take $K=\max\{K_0,\,1+ b\,K_0\,\}$.

 \end{proof}

 \bigskip

 \subsection{Palais-Smale sequences}\quad

  \bigskip

  During the rest of this section $(x_n,T_n)$ will be a Palais-Smale 
  sequence. This is, $(x_n,T_n)$ will be a sequence in a fixed 
  connected component of  $\La_M$ or $\Om_M(q_0,q_1)$ such that
  $$
  \lv \cA_k(x_n,T_n)\rv<A_1
   \qquad\text{ and }\qquad
   \lV d_{(x_n,T_n)}\cA_k\rV^{\cH^1(M)\times\re^+}_{(x_n,T_n)} 
   < \tfrac 1n.
  $$
  Also, $L$ will be a convex lagrangian on a compact manifold $M$,
  Riemannian at infinity.

  Write $y_n(t):=x_n(t/T_n)$, $0\le t\le T_n$. As in~\eqref{deriv}
  we compute
 \begin{align}
 d_{(x_n,T_n)}\cA_k(\xi,\a) 
          &= \int_0^{1} \big[L_x(x_n,\tfrac{\dx_n}{T_n})\,\xi
            +L_v(x,\tfrac{\dx_n}{T_n})\,\tfrac\dxi{T_n}\,\big]\;T_n\;ds
            +\a\,\int_0^1\big[k-E(x_n,\tfrac{\dx_n}{T_n})\big]\;ds
            \notag  \\
            &=\int_0^{T_n} \big[\,L_x(y_n,\dy_n)\,\zeta
              +L_v(y_n,\dy_n)\,\dzeta\,\big]\; dt
              +\frac{\a}{T_n}\,\int_0^{T_n}\big[k-E(y_n,\dy_n)\big]
              \;dt,   
 \label{derivative}
 \end{align}
 where $\zeta(t):=\xi(t/T_n)$.

  \bigskip

  \begin{Lemma}\label{lt0}
  There exists 
  $B=B(k,A_1,A_2)>0$ such that if 
   $x_n\in \cH^1(M)$, \linebreak 
   $\cA_k(x_n,T_n)\le A_1$ and $T_n\le A_2$, then 
  $$
  \frac{\ell_n^2}{T_n} \le 
  \int_0^{T_n} |\dy_n|^2\;dt
  =\frac 1{T_n}\,\int_0^{1} |\dx_n|^2\;ds
  < B,
  $$
  where $\ell_n:=\length(x_n)$ and $y_n(t)=x_n(t/T_n)$.
   In particular if $T_n\to 0$, then $\lim_n\ell_n =0$.
  \end{Lemma}

  \begin{proof}
  Using \eqref{a1a2},
    \begin{alignat*}{2}
  A_1&\ge \cA_k(x_n,T_n)=A_{L+k}(y_n)
   &&\ge a_1\int_0^{T_n} \lv{\dy_n}\rv^2\;ds - (a_2+|k|)\,T_n, 
   \end{alignat*}
  Let $\ell_n:=\length(y_n)$. By the Cauchy-Schwartz inequality,
  \begin{equation*}
  \ell_n^2 =\bigg(\int_0^{T_n}|\dy_n|\;dt\bigg)^{\!\! 2}
  \le T_n \cdot \int_0^{T_n}|\dy_n|^2\,dt.
  \end{equation*}
 The inequalities above  
 imply the lemma with 
 $B=1+\frac 1{a_1}\big[A_1+( a_2+|k|) A_2\big]$.
 
 \quad
 \end{proof}

  \begin{Corollary}\label{Om0}
    If $(x_n,T_n)\in\Om_M(q_0,q_1)$, $q_0\ne q_1$ and
  $\cA_k(x_n,T_n)<A_1$,
  then $T_n$ is bounded away from zero.
  \end{Corollary}

  \begin{Corollary}\label{La1}
  If $\La_1$ is a connected component of $\La_M$ with a 
  non-trivial free homotopy class, then for all $A_1>0$,
  $$
  \inf\big\{\,T>0\,\big\vert\,(x,T)\in\La_1,\;\cA_k(x,T)<A_1\,\big\}>0.
  $$
   \end{Corollary}
  \begin{proof}\quad
   
    \noindent
    Since $M$ is compact, 
    $\inf\{\,\length(x)\,|\, (x,T)\in\La_1\,\}$
    is positive.
    Now use Lemma~\ref{lt0}.
    
  \end{proof}
  
  \bigskip
  
   \begin{Proposition}\label{La0}\quad
  
  \noindent 
 If a sequence $(x_n,T_n)\in \La_M$ satisfies
 $ \cA_k(x_n,T_n)<A_1$, $\lV d_{(x_n,T_n)}\cA_k\rV<\tfrac 1n$
 and $T_n\to 0$, then there is $q_0\in M$ and a subsequence $x_{n_i}$ such that
              $q_0=\lim_i x_{n_i}(s)$ for all $s\in[0,1]$ and
\renewcommand{\theenumi}{\roman{enumi}}
 \begin{enumerate}
 \item \label{La0i} $\lim_i\cA_k(x_{n_i},T_{n_i})=0$.
 \item\label{La0ii}  $(q_0,0)$ is a singularity of the Euler-Lagrange flow.
 \item\label{La0iii} $E(q_0,0)=k$.   In particular $k\le e_0(L)$.
 \item \label{La0iv}
              $
	      \displaystyle{
              \lim _i\frac 1{T_{n_i}^2}\int_0^1\lv\dx_{n_i}(s)\rv^2\;ds =0.
	      }
              $        
 \end{enumerate}

 In particular, the Dirac probability measure supported on $(q_0,0)$
 is an invariant measure, supported on the energy level $E^{-1}\{k\}$,
  whose $(L+k)$-action is zero and has trivial
 homology.
 
 Also, the (singular) energy level $E=k$ does not
 satisfy the Palais-Smale condition.
 \end{Proposition}

 \begin{Remark}\label{La0b}
 Proposition~\ref{La0} will be applied to sequences in
 $\La_M$ and also to sequences in $\Om_M(q_0,q_1)$ with $q_1=q_0$.
 This means that $\lV d\cA_k\rV$ is to be understood
 as the norm of the derivative $d_{(x_n,T_n)}\cA_k$
 restricted to the subspace $T_{(x_n,T_n)}\Om(x_n(0),x_n(1))
 \subset T_{(x_n,T_n)}\La_M$ given by variational vector fields
 which are zero at the endpoints.
 
 \end{Remark}

 \begin{proof}\quad

   Assume that $1\ge T_n\to 0$.     
   Let $\ell_n:=\length(y_n)$. By Lemma~\ref{lt0}, we have that
   $\lim_n\ell_n=0$. 
  Since $M$ is compact, taking a subsequence, we can assume 
  that $\lim_n y_n(0)=q_0\in M$. 

 Since $\lim_n\ell_n=0$ and $\lim_ny_n(0)=q_0$, we can assume that
  all the curves $y_n$ are in the domain $U_i\subset\re^m$, $m=\dim M$
  of a bounded chart $U_i\in\cU$.
  By lemma~\ref{equivM} we can assume that
  on the  chart $U_i$ we have
  $\lV d_{(x_n,T_n)}\cAk\rV_{\cH^1(\re^m)\times\re^+}<\tfrac 1n$   
  for all $n$.
  From now on we work on the chart $U_i$ as if $M=\re^m$.

  Let $\xi(s):=x_n(s)-x_n(0)$. Then $\xi(0)=\xi(1)=0$.
   Observe that\footnote{Since $\xi(0)=\xi(1)=0$
   this tangent vector is also in $T_{(x_n,T_n)}\Om_M(q,q)$, where
   $q=x_n(0)=x_n(1)$.} 
   $(\xi,0)\in T_{(x_n,T_n)}\La_M$. 
   Let $\zeta(t):=\xi(t/T_n)$, $t\in[0,T_n]$.
  Then $\zeta(0)=\zeta(T_n)=0$, $\dzeta(t)=\dy_n(t)$. 
  Using~\eqref{NormForm}, we have that
  \begin{align*}
   \Big|\; d_{(x_n,T_n)}\cA_k\cdot(\xi,0)\;\Big|
   \le \frac {1}n\, \left[T_n^3\int_0^{T_n}|\dzeta|^2 \;dt\right]^{\frac
   12}
   \le\frac {T_n}n\, \left[\int_0^{T_n}|\dy_n|^2 \;dt\right]^{\frac
   12}.
  \end{align*}

    Using Lemma~\ref{III.8}.\eqref{lv2} with   $w=0$, we get that
  \begin{equation}\label{LTh}
  L_v(x,v)\cdot v
  \ge 
  L_v(y,0)\cdot v - b_3\,|v|\, |x-y| - b_3\,|v|^2\,|x-y|
  +a_0\,|v|^2,
  \end{equation}
  for all $(x,v)\in TU_i$.
  Using inequality~\eqref{LTh} with $(x,v)=(y_n,\dyn)$ and $y=y_n(0)$,
  we get
  \begin{align*}
  d_{(x_n,T_n)}\cA_k&\cdot(\xi,0)
  = \int_0^{T_n}
  \Big[\,L_x(y_n,\dy_n)\cdot\zeta
  +L_v(y_n,\dy_n)\cdot\dzeta\,\Big]\;dt\,,
  \end{align*}
  \begin{alignat*}{2}
  \frac{T_n}n \Big[\int_0^{T_n}|\dyn|^2\;dt\Big]^{\frac 12} 
  &\ge -b_2 &&\int_0^{T_n}\big[1+|\dy_n|^2\big]\,|y_n-y_n(0)|\;dt
  +\theta_{y_n(0)}\cdot\Big(\int_0^{T_n}\dyn\;dt\Big)\;\; +
  \\
  &&&-b_3\int_0^{T_n}|\dyn|\,|y_n-y_n(0)|\;dt
  -b_3\int_0^{T_n}|\dyn|^2\,|y_n-y_n(0)|\;dt\;+
  \\
  &&&+ a_0\,\int_0^{T_n}|\dyn|^2\;dt\,,
   \\
   &\ge -b_2\,&&\ell_n\,T_n + 0 - b_3\,\ell_n^2
  -(b_3+b_2)\,\ell_n \int_0^{T_n}|\dyn|^2\;dt
  +a_0\,\int_0^{T_n}|\dyn|^2\;dt\,,
  \end{alignat*}
  where $\theta_x:=L_v(x,0)$ and $b_2$ is from~\eqref{b2}.
  
  Dividing the last inequality by $T_n$   we have that
  \begin{align}\label{1Tdyn}
  -b_2\,\ell_n  - b_3\,\frac{\ell_n^2}{T_n}
  +\big[a_0-(b_3+b_2)\,\ell_n\big]
  \Big[\frac 1{T_n}\, \int_0^{T_n}|\dyn|^2\;dt\Big]
  &\le 
  \frac 1n \Big[\int_0^{T_n}|\dyn|^2\;dt\Big]^{\frac 12} .
  \end{align}
  Since $\cA_k(x_n,T_n)<A_1$, from lemma~\ref{lt0} 
  in the chart $U_i$ we get
      \begin{equation*}\label{l1t}
  \limsup_n \frac{\ell_n^{\,2}}{T_n}\le  
  \limsup_n \int_0^{T_n}|\dy_n|^2\;dt <+\infty.
  \end{equation*} 
  From~\eqref{1Tdyn} and lemma~\ref{lt0} we get that
  \begin{equation}\label{l2t}
  \limsup_n \frac{\ell_n^{\,2}}{T_n^2}\le  
  \limsup_n \frac 1{T_n}\int_0^{T_n}|\dy_n|^2\;dt <+\infty.
  \end{equation}
  Since $\lim_n T_n=0$,  we get that
  $$
  \lim_n\ell_n=0,\;\; \lim_n\,\frac{\ell_n^2}{T_n}=0,\;\; 
  \lim_n\int_0^{T_n}|\dyn|^2\;dt=0
  \text{ and }\frac1{T_n}\int_0^{T_n}\!\!|\dyn|^2\,dt \text{ is bounded.}
  $$
  Hence, from inequalities~\eqref{l2t}  and~\eqref{1Tdyn}, we get that
  \begin{equation}\label{l2t2}
  \limsup_n \left[\frac 1{T_n}\int_0^{T_n}|\dyn|\;dt\right]^2
  =\limsup_n \frac{\ell_n^2}{T_n^2}
  \le\lim_n \frac 1{T_n}\int_0^{T_n}|\dyn|^2\,dt =0.
  \end{equation}
  Changing variables in the integral, this proves item~\eqref{La0iv}.  

  \eqref{La0i}. By lemma~\ref{Lo}, for all $(x,v)\in TM$
  $$
  \lv L(x,v)+k\rv
  \le\tfrac12\,A_0\,|v|_x^2+\lv\theta_x(v)\rv+\lv\psi(x)\rv+|k|.
  $$
  Then 
  $$
  \lv\cAk(x_n,T_n)\rv \le
  \tfrac 12\,A_0\,\int_0^{T_n}\lv\dyn\rv^2_{y_n}\;dt
  +\ell_n\,\sup_{x\in M}\lV\theta_x\rV
  +T_n\,\left[|k|+\sup_{x\in M}\lv\psi(x)\rv\right].
  $$
  Hence $\lim_n\cAk(x_n,T_n)=0$.
  
  \eqref{La0ii}. Let $h:[0,1]\to [0,2]$ be a smooth function such that
  $h(0)=h(1)=0$ and $\int_0^1 h(s)\;ds =1$. Let
  $\xi(s):= h(s)\, d\psi(x_n(0))\in\re^m$, $s\in[0,1]$ 
  and $\zeta(t):=\xi(t/T_n)$.
  We have that
   \begin{align}
  d_{(x_n,T_n)}\cA_k&\cdot(\xi,0)
  = \int_0^{T_n}\left[
  L_x(y_n,\dyn)\cdot\zeta
  + L_v(y_n,\dy_n)\cdot\dzeta\right]\;dt
  \le \tfrac 1n \lV\xi\rV_{(x_n,T_n)}.
  \label{LxLv}
  \end{align}
    Using~\eqref{b3},
  \begin{align}
  L_x(x,v) 
  &= L_x(x,0)+\int_0^1 \tfrac{d\,}{ds}L_x(x,sv)\;ds
  \notag\\
   L_x(x,v)\cdot\zeta&\ge  d\psi(x)\cdot\zeta
   -b_3\,\big(1+|v|_x\big)\,|v|_x\,|\zeta|.
   \label{Lx}
  \end{align}
  Write $\theta_{q_0}=L_v(q_0,0)$.
  Using~\eqref{Lx} and Lemma~\ref{III.8}.\eqref{lv1} 
  with $(x,v)=(y_n,\dy_n)$,  $(y,w)=(q_0,0)$ 
   in inequality in~\eqref{LxLv}, we get that
   \begin{gather*}
   T_n\int_0^1 d\psi(x_n(s))\cdot h(s)\cdot d\psi(x_n(0))\;ds
   - b_3\, \lV\zeta\rV_\infty\int_0^{T_n}(\lv\dyn\rv+\lv\dyn\rv^2)\;dt
   \\
   +\,\theta_{q_0}\left(\int_0^{T_n} \dzeta\; dt\right)
   -b_3\,\Vert\dzeta\Vert_\infty\,\ell_n\,\int_0^{T_n}
   \left(1+\lv\dyn\rv\right)\;dt
   -A_0\,\Vert\dzeta\Vert_\infty\int_0^{T_n}\lv\dyn\rv\;dt
   \\
   \le \frac{T_n}n \lv d\psi(q_0)\rv\, \Vert \dot{h}\Vert_{\cL^2}.
   \end{gather*}
   In the inequality above the third term is zero. Dividing by $T_n$,
   letting $n\to+\infty$ and using~\eqref{l2t2}, we get 
   $$
   0\le \lv d\psi(q_0)\rv^2=\lim_{n\to+\infty}\int_0^1
   d\psi(x_n(s))\cdot h(s)\cdot d\psi(x_n(0))\;ds
   \le 0.
   $$
   Hence $(q_0,0)$ is a singularity of the Euler-Lagrange flow.
   
   \eqref{La0iii}. We now see that $E(q_0,0)=k$.
    From~\eqref{b1},
  $$
  \lv\psi\big(y_n(t)\big) -\psi\big(y_n(0)\big)\rv
  \le b_1\,\ell_n \qquad \text{for all }\;t\in[0,T_n].   
  $$
  The hypothesis $\lV d_{(x_n,T_n)}\cA_k\rV<\tfrac 1n$ implies
  that
  \begin{align}
  \frac 1n 
  \ge \lv\frac{\partial\cA_k}{\partial T}\Big\vert_{(x_n,T_n)}\rv
  = \frac 1{T_n}\;\lv\int_0^{T_n}
  \Big[\, E(y_n,\dy_n)-k\,\Big]\;dt\;\rv.
  \label{mEn}
  \end{align}
  Using  Lemma~\ref{Lo} we get that
  \begin{align*}
   &-\big[\,\psi\big(y_n(0)\big)+k
  \,\big]-b_1\,\ell_n
  +\frac {a_0}{2\,T_n}\int_0^{T_n}|\dy_n|^2\;dt
  \le \frac 1n\,,
  \\
  -\frac 1n \le&-\big[\,\psi\big(y_n(0)\big)+k
  \,\big]+b_1\,\ell_n
  +\frac {A_0}{2\,T_n}\int_0^{T_n}|\dy_n|^2\;dt.
  \end{align*}
    Then from inequality~\eqref{l2t2}, we get
  $$
  E(q_0,0)=-\psi(q_0)=\lim_n-\psi\big(y_n(0)\big)=k.
  $$

     We have that $E(q_0,0)=k=-\psi(q_0)$ and $d\psi(q_0)=0$, hence
  the point $(q_0,0)\in TM$ is a singularity of the Euler-Lagrange
  flow in the energy level $E=k$. The Dirac measure supported on
  $(q_0,0)$ is an invariant measure whose $(L+k)$-action is zero
  and has trivial homology.

      This (singular) energy level $E=k$ does not satisfy the
  Palais-Smale condition 
  because the curves $(x_n,T_n)$, where  $x_n(t)\equiv q_0$,  $T_n=n$,
  are in the same connected component in $\La_M$
  of closed curves with trivial homotopy class,
  they satisfy $\cA_k(x_n,T_n)\equiv 0$ and $d\cA_k(x_n,T_n)\equiv 0$
  but they do not have an accumulation point in the topology of $\La_M$.
  
  In the case $(x_n,T_n)\in\Om_M(q_0,q_0)$, 
  the same choice $x_n(s)\equiv q_0$,
  $T_n=n$ is an unbounded Palais-Smale sequence in
  $\Om_M(q_0,q_0)$.
 \newline
 \end{proof}

 
 \begin{Corollary}
  If $(q_0,0)\in T_{q_0}M$ is not a singularity of the Euler-Lagrange
 flow and 
 \linebreak
 a sequence $(x_n,T_n)\in\Om_M(q_0,q_0)$ satisfies
  $\cA_k(x_n,T_n)<A_1$ and $\lV d_{(x_n,T_n)}\cA_k\rV<\tfrac 1n$, 
 \linebreak
  then \;$\liminf_n T_n>0$.
 \end{Corollary}

 \begin{Remark}\label{NCt} \quad
 
 Observe that the Hilbert manifolds $\La_M$ and
 $\Om_M(q_0,q_1)$ are not complete with our riemannian
 metric~\eqref{metric} because they do not contain 
 the points $(x,0)\in \cH^1(M)\times\{0\}$ that would be at finite
 distance from $(x,1)$. The discussion above shows that in order
 to prevent a Palais-Smale sequence $(x_n,T_n)$
 from leaving the space at
 $\cH^1(M)\times \{0\}$ we can either
 \begin{itemize}
 \item work on a connected component $\La_1$ of $\La_M$ or
 $\Om_M(q_0,q_1)$ with a non-trivial homotopy class.
 \item work on $\Om_M(q_0,q_1)$ with $q_0\ne q_1$.
 \item work on $\Om_M(q_0,q_0)$ where $(q_0,0)$
       is not a fixed point of the Euler-Lagrange flow.
 \item ask that $E^{-1}\{k\}$ is not a singular energy level.
 \item ask that $\lim_n\cA_k(x_n,T_n)\ne 0$.      
 \end{itemize}
 On a given connected component $\La_1$ of
 $\La_M$ or $\Om_M(q_0,q_1)$ 
 a singular energy level may not satisfy the Palais-Smale
 condition with a counter-example made with sequences
 of curves $(x_n,T_n)$ with $\lim_n T_n=+\infty$
 which spend long time near the singularity. 
 For example in $\Om_M(q_0,q_1)$ when the singularity 
 is hyperbolic and $q_0$ and $q_1$ are respectively
 in the projections of the unstable and stable manifolds 
 of the singularity.
 In such an example
 theorem~\ref{NPS} says that the measure $\mu_n$ defined 
 in page~\pageref{mun} converges to the Dirac measure
 at the singularity. 
 
  \end{Remark}

  \bigskip

  \begin{Proposition}\label{Tfinito}\quad
  
   If a sequence $\{(x_n,T_n)\}_{n\in{\mathbb N}}\subset \Om_M(q_0,q_1)$ 
  or $\{(x_n,T_n)\}_{n\in{\mathbb N}}\subset \La_M$
  satisfies
 $$
  \cA_k(x_n,T_n)<A_1,\;\;\; \lV d_{(x_n,T_n)}\cA_k\rV<\tfrac 1n
 \;\;\text{ and }\;\;0<\liminf_nT_n< +\infty;
 $$ 
 then there exists a convergent subsequence.
 \end{Proposition}

\begin{proof}

 Since $M$ is compact, if $(x_n,T_n)\in \La_M$, taking a subsequence,
 we can assume that $q_0:=\lim_n x_n(0)=\lim_nx_n(1)$ exists.
 In this case write $q_1:=q_0$. Thus in both cases,
 in $\La_M$ and $\Om_M(q_0,q_1)$,
 we have that $q_0=\lim_n x_n(0)$
 and $q_1=\lim_nx_n(1)$.

 Taking a subsequence we can assume that
  $T=\lim_n T_n\in\re^+$ exists. 
 We will extract a Cauchy sequence from $\{(x_n,T_n)\}_{n\in\na}$. 
 Since $T>0$, such Cauchy sequence has a limit in
 $\La_M$ (resp. in $\Om_M(q_0,q_1)$).
 By Lemma~\ref{regu}, there are smooth curves $\tx_n$
 such that $d[(x_n,T_n),(\tx_n,T_n)]<\tfrac 1n$.
 Hence we can assume that the curves $x_n$ are $C^\infty$,
 for if $\{(\tx_n,T_n)\}_{n\in\na}$ is a Cauchy sequence, so is
 $\{(x_n,T_n)\}_{n\in\na}$ and since $\cA_k$ is $C^1$,
 also $\lim_n\lV d\cA_k(\tx_n,T_n)\rV=0$. Similarly, since  
 $d[(x_n,T_n),(x_n,T)]\le |T_n-T|\overset{n}\to 0$, we can assume
 that $T_n=T$ for all $n$. Also, since we are assuming
 that $T_n=T$ is fixed, it is equivalent to use the
 metric~\eqref{metricM} with $f(T)=g(T)=1$.

  Let $y_n(t):=x_n(t/T)$ and
 let $\a_n:[0,1]\to M$ and
 $\be_n:[T+1,T+2]\to M$ be minimal geodesics joining
 $\a_n(0)=q_0$, $\a_n(1)=y_n(0)$; $\be_n(T+1)=y_n(T)$, $\be_n(T+2)=q_1$.
 Taking a subsequence we can assume that $d(x_n(0),q_0)<1$ and 
 $d(x_n(1),q_1)<1$. Then $|\da_n|\le 1$ and $|\dbe_n|\le 1$.
  Define  
 $$
 w_n(t)=
 \begin{cases}
 \a_n(t) &\text{ if}\quad \phantom{T+\;\,}0\le t\le 1,   \\
 y_n(t-1) &\text{ if}\quad \phantom{T+\;\,}1\le t\le T+1, \\
 \be_n(t) &\text{ if}\quad T+1\le t\le T+2.
 \end{cases}
 $$
 Then all the curves $w_n:[0,T+2]\to M$  join $q_0$ to $q_1$.
 Their action is uniformly bounded because 
 \begin{align*}
 A_{L+k}(w_n) &=\cA_k(x_n,T)+A_{L+k}(\a_n)+A_{L+k}(\be_n) \\
            &\le A_1
            +2\cdot\sup_{\lv v\rv\le 1}[L(x,v)+k]=:A_2.
 \end{align*}
 
 By Cauchy-Schwartz inequality and Lemma~\ref{lt0},
 \begin{align*}
 d\big(w_n(t_1),w_n(t_2)\big)
 &\le\int_{t_1}^{t_2}\lv\dwn(s)\rv\;ds
 \le\sqrt{|t_2-t_1|}\left[\int_0^{T}\lv\dwn\rv^2\;ds\right]^{\frac 12}
 \\
 &\le B(k,A_2,T+2)^{\frac 12}\;|t_2-t_1|^{\frac 12}.
 \end{align*}
 Then the family $\{w_n\}$ is equicontinuous.
 By Arzel\'a-Ascoli Theorem there is
 a convergent subsequence of $\{w_n\}$ in the $C^0$ topology. 
 This implies that also $\{x_n\}$ has a
 convergent subsequence in the $C^0$ topology.
 For the sequel we work with a convergent subsequence of $\{x_n\}$.

  We can assume that the injectivity radius of $M$ is larger than
 2. For $n$, $m$ large enough $d(x_n(s),x_m(s))<1$ for all
 $s\in[0,1]$. Let $\ga_s:[0,1]\to M$ be the minimizing geodesic
 joining $\ga_s(0)=x_n(s)$ to $\ga_s(1)=x_m(s)$. Let
 $\Ga:[0,1]\times[0,1]\to M$ be defined by 
 $\Ga(s,r):= \ga_s(r)$. Then 
 \begin{equation}\label{disr}
 \lv\frac{\partial\Ga}{\partial r}(s,r)\rv
 =\lv\dga_s(r)\rv =d_M\big(x_n(s),x_m(s)\big)
 \le\dnm,
 \end{equation}
 where $\dnm=\sup_{s\in[0,1]}d(x_n(s),x_m(s))$.
 Observe that $J(r):=\frac{\partial \Ga}{\partial s}(s,r)$
 is a Jacobi field along $\ga_s$ with $J(0)=\dx_n(s)$ and
 $J(1)=\dx_m(s)$. Since $|\dga_s|\le 1$,
  by Lemma~\ref{Jacobi},
 \begin{align}
 \lv \frac{\partial\Ga}{\partial s}(s,r)\rv
 &=|J(r)|\le K\,\Big[|\dx_n(s)|+|\dx_m(s)|\Big],
 \label{Gr}\\
 \lv \frac{D}{ds}\frac{\partial\Ga}{\partial r}(s,r)\rv
 &= \lv \frac{D}{dr}\frac{\partial\Ga}{\partial s}(s,r)\rv
 =|J'(r)|
 \le K\,\Big[|\dx_n(s)|+|\dx_m(s)|\Big],
 \label{Gr'}
  \end{align}
 
 By Lemma~\ref{lt0},
 \begin{gather}\label{Lnorm}
 \lV \dxn\rV_{\cL^1}\le \lV \dxn\rV_{\cL^2}
 \le T\,B(k,A_1,T) = :B_1,
 \\
 \lV \tfrac{\dxn}T\rV_{\cL^2}
 = \lV \dyn\rV_{\cL^2}
 \le B(k,A_1,T)=:B_2.
 \notag
 \label{Lnormy}
 \end{gather}

 Let $\etanm(s):=\frac{\partial \Ga}{\partial r}(s,1)\in T_{x_m(s)}M$ and
 $\xinm(s):=\frac{\partial \Ga}{\partial r}(s,0)\in T_{x_n(s)}M$. 
  We have  that
  \begin{align*}
  \lV\etanm\rV_{\cH^1(M)}^2
  &=|\etanm(0)|^2+\int_0^1|\deta(s)|^2\;ds
  \\
  &=\lv\frac{\partial\Ga}{\partial r}(0,1)\rv^2
  +\int_0^1\lv\frac{D}{ds}\,\frac{\partial\Ga}{\partial r}(s,1)\rv^2
  \;ds.
  \end{align*}

 From~\eqref{disr}, ~\eqref{Gr'} and ~\eqref{Lnorm}, 
 \begin{align*}
 \lV\etanm\rV_{\cH^1(M)}^2
  &\le d(x_n(0),x_m(0))
  + K\, \big( 2\, \lV x_n\rV_{\cH^1(M)}^2+
  2\,\lV x_m\rV_{\cH^1(M)}^2\big)
  \\
  &\le 1 +4\,K\,B_1 =: K_1.
 \end{align*}
 Similarly,
 $$
 \lV\xinm\rV_{\cH^1(M)}^2\le K_1.
 $$
 Also
 $$
 \lV\etanm\rV_\infty < \dnm
 \qquad\text{ and }\qquad
 \lV\xinm\rV_\infty < \dnm.
 $$
 
 \medskip
 
  Since
 $
 \lim_n \lV\, d_{(x_n , T)} \cA_{k}\rV_{(x_n,T)}= 0,
 $
 for the product norm
 $\lV \,\cdot\,\rV_{\cH^1(M)\times\re}$, with $f(T)\equiv g(T)\equiv 1$,
 we also have that
 \begin{equation}\label{limHN}
 \lim_n \lV\, d_{(x_n , T)} \cA_{k}\rV_{\cH^1(M)\times\re}= 0,
 \end{equation}
 where the derivative is restricted to the tangent space
 $T_{(x_n,T_n)}\Om_M(q_0,q_1)$ [resp. $T_{(x_n,T_n)}\La_M$].
 Therefore
 given $\e>0$ there is  $N>0$ 
 such that 
 \[
 \lV \partial_{x_n}\cA_{k}\vert_{(x_n,T)}\cdot\eta\rV_{\cH^1(M)} 
 <\tfrac 12\,\varepsilon 
 \]
 for every $n\geq N$ and $\lV \eta\rV_{\cH^1(M)} \leq K_1$ with
 $\eta(0)=\eta(1)=0$ when $(x_n,T)\in\Om_M(q_0,q_1)$
 and $\eta(0)=\eta(1)$ when $(x_n,T)\in\La_M$.
 We can take $\eta=\etanm$ and $\eta=\xinm$ defined
 above over $(x_m,T)$ and $(x_n,T)$ respectively. Therefore
  \[
 \lV \partial_{x_m}\cA_{k}\vert_{(x_m,T)}\cdot\etanm 
 -\partial_{x_n}\cA_{k}\vert_{(x_n,T)}\cdot\xinm\rV_{\cH^1(M)}
 <\varepsilon
 \]

  From formula~\eqref{deriv} for $\partial_x\cA_k$, we have that 
  \begin{align}
  \Big|\int_{0}^{1} T\,\Big[\langle\nabla_x L&(x_m,\tfrac{\dx_m}T),
  \,\etanm\rangle
  -  \langle\nabla_x L(x_n, \tfrac{\dx_n}T),
  \,\xinm\rangle\Big] \; ds \;+
 \notag\\
 &+ \int_{0}^{1}\Big[\langle\nabla_v L(x_m,\tfrac{\dx_m}T),
 \, \detanm\rangle- 
 \langle\nabla_v L(x_n, \tfrac{\dx_n}T),
 \, \dxinm\rangle\Big] \; ds\Big| 
 < \varepsilon
 \label{2terms}
  \end{align}
 for $m$, $n$ $>N$. Since $L$ is quadratic at infinity,
 $$
 b_2:=\sup_{(x,v)\in TM}\frac{\lv \nabla_x
 L(x,v)\rv_x}{1+|v|_x^2}<+\infty.
 $$
 The first term in~\eqref{2terms} is bounded by
 $2\,T\, b_2\,(1+B_2^2)\, \dnm$.
  Consequently the second integral in~\eqref{2terms} 
 is small for big $m$, $n$.

 Observe that the integrand in the second term of~\eqref{2terms} is
 \begin{align}
 \langle&\nabla_v L(x_m,\tfrac{\dx_m}T),\, \detanm\rangle
 - 
 \langle\nabla_v L(x_n, \tfrac{\dx_n}T),\, \dxinm\rangle
 \notag\\
 &= \left\langle\nabla_v L\left(\Ga(s,r),\tfrac 1T\,
 \tfrac{\partial\Ga}{\partial s}(s,r)\right),
 \;\tfrac{D}{ds}\,\tfrac{\partial\Ga}{\partial r}(s,r)\right\rangle
 \Big|_{r=0}^{r=1}
 \notag\\
 &=\int_0^1 \frac{D}{dr}\;
 \left\langle\nabla_v L\left(\Ga(s,r),\tfrac 1T\,
 \tfrac{\partial\Ga}{\partial s}(s,r)\right),
 \;\tfrac{D}{ds}\,\tfrac{\partial\Ga}{\partial r}(s,r)\right\rangle
 \; dr
 \notag\\
 &=\int_0^1 
 \Big\langle\partial_x
 \nabla_v L\left(\Ga,\tfrac 1T
 \tfrac{\partial\Ga}{\partial s}\right)
 \cdot \tfrac{\partial\Ga}{\partial r}
 +
 \partial_v\nabla_v L\left(\Ga,\tfrac 1T
 \tfrac{\partial\Ga}{\partial s}\right)
 \cdot \tfrac 1T \tfrac{D}{d r}
 \tfrac{\partial\Ga}{\partial s}\; ,
 \;\tfrac{D}{ds}\tfrac{\partial\Ga}{\partial r}\Big\rangle
 \; dr
 \notag\\
 &=\int_0^1 
 \Big\langle\partial_x
 \nabla_v L\left(\Ga,\tfrac 1T
 \tfrac{\partial\Ga}{\partial s}\right)
 \cdot \tfrac{\partial\Ga}{\partial r}
 \, ,
 \;\tfrac{D}{ds}\tfrac{\partial\Ga}{\partial r}\Big\rangle
 \; dr
 +\int_0^1\Big\langle
 \partial_v\nabla_v L\left(\Ga,\tfrac 1T
 \tfrac{\partial\Ga}{\partial s}\right)
 \cdot \tfrac 1T \tfrac{D}{d r}
 \tfrac{\partial\Ga}{\partial s}\; ,
 \;\tfrac{D}{ds}\tfrac{\partial\Ga}{\partial r}\Big\rangle
 \; dr,
 \label{1ter}
  \end{align}
  here $\partial_x\nabla_vL$ and $\partial_v\nabla_vL$
  are the partial derivatives of the second component
  of $TM\ni(x,v)\mapsto\big(x,\nabla_vL(x,v)\big)\in TM$ with respect 
  to the splitting $T_{(x,v)}TM=H\oplus V$ described in
  page~\pageref{splitting}. The partial derivative
  $\partial_v\nabla_vL(x,v)$ coincides with the second derivative
  of $v\mapsto L(x,v)\in T_xM$ in the vector space $T_xM$.
  
 Since $L$ is quadratic at infinity,
 $$
 b_3:=\sup_{(x,v)\in TM}\frac{\lV \partial_x\,\nabla_v
 L(x,v)\rV}{1+|v|_x}<+\infty.
 $$
 Then, by~\eqref{Gr},  \eqref{disr} and \eqref{Gr'},
 \begin{align*}
 |\text{ first term in~\eqref{1ter} }|
 &\le \int_0^1 b_3 \left[\,1
 +\lv\tfrac 1T \tfrac{\partial\Ga}{\partial s}\rv\,\right]
 \lv \tfrac{\partial\Ga}{\partial r}\rv
 \,\lv\tfrac{D}{ds}\tfrac{\partial\Ga}{\partial r}\rv
 \;dr
 \\
 & \le b_3 \left[\,1+ \tfrac1T\, K \big(|\dx_n(s)|+|\dx_m(s)|\big)\right]
   \dnm \,K\big(|\dx_n(s)|+|\dx_m(s)|\big).
 \end{align*}
 By~\eqref{Lnorm} and Cauchy-Schwartz inequality,
  \begin{align*}
 \int_0^1 |\text{ first term in~\eqref{1ter} }|  \;ds
 &\le b_3 \,( 2\,K B_1+ \tfrac 1T\,4\,K^2 B_1^2)\;\dnm
 \overset{n,m}\longrightarrow 0.
  \end{align*}
 
 Since $\tfrac{D}{d r} \tfrac{\partial\Ga}{\partial s}
 =\tfrac{D}{ds}\tfrac{\partial\Ga}{\partial r}$,
 from~\eqref{a_0}
 we have that
 \begin{align*}
 \int_0^1 |\text{ second term in~\eqref{1ter} }|  \;ds
 &\ge \int_0^1 \int_0^1 a_0\,\frac 1T\,
 \lv\frac{D}{ds}\frac{\partial\Ga}{\partial r}\rv^2\;dr\;ds.
 \end{align*}
 The integral of~\eqref{1ter} corresponds to the 
 second term in the left of~\eqref{2terms}.
 Since the first term in~\eqref{2terms} is small, we get
 that
 \begin{equation}\label{Dd} 
 \lim_{n,m\to+\infty}\int_0^1 \int_0^1 
 \lv\frac{D}{ds}\frac{\partial\Ga}{\partial r}\rv^2\;dr\;ds
 =0.
 \end{equation}

   Using~\eqref{disr}, in $\cH^1(M)\times\re^+$ we have that
 \begin{align*}
 d\big[(x_n,T),\,(x_m,T)\big]^2
 &\le\int_0^1\lV\frac{\partial\,\Ga}{\partial
 r}\rV_{(\Ga(\cdot,r),T)}^2\; dr
 \\
 &=\int_0^1\lv\frac{\partial\Ga}{\partial r}(0,r)\rv^2dr
 +\int_0^1\int_0^1\lv\frac{D}{ds}\frac{\partial\Ga}{\partial
 r}\rv^2\;ds\;dr
 \\
 &\le d_{n,m}^{\;2} 
 +\int_0^1\int_0^1\lv\frac{D}{ds}\frac{\partial\Ga}{\partial
 r}\rv^2\;ds\;dr.
 \end{align*}
 By~\eqref{Dd}, $\{(x_n,T)\}$ is a Cauchy sequence. 
 
 \end{proof}
 
 \bigskip

 \begin{Proposition}\label{Tinf}
 
  Suppose that $L$ is Riemannian at infinity.
   Let $\La_1$ be a connected component of $\Om_M(q_0,q_1)$
  or $\La_M$. 
    If a sequence 
  $\{(x_n,T_n)\}_{n\in{\mathbb N}}\subset \La_1$
  satisfies
 $$
 \lv \cA_k(x_n,T_n)\rv<A_1,\;\;\; \lV d_{(x_n,T_n)}\cA_k\rV<\tfrac 1n
 \;\;\text{ and }\;\;\limsup_nT_n= +\infty,
 $$ 
  then there exists a Borel probability measure $\mu$, invariant
 under the Euler-Lagrange flow, supported on a connected component
 of the energy level $E\equiv k$, which has homology $\rho(\mu)=0$
 and whose $(L+k)$-action is zero.
 \end{Proposition}

 In the proof of this proposition we can not
 use Lemma~\ref{equivM}.\eqref{equivCH} on the equivalence of
 the metric of $\cH^1(M)\times\re^+$ to the metric
 of $\cH^1(\re^m)\times\re^+$ on local charts 
 because the times $T_n$ are not bounded.
 Here we shall use strongly that the lagrangian
 is  Riemannian at infinity and not only quadratic
 at infinity.

 \label{splitting}
 We first fix the notation used in 
 the proof of Proposition~\ref{Tinf}. 
 There is a canonical splitting of the tangent space
 $$T_{\theta}TM=H(\theta)\oplus V(\theta)$$
 where the {\it vertical subspace} $V(\theta)$ is the kernel of
 the derivative $d_\theta\pi$ of the projection $\pi:TM\to M$ and 
 the {\it horizontal subspace} $H(\theta)$ is the kernel of the connection 
 map $K:T_\theta TM\to T_{\theta}M$.
 Both subspaces are naturally  identified with 
 $T_{\pi(\theta)}M\approx H(\theta)\approx V(\theta)$
 in the following way: a tangent vector $\zeta=(h,v)\in H\oplus V$
 has horizontal and vertical components given by
 $h=d_\theta\pi(\zeta)\in T_\theta M\approx H(\theta)$
 and $v=K(\zeta)\in T_\theta M\approx V(\theta)$.
 
 The Sasaki metric on $TM$ is given by
 \begin{align*}
 \langle \zeta_1,\zeta_2\rangle_\theta
 :&=\langle d_\theta\pi(\zeta_1), d_\theta\pi(\zeta_1)\rangle_{\pi(\theta)}
 +\langle K(\zeta_1),K(\zeta_2)\rangle_{\pi(\theta)}
 \\
 &=\langle h_1,h_2\rangle_{\pi(\theta)}
  +\langle v_1,v_2\rangle_{\pi(\theta)},
 \end{align*}
 where $\zeta_i=(h_i,v_i)\in H(\theta)\oplus V(\theta)$, $i=1,2$.

 The identification $TM\longleftrightarrow T^*M$
 induced by the Riemannian metric 
 $(x,v)\leftrightarrow \langle v,\cdot\,\rangle_x=p$
 preserves the norm on each fiber
 and the canonical symplectic form $\om=dp\wedge dx$ on $T^*M$
 is sent to the form
 \begin{align*}
 \om_{\theta}(\zeta_1,\zeta_2)
 &=\langle K(\zeta_1),d\pi(\zeta_2)\rangle_{x}
 -\langle K(\zeta_2),d\pi(\zeta_1)\rangle_{x}
 \\
 &=\langle v_1,h_2\rangle_{\pi(\theta)}
 -\langle v_2,h_1\rangle_{\pi(\theta)}.
 \end{align*} 
 We shall ambiguously use this identification 
 along the rest of this section.

 Let $H:TM\to\re$  be the hamiltonian associated to $L$:
 \begin{align*}
 H(x,p) :&= \max_{v\in T_xM} \langle p,v\rangle_x -L(x,v), 
 \qquad p\in T_xM.
 \end{align*}
  The hamiltonian vector field $X$  on $TM$ is given by
 $i_X\om =- dH$. Its flow lines solve the hamiltonian equations
  \begin{equation*}
  \dx = \nabla_pH(x,p)\,, 
   \qquad \tfrac{D}{dt}\, p = -\nabla_x H(x,p). 
  \end{equation*}
 When seen in $TM$, $\nabla_x H$ and $\nabla_p H$ 
 are the projections in the horizontal and vertical subspaces
 of the gradient of the hamiltonian  $H$ with respect to 
 the Sasaki metric. The hamiltonian flow of $H$ is conjugated to 
 the Euler-Lagrange flow of $L$ by the {\it Legendre transform}
 $\cL(x,v)=(x,\nabla_v L(x,v))=(x,p)$.
 This is, 
 $\langle p,\cdot\,\rangle_x=\tfrac{\partial}{\partial v}L(x,v)$,
 in the vector space $T_xM$.

 Observe that
 $$
 \nabla_xL(x,v) =-\nabla_xH(x,p),
 \qquad\text{ if }\quad p=\nabla_v L(x,v).
 $$
 If $(x,T)\in\cH^1(M)\times\re^+$, $x\in C^\infty([0,1],M)$
  and $\xi\in T_{(x,T)}\Om_M(x(0),x(1))$,
 the partial derivative of the action functional is given by
 \begin{equation}\label{dAkp}
 d_{(x,T)}\cA_k(\xi,0)
 =\int_0^{T}\langle p,\tfrac{D}{dt}\zeta\rangle_{y}
    -\langle\nabla_xH(y,p),\zeta\rangle_{y} \;\;dt,
 \end{equation}
 where $y(sT):=x(s)$, $\zeta(sT):=\xi(s)$
 and $p(t)=\nabla_vL(y(t),\dy(t))$.
 Let $\cT(t,s):T_{y(s)}M\to T_{y(t)}M$ be the parallel 
 transport along $y(t)$.
 For $t\in[0,T]$ let
 $$
 \Hx(t):=\cT(t,0)\cdot a+\int_0^t \cT(t,s)\cdot\nabla_xH(y(s),p(s))\;ds,
 $$ 
 where the constant $a\in T_{y(0)}M$ is chosen such that
 \begin{equation}\label{defa}
  \int_0^{T} \cT(T,t)\cdot\big[\,p(t)+\Hx(t)\big]\;dt = 0.
 \end{equation}
 Let 
 \begin{equation}\label{defrho}
 \rho(t):=p(t)+\Hx(t).
 \end{equation}
 Integrating by parts in~\eqref{dAkp} we have that
 \begin{equation}\label{dAkr}
 d_{(x,T)}\cA_k(\xi,0)
 =- \big\langle \Hx,\zeta\big\rangle_{y}\Big\vert_0^T 
 +\int_0^T\big\langle p(t)+\Hx(t),
 \tfrac{D}{dt}\zeta(t)\big\rangle_{y(s)}\;dt. 
 \end{equation}

 Define $\zeta_1(t)$ by
 $$
 \zeta_1(t):=\int_0^t \cT(t,s)\cdot\big[\,p(s)+\Hx(s)\big]\;ds.
 $$ 
 By~\eqref{defa}, $\zeta_1(0)=\zeta_1(T)=0$.
 Then if $\xi_1(s):=\zeta_1(sT)$ we have that
 $(\xi_1,0)\in T_{(x,T)}\Om_M(x(0),x(1))$ and also
 $(\xi_1,0)\in T_{(x,T)}\La_M$ if $x(0)=x(1)$.
 Observe that 
 $$
 \tfrac{D}{dt}\,\zeta_1(t)=\rho(t)=p(t)+\Hx(t).
 $$
 Applying~\eqref{dAkr}, we get 
 \begin{align*}
 d_{(x,T)}\cA_k(\xi_1,0)
 =\int_0^T\lv\rho(t)\rv_{y(t)}^2\;dt
 &\le \lV d_{(x,T)}\cA_k\rV\cdot\lV\xi_1\rV_{(x,T)}
 \\
 &\le  \lV d_{(x,T)}\cA_k\rV\cdot
  e^{-2T^2}
 \left[\int_0^T\lv\tfrac{D}{dt}\zeta_1\rv_{y(t)}^2\right]^{\frac 12}
 &\text{ if } T\ge 10,
 \\
 &\le  \lV d_{(x,T)}\cA_k\rV\cdot
  e^{-2T^2}
 \left[\int_0^T\lv\rho(t)\rv_{y(t)}^2\right]^{\frac 12}
&\text{ if } T\ge 10.
 \end{align*}
 Therefore 
 \footnote{ Inequality~\eqref{rhoL2} 
 is our fundamental estimate for the rest of the
 section. Observe that since, by definition, 
 $p(t)=\nabla_vL(y,\dy)$,
 the first hamiltonian equation 
 \quad $\dy=\nabla_pH(y,p)$ \quad 
 follows from the Legendre transform
 of $H$. Thus the quantity $\rho(t)=p(t)+\H_x(t)$
 measures the deviation of $(y(t),p(t))$ from being
 a solution of the second hamiltonian equation \, $\dot p=-\nabla_x H(x,p)$.
 Our problem now is to obtain a true invariant measure from
 this $\cL^2$ estimate.}
 \begin{equation}\label{rhoL2}
 \lV\,\rho\,\rV_{\cL^2([0,T])}\le e^{-2T^2}
 \lV d_{(x,T)}\cA_k\rV\,,
 \qquad\text{ if } T\ge 10.
 \end{equation}
 
 \bigskip

 \begin{proof}
 [\bf Proof of Proposition~\ref{Tinf}:]\quad
 
 We can assume that $T_n\to+\infty$. 
 Moreover, we can assume that
 $$
 T_n\ge n\ge 10.
 $$
  Since $\cA_k$ and $d\cA_k$ are continuous on $\Om_M(q_0,q_1)$
  and $\La_M$, by Lemma~\ref{regu}
 we can assume that $x_n:[0,1]\to M$ is
 $C^\infty$ for all $n$.
 Observe that if $q_0=q_1$,
 $T_{(x,T)}\Om_M(q_0,q_1)\subset T_{(x,T)}\La_M$.
 In the sequel we shall only use tangent vectors in
 $T_{(x,T)}\Om_M(q_0,q_1)$, so that the arguments apply for both
 $\Om_M(q_0,q_1)$ and $\La_M$.

 Let $y_n(t):=x_n(t/T_n)$, $t\in[0,T_n]$ and 
 $$
 p_n(t):=\nabla_vL(\yn(t),\dyn(t))\in T_\yn M. 
 $$
 Let $\rho_n(t)$ be defined as in~\eqref{defrho}
 for the path $(y_n(t),p_n(t))$:
 \begin{equation*}
 \rho_n(t):=p_n(t)+\Hx(t).
 \end{equation*}
 In particular
 \begin{equation}\label{drho}
 \tfrac{D}{dt}\,\rho_n
 =\tfrac{D}{dt}\,p_n+\nabla_xH(y_n,p_n).
 \end{equation}
 From~\eqref{rhoL2} we have that
 \begin{equation}\label{rhon}
 \lV\,\rho_n\rV_{\cL^2([0,T_n])}
 =\Big[\int_0^{T_n}\lv\rho_n(t)\rv_{y_n(t)}^2\;dt\Big]^{\frac12}
 \le
 \frac{e^{-2 T_n^2}}{n}.
 \end{equation}

  Let 
  $$
  A_n:=\big\{\,t\in[0,T_n]\,\big\vert\,
  \lv\rho_n(t)\rv_{y_n(t)}<e^{-\frac 32 T_n^2}\,\big\}.
  $$ 
  If $m$ is the Lebesgue measure on $[0,T_n]$, we have that
  \begin{align*}
  e^{-3T_n^2}\,m(A_n^c)\le \int_0^{T_n} \lv\rho_n(t)\rv^2\;dt
  < \tfrac 1{n^2}\; { e^{-4T_n^2}}.
  \end{align*}
  Thus 
  $ m(A_n^c)< \tfrac 1{n^2} \, e^{-T_n^2}$ and hence
  \begin{equation}\label{mAn}
  m(A_n)> T_n -\tfrac 1{n^2} \, e^{-T_n^2}.
  \end{equation}

  We assume that the Riemannian metric on $M$
  has injectivity radius larger than 2.
  Since the lagrangian is Riemannian at infinity then the hamiltonian
  is also Riemannian at infinity: $H(x,p)=\tfrac 12\;\lv p\rv_x^2$
  when $|p|_x\ge R$.  Define $R_1\ge R$ by
  $$
  \tfrac 12\,R_1^2 = \sup\,\big\{\,H(x,p)\,,\,
  \tfrac 12\,|\nabla_pH(x,p)|_x^2
  \;\;\big |\;\; |p|_x\le R\,\big\}
  \ge \tfrac 12\,R^2.
  $$

   Let $d_1>|k|$ be such that 
   \begin{equation}\label{a1}
   H(x,p)>\tfrac 12\,\lv p\rv_x^2- d_1,\qquad\text{for all } (x,p)\in TM.
   \end{equation}
   Choose $R_0>0$ such that
   \begin{equation}\label{R0}
   R_0\gg \max\big\{\,R_1\,,\;10\big(|k|+d_1+1\big)\,\big\}.
   \end{equation}
   In particular
   $$
   H(x,p)=\tfrac 12 \, \lv p\rv_x^2 \quad\text{ if}\quad \lv p\rv_x>R_0.
   $$
   Also
   $$
   \nabla_pH(x,p)=p\quad \text{ and }\quad\nabla_xH(x,p)=0
   \quad\text{ if }\; |p|_x>R_0.
   $$

   \begin{Lemma}\label{R1}
   
   \begin{equation*}
   |p|_x<R_0\; \Longleftrightarrow\;
   H(x,p)<\tfrac 12\,R_0^2
   \; \Longleftrightarrow\;
   |v|_x<R_0,\text{ where } v=\nabla_pH(x,p).
   \end{equation*}
   \end{Lemma}
   
   \begin{proof}
   
   If $H(x,p)\ge\tfrac 12\,R_0^2>\tfrac 12\,R_1^2$ then
   $|p|_x>R$, $v=\nabla_pH(x,p)=p$ and $H(x,p)=\tfrac 12\,|p|_x^2$.
   Hence $|v|_x=|p|_x\ge R_0$. 
   If $|p|_x\ge R_0>R$ then 
   $H(x,p)=\tfrac 12\,|p|_x^2\ge \tfrac 12\,R_0^2$.
   If $|v|_x=|\nabla_pH(x,p)|_x\ge R_0>R_1$
   then $|p|_x>R$, $v=\nabla_pH(x,p)=p$
    and $H(x,p)=\tfrac 12\,|p|_x^2=
    \tfrac 12\,|v|_x^2\ge \tfrac 12\,R_0^2$.
   \end{proof}
   
   We start by estimating the difference between $(y_n,p_n)$ and an
   orbit of the hamiltonian flow. 
 
 \begin{Lemma}\label{ApproxSol}\quad
 
 \noindent
 Given $t_0\in[0,T_n]$, let 
 $\big(x(t),q(t)\big)$ be the solution of the hamiltonian equations
 $$
 \dx=\nabla_pH(x,q)\;,\qquad\quad \dot{q}=-\nabla_xH(x,q),
 $$
 with initial conditions $x(t_0)=y_n(t_0)$, $q(t_0)=p_n(t_0)$. 
%
%

  There is $n_0>0$ such that if
  $n>n_0$, $t_0\in A_n$ and $\lv p_n(t_0)\rv\le R_0$, then 
  for all $t\in[0,T_n]$,
  \begin{gather}
   d_M(x(t),y_n(t))<1,
   \label{bdsol2}
   \notag\\
   d_{TM}\big[(x(t),q(t)),(y_n(t),p_n(t))\big]
   \le \lv \rho_n(t)\rv_{y_n(t)}
          + e^{-T_n^2}.
  \label{bdsol3}
  \end{gather}
 \end{Lemma}

 \begin{Remark}\label{RApproxSol}
 The bounds in Lemma~\ref{ApproxSol} are actually made for
 $z_t$ defined in~\eqref{zt} instead of the
 distance between $(x(t),q(t))$ and $(y_n(t),p_n(t))$.
 \end{Remark} 
 \medskip 

 \begin{proof}
 We only prove the estimates for $t>t_0$. The case $t<t_0$ is similar.
 
 Recall that using Lemma~\ref{regu},
 we are assuming that $y_n$ is $C^\infty$.
 Let $\ga_t:[0,1]\to M$ be a geodesic joining $x(t)$ to
 $y_n(t)$ such that $\ga_{t_0}(s)\equiv y_n(t_0)$
  for all $s\in[0,1]$ and that $f(s,t):=\ga_t(s)$ is
  $C^\infty$. Let
 \begin{align*}
  e_t:&=\lv\dga_t\rv=\length(\ga_t)\ge d_M(x(t),y_n(t)).
  \end{align*}
   Let $I_n$ be the maximal interval in $[0,T_n]$ containing 
   $t_0$ such that $|e_t|<1$ for all $t\in I_n$.
  
   We first prove that  there is $B=B(L,R_0)>0$ such that  
   for all $n$, if $|p_n(t_0)|\le R_0$ and $t\in I_n$ then
  \begin{align}
  d_{TM}\big[(x(t),q(t)),(y_n(t),p_n(t))\big]
  \le
  \lv \rho_n(t)\rv_{y_n(t)}
          +\big(1&+e^{B(t-t_0)}\big)\,\lv \rho_n(t_0)\rv_{y_n(t_0)}
        \notag  \\
          &\quad+B\,e^{B(t-t_0)}\,
          \lV\rho_n\rV_{\cL^2([0,T_n])}.
          \label{bdsol}
  \end{align}

  Let $A_t(s_2,s_1):T_{\ga_t(s_1)}M\to T_{\ga_t(s_2)}M$ be the
  parallel transport along $\ga_t$.
  Let
  \begin{equation}\label{zt}
  z_t:=\length(\ga_t)
  +\lv p_n(t)-A_t(1,0)\cdot q(t)\rv_{y_n(t)}.
  \end{equation}
  Let $\La_t:[0,2]\to TM$ be the curve defined by
  \begin{equation}\label{CurLa}
  \La_t(s)=\begin{cases}
  \big(\ga_t(s),A_t(s,0)\cdot q(t)\big) &\text{ if }0\le s\le 1,
  \\
  \big(y_n(t),(s-1)\,p_n(t)+(2-s)\,A_t(1,0)\cdot q(t)\big)
  &\text{ if }1\le s\le 2.
   \end{cases}
  \end{equation}
  Then 
  \begin{align*}
  d_{TM}\big((x(t),q(t))&,(y_n(t),p_n(t))\big)
  \le \length(\La_t) \\
  &\le\int_0^1 \lv\dga_t\rv \;ds+ \int_1^2 \lv p_n(t)-A_t(1,0)\cdot
  q(t)\rv_{y_n(t)}\;ds \\
  &\le z_t.
  \end{align*}
  Recall that $f(s,t)=\ga_t(s)$. Then
  \begin{align*}
  e_t^2=\int_0^1 e_t^2\;ds=\int_0^1\lv\dga_t\rv^2\;ds
  =\int_0^1\left\langle\frac{\partial f}{\partial s},
  \frac{\partial f}{\partial s}\right\rangle_{f(s,t)}\;ds.
  \end{align*}
  Since $\frac{D}{ds}\frac{\partial f}{\partial
  s}=\frac{D}{ds}\dga_t=0$, we have that
  \begin{align}
  \frac 12\;
  \frac{d \,e_t^2}{dt}  &=
  \frac 12\int_0^1 \frac{D}{dt}
  \left\langle\frac{\partial f}{\partial s},
  \frac{\partial f}{\partial s}\right\rangle_{f(s,t)} ds
  =\int_0^1 \frac{D}{ds}
  \left\langle\frac{\partial f}{\partial t},
  \frac{\partial f}{\partial s}\right\rangle_{f(s,t)} ds
  \notag\\
  &=\left\langle\frac{\partial f}{\partial t}(1,t),
  \frac{\partial f}{\partial s}(1,t)\right\rangle_{y_n(t)}
  -
  \left\langle\frac{\partial f}{\partial t}(0,t),
  \frac{\partial f}{\partial s}(0,t)\right\rangle_{x(t)}
  \notag\\
  &=\big\langle\dy_n(t),\dga_t(1)\big\rangle_{y_n(t)}
  -\big\langle A_t(1,0)\cdot \dx(t),\, A_t(1,0)\cdot\dga_t(0)
  \big\rangle_{y_n(t)}
  \notag\\
  &= \big\langle\dy_n(t)-A_t(1,0)\cdot \dx(t),
  \dga_t(1)\big\rangle_{y_n(t)}
  \notag\\
  &\le \big\vert\dy_n(t)-A_t(1,0)\cdot \dx(t)\big\vert_{y_n(t)}
  \; e_t.
  \label{de2t}
  \end{align}
  Since $p_n=\nabla_vL(y_n,\dy_n)$, 
  $\dyn=\nabla_pH(y_n,p_n)$.
 Since $\nabla_pH(y,p)=p$ when
  $|p|_y>R_0$,  $(x,p)\mapsto \nabla_pH(x,p)$ has bounded derivative
  on $TM$.
  Let $K_1>1$ be a bound for its derivative. Then
  \begin{equation}\label{difA}
  \begin{aligned}
  \lv\dy_n-A_t(1,0)\cdot\dx\rv
  &=\lv\nabla_pH(y_n,p_n)-A_t(1,0)\cdot \nabla_pH(x,q)\rv
  \\
  &\le\int_0^1\lv\tfrac{D}{ds}\left[A_t(1,\min\{2s,1\})\cdot
  \nabla_p H\big(\La_t(2s)\big)\right] \rv\;ds \\
  &=\int_0^1\lv\tfrac{D}{ds}\,
  \nabla_p H\big(\La_t(2s)\big) \rv\;ds \\
  &\le K_1\; \length(\La_t)
  = K_1\; z_t. 
  \end{aligned}
  \end{equation}
  Thus, from~\eqref{de2t},
   \begin{align*}
  \frac 12\;\frac{d\,e_t^2}{dt}
  &\le K_1\; z_t\; e_t.
  \end{align*}
  Since $e_{t_0}=0$, 
  \begin{equation}\label{etb}
  e_t = \int_{t_0}^t \frac 12\,\frac 1{e_t}\,\frac{d\,e_t^2}{dt}\;dt
  \le \int_{t_0}^t K_1 \, z_t\;dt.
  \end{equation}
   
  Let $T(t_2,t_1):T_{y_n(t_1)}M\to T_{y_n(t_2)}M$ 
  be the parallel transport along $y_n(t)$.
   Since $p_n(t_0)=q(t_0)$, we have that
  \begin{alignat}{3}
  p_n(\tau)&-A_\tau(1,0)\,\cdot  q&&(\tau)
  =&&\!\!\!\int_{t_0}^\tau T(\tau,t)\left[\tfrac{D}{dt}p_n(t)-
    \tfrac{D}{dt}A_t(1,0)\cdot q(t)\right]\;dt,
  \notag\\
   &=\int_{t_0}^\tau T(\tau,t)&&\Big[\tfrac{D}{dt}p_n
   \,+&&\,\nabla_xH(y_n,p_n)\Big]dt
   \notag \\ 
   &&&+\int_{t_0}^\tau &&\!\!\!T(\tau,t)\Big[A_t(1,0)\cdot\nabla_xH(x,q)
   -\nabla_xH(y_n,p_n)\Big]dt
   \label{3sum}\\
   &&&&&·+\int_{t_0}^\tau T(\tau,t)\Big[-A_t(1,0)\cdot\nabla_xH(x,q)
   -\tfrac{D}{dt}A_t(1,0)\cdot q(t)\Big]dt.
   \notag
   \end{alignat}
  Since $\nabla_xH(x,p)=0$ if $|p|_x>R_0$, the function $\nabla_xH$ 
  has bounded derivative on $TM$. Then, as in~\eqref{difA}, 
  $$
  \text{the norm of the second term in~\eqref{3sum} }
  \le \int_{t_0}^\tau K_2\;z_t\,dt,
  $$
  where $K_2$ is a bound for the derivative of $\nabla_xH$.
  We estimate the third term.
  Since $(x,q)$ is a solution of the hamiltonian equations,
  then
  $$
  \tfrac{D}{dt}\;q=-\nabla_xH(x,q).
  $$
  Let
  $F(s,t):=A_t(s,0)\cdot \tfrac{D}{dt}\,q(t)
  -\tfrac{D}{dt}\,A_t(s,0)\cdot q(t)\in T_{f(s,t)}M$, then 
  $F(0,t)\equiv 0$ and 
  \begin{align*}
  \tfrac{D}{ds} F(s,t)
  &=
  \tfrac{D}{ds}\Big[A_t(s,0)\cdot \tfrac{D}{dt}\,q(t)
  -\tfrac{D}{dt}\,A_t(s,0)\cdot q(t)\Big]
  \\
  &=0-\tfrac{D}{ds}\,\tfrac{D}{dt}\,A_t(s,0)\cdot q(t)
 \\
 &= -\tfrac{D}{dt}\,\tfrac{D}{ds}\,\big[A_t(s,0)\cdot q(t)\big]
 +R\Big(\tfrac{\partial f}{\partial s},\tfrac{\partial f}{\partial t}\Big)
 \big[A_t(s,0)\cdot q(t)\big]
 \\
 &=R\Big( \dga_t(s),\tfrac{\partial f}{\partial t}\Big)
 \big[A_t(s,0)\cdot q(t)\big],
  \end{align*}
 where $R$ is the curvature tensor.
 Let $K_3>1$ be such that $\lv R(u,v)w \rv_x\le K_3\,|u|_x \,|v|_x\,|w|_x$
 for all $x\in M$, $u,v,w\in T_xM$.

 Observe that $J(s)=\frac{\partial f}{\partial t}(s,t)$ is a
 Jacobi field along the geodesic $\ga_t$  with 
 $J(0)=\dx(t)$, $J(1)=\dy_n(t)$ and that if $t\in I_n$ then
 $|\dga_t|<1$. By Lemma~\ref{Jacobi} there is $K_4>1$ such that
 $$
 \lv \frac{\partial f}{\partial t}\rv_{f(s,t)}
 \le K_4\,\Big[\lv\dx(t)\rv+\lv\dy_n(t)\rv\Big]
 \qquad\text{ for all }\quad
 (s,t)\in[0,1]\times I_n.
 $$


 Then, using~\eqref{difA},
 \begin{align*}
 \lv\tfrac{D}{ds}F(s,t)\rv
 &\le K_3\, e_t\,  K_4\,\Big[\lv\dx(t)\rv+\lv\dy_n(t)\rv\Big]
 \,\lv q(t)\rv
 \\
 &\le  K_3\, e_t\,  K_4\,\Big[\lv\dx(t)\rv
 +\lv\dy_n(t)-A_t(1,0)\cdot \dx(t)\rv
  +\lv\dx(t)\rv\Big]
 \,\lv q(t)\rv
 \\
 &\le  K_3\,  K_4\, e_t\,\Big[2\,\lv\dx(t)\rv
 +K_1\,z_t  \Big] \,\lv q(t)\rv.
 \end{align*}

 Since, by hypothesis, $\lv q(0)\rv\le R_0$,
 by Lemma~\ref{R1}, 
 $E(x,\dx)=H(x,q)\le \tfrac 12 \,R_0^2$
 and hence, by Lemma~\ref{R1},
  $\lv\dx(t)\rv\le R_0$ and $\lv q(t)\rv\le R_0$ for all $t$.
  Let $K_5:=2\, K_1 K_3 K_4$. If $t\in I_n$ then $|e_t|\le 1$ and
  hence
 \begin{align*}
 \lv F(1,t)\rv_{y_n(t)}  
      &= \lv \,0 +\int_0^1 A_t(1,s)\cdot \tfrac{D}{ds}F(s,t)\; ds\,
      \rv_{y_n(t)}
       \\
       &\le \int_{0}^1\big[ K_5\, R_0^2\, e_t + K_5\, R_0\, z_t\big]\;ds
        \\
       &\le 2\, K_5\, R_0^2 \, z_t
       \qquad\qquad \text{ for all }t\in I_n.
 \end{align*}
 Thus, when $\tau\in I_n$,
 \begin{align*}
 \lv\text{third term in~\eqref{3sum}}\rv
 \le \int_{t_0}^\tau  \lv F(1,t)\rv_{y_n(t)}\;dt
 \le 2\, K_5\, R_0^2\int_{t_0}^\tau z_t \, dt.
 \end{align*}
 From~\eqref{zt}, \eqref{etb}, \eqref{3sum} and~\eqref{drho}
  when $\tau\in I_n$, we get that
 $$
 z_\tau
  \le K_1\int_{t_0}^\tau z_t \;dt
    +\lv\int_{t_0}^\tau \left[T(\tau,t)\cdot 
     \tfrac{D}{dt}\rho_n(t)\right]\, dt\rv_{y_n(\tau)}
    + K_2\int_{t_0}^\tau z_t\; dt
    + 2\,K_5\,R_0^2\int_{t_0}^\tau z_t\; dt
 $$ 
 \begin{align}
  z_\tau
   &\le \lv \rho_n(\tau) -T(\tau,t_0)\cdot\rho_n(t_0)\rv_{y_n(\tau)} 
    +B \int_{t_0}^\tau z_t\; dt
  \notag\\
 z_\tau &\le  \lv \rho_n(\tau)\rv_{y_n(\tau)}
          +\lv \rho_n(t_0)\rv_{y_n(t_0)}
          + B\int_{t_0}^\tau z_t\;dt
	  \qquad\text{ when }\tau\in I_n.
 \label{utau}
 \end{align}
 where $B:=\max\{1, K_1+K_2+2\,K_5\,R_0^2\}$.
 Let $u(\tau):=\int_{t_0}^\tau z_t\;dt$.
 Then, using~\eqref{utau}, we have that
 \begin{align*}
 \tfrac{d}{dt}\big( e^{-B(t-t_0)}\, u(t)\big)
 &=e^{-B(t-t_0)}\big( z_t - B \,u(t)\big)
 \\
 &\le e^{-B(t-t_0)}\big(
   \lv \rho_n(t)\rv_{y_n(t)}
          +\lv \rho_n(t_0)\rv_{y_n(t_0)}\big)
	  \qquad\text{ for } t\in I_n. 
 \end{align*}
 Since $u(t_0)=0$, 
 \begin{align*}
 u(\tau) &\le
  e^{B(\tau-t_0)}\int_{t_0}^\tau e^{-B(t-t_0)} \lv \rho_n(t)\rv_{y_n(t)}\;dt
  +
   e^{B(\tau-t_0)}\int_{t_0}^\tau e^{-B(t-t_0)} 
   \lv \rho_n(t_0)\rv_{y_n(t_0)}\;dt
 \\
 &\le \frac{e^{B(\tau-t_0)}}{\sqrt{2B}}\;\lV\rho_n\rV_{\cL^2([0,T_n])}
 + \frac{e^{B(\tau-t_0)}}{B}\;{\lv \rho_n(t_0)\rv_{y_n(t_0)}}
 \qquad\text{ when }\tau\in I_n.
 \end{align*}
 Then from~\eqref{utau}, if $\tau\in I_n$, 
 \begin{align*}
 z_\tau &\le \lv \rho_n(\tau)\rv_{y_n(\tau)}
          +\lv \rho_n(t_0)\rv_{y_n(t_0)}
          + B\, u(\tau)
 \\
 &\le  \lv \rho_n(\tau)\rv_{y_n(\tau)}
          +\lv \rho_n(t_0)\rv_{y_n(t_0)}
          +\sqrt{\frac B2}\,e^{B(\tau-t_0)}\,
          \lV\rho_n\rV_{\cL^2([0,T_n])}
          + e^{B(\tau-t_0)}\;
           {\lv \rho_n(t_0)\rv_{y_n(t_0)}}.
 \end{align*}
  Since $B>1$, $\sqrt{\frac B2}<B$. This completes the proof
  of~\eqref{bdsol}.

  \bigskip

  Since $\lim_n T_n=+\infty$, by~\eqref{rhon},
  there exists $n_0>0$ such that if $n>n_0$ then
  \begin{gather*}
  B\,e^{BT_n}\,\lV\rho_n\rV_{\cL^2([0,T_n])}
  \le B\, e^{BT_n}\,\tfrac 1n\,e^{-2 T_n^2}
  <\tfrac 12 \, e^{-T_n^2},
  \\
  (1+ e^{BT_n})\, e^{-\frac 32 T_n^2}
  <\tfrac 12 \,e^{-T_n^2}
   \end{gather*}
   and
   $$
   K_1\,\big(\tfrac 1n\,\sqrt{T_n}\,e^{-2T_n^2}+T_n\,e^{-T_n^2}\big)<1.
   $$
   If $t_0\in A_n$ then $|\rho_n(t_0)|<e^{-\frac 32 T_n^2}$.
   Thus if $n>n_0$ each of the last two terms in~\eqref{bdsol}
   is bounded by $\tfrac 12\,e^{-T_n^2}$.
   This validates~\eqref{bdsol3} when $t\in I_n$.
   
   It remains to prove that if $t_0\in A_n$ then
   $I_n=[0,T_n]$. 
   Let $I_n=[a,b]$. From~\eqref{etb} and~\eqref{bdsol3} we have that 
   for $\tau\in I_n$,
   \begin{align*}
   d_M(x(t),y_n(t))\le e_t &\le \int_{t_0}^\tau K_1\,z_t\;dt
   \le K_1\int_{t_0}^\tau \lv\rho_n(t)\rv\;dt
   +K_1\int_{t_0}^\tau e^{-T_n^2}\;dt,
   \\
   &\le K_1\,\lV\rho_n\rV_{\cL^2([0,T_n])}\,\sqrt{|\tau-t_0|}
   +K_1\,|\tau-t_0|\,e^{-T_n^2}
   \\
   &\le K_1\,\big(\tfrac
   1n\,e^{-2T_n^2}\,\sqrt{T_n}+T_n\,e^{-T_n^2}\big)
   \\
   &<1.
   \end{align*}
  If $b<T_n$, since $t\mapsto e_t$ is continuous,
  $I_n$ could be extended. Then $b=T_n$. Similarly, $a=0$. 
 \end{proof}

 \bigskip

 \begin{Lemma}\label{Hkn}
 If $t\in A_n$ and $n$ is large enough then
 \begin{equation*}
 \big| H(y_n(t),p_n(t))-k\big| < \tfrac 3n
 \qquad \text{ and }\qquad
  \lv p_n(t)\rv\le R_0.
  \end{equation*} 
 \end{Lemma}

 \begin{proof}\quad

 \noindent {\bf Claim 1.}
    $A_n\cap 
   \big\{\,t\in[0,T_n]\,\big|\,|p_n(t)|<R_0\,\big\}
   \ne\emptyset$.
 
 \begincproof
 Let $d_1>|k|$ be from~\eqref{a1}.
 Since $\lV d_{(x_n,T_n)}\cA_k\rV<\tfrac 1n$, 
 $$
 -\left.\frac{\partial\cA_k}{\partial T}\rv_{(x_n,T_n)}
 =\frac 1{T_n}\int_0^{T_n}\big[ H(y_n,p_n)-k\big]\,dt <\frac 1n.
 $$
 Then
 \begin{align*}
 \tfrac 12\, R_0^2\cdot m\big(\big[|p_n|\ge R_0\big]\big)
 &\le\int_0^{T_n}\tfrac 12\,|p_n(t)|^2\;dt
  \le \int_0^{T_n}\big[H(y_n,p_n) +d_1\big]\,dt
 \\
 &\le \big[ k+d_1+\tfrac 1n\big]\,T_n.
 \\
  m\big(\big[|p_n|\ge R_0\big]\big)
 &\le \frac{2\,\big[ k+d_1+\tfrac 1n\big]}{R_0^2}\;T_n
 \\
 &\le \frac 1{10}\;T_n,
 \qquad\text{from \eqref{R0}.}
 \end{align*}
 Therefore, using~\eqref{mAn},
 \begin{align*}
 m(A_n\cap[|p_n|<R_0])
 &= T_n - m\big(A_n^c\cup[|p_n|\ge R_0]\big)
 \\
 &\ge T_n -\tfrac 1n - \tfrac1{10}\;T_n
 \\
 &>0.
 \end{align*}
 This proves Claim~1.
 
 \endcproof

 \noindent{\bf Claim 2.} {\it There exist $K=K(R_0)>0$ and $n_1>0$ 
 such that 
 $$
 \lv p_n(t)\rv_{y_n(t)} \le
 K(R_0)\,\big[\lv\rho_n(t)\rv_{y_n(t)}+1\big],
 \qquad\forall\, t\in[0,T_n], \quad\forall\, n>n_1. 
 $$
 }\indent
 \begincproof
 By Claim~1 there exists $t_0\in A_n$ such that $\lv p_n(t_0)\rv<R_0$.
 Let $(x(t),q(t))$ be the solution of the hamiltonian equations
 with initial conditions $x(t_0)=y_n(t_0)$, $q(t_0)=p_n(t_0)$.
 By Lemma~\ref{ApproxSol}, $d_M\big(x(t),y_n(t)\big)<1$ 
 for all $t\in[0,T_n]$. Given $t\in[0,T_n]$, 
 let $\ga_t:[0,1]\to M$ be the minimizing
 geodesic joining $x(t)$ to $y_n(t)$. 
 Let $\La_t:[0,2]\to TM$ be defined by~\eqref{CurLa},
 let $e_t:=d_M(x(t),y_n(t))$ 
 and let $z_t:=e_t+\lv p_n(t)-A_t(1,0)\cdot q(t)\rv$
 be as in Lemma~\ref{ApproxSol}. Then 
 \begin{align*}
 H(y_n(t),p_n(t))-H&(y_n(t_0),p_n(t_0))
 =H(y_n(t),p_n(t))-H(x(t),q(t))
 \\
 & = \int_0^1\big\langle\nabla_xH(\La_t(s))\,,\,\dga_t(s)
    \big\rangle_{\ga_t(s)}\;ds
 \\
  &\qquad\qquad+\int_1^2\big\langle\nabla_pH(\La_t(s))\,,\,
   p_n(t)-A_t(1,0)\cdot q(t)\big\rangle_{y_n(t)}
  \;ds.
 \end{align*}
 Since $\nabla_xH(x,p)=0$ when $|p|_x>R_0$, there is $d_2>0$ such that
 $$
 \lv \nabla_xH(x,p)\rv_x <d_2\qquad\text{ for all }(x,p)\in TM.
 $$
 Since $\nabla_pH(x,p)=p$ when $|p|_x>R_0$,  there is
 $d_3>0$ such that
 \begin{equation*}
 \lv\nabla_pH(x,p)\rv_x\le |p|_x+d_3\qquad\text{ for all }
 (x,p)\in TM.
 \end{equation*}
 Then
 $$
 \big|H(y_n(t),p_n(t))-H(y_n(t_0),p_n(t_0))\big|
 \le \int_0^1 d_2\, e_t\; ds
 +\int_1^2\big[|p_n(t)|_{y_n(t)}+|q(t)|_{x(t)}+d_3\big]\,z_t\;ds.
 $$
 
 Since $|p_n(t_0)|<R_0$, by Lemma~\ref{R1},  
 $H(y_n(t_0),p_n(t_0))=H(x(t),q(t))<\tfrac 12\, R_0^2$
 and   $|q(t)|<R_0$ for all $t$.
 Let $d_4:=\max\{d_2,\,d_3\}$.
 Then
 \begin{equation}\label{eqH0}
 \big| H(y_n(t),p_n(t))-H(y_n(t_0),p_n(t_0))\big|
 \le \big[|p_n(t)|+R_0+2\, d_4\big]\;z_t.
 \end{equation}

 Suppose first that $|p_n(t)|>R_0$. 
 Then $H(y_n(t),p_n(t))=\tfrac 12\,\lv p_n(t)\rv^2$.
 In this case we have that
 $$
 \tfrac 12\, \lv p_n(t)\rv^2\le H_0+\lv p_n(t)\rv\, z_t
 + [R_0+2\,d_4]\,z_t
 $$
 where  $H_0:=H(y_n(t_0),p_n(t_0))$.
  Using that $H_0<\tfrac 12\,R_0^2$, we get 
 \begin{align*}
 \tfrac 12\,\big[\lv p_n(t)\rv-z_t\big]^2
 &\le H_0 +[R_0+2\,d_4]\,z_t+\tfrac 12\, z_t^2
 \\
 &\le \tfrac 12\,R_0^2+R_0\,z_t+\tfrac 12\, z_t^2+2\,d_4\,z_t
 \\
 &\le\tfrac 12\,[R_0+z_t]^2+ 2\,d_4\,z_t.
 \end{align*} 
 \begin{equation}\label{Npn}
 \lv p_n(t)\rv\le z_t+\sqrt{ [R_0+z_t]^2+4\,d_4\,z_t}\;.
 \end{equation}
 The other case is when $\lv p_n(t)\rv\le R_0$.
 Since the right hand side in~\eqref{Npn} is $\ge R_0$, 
 the inequality~\eqref{Npn} is valid for all $t\in [0,T_n]$.
 Using the identity $(a+b)^2\le 3\,(a^2+b^2)$, we have that
 \begin{align*}
 \lv p_n(t)\rv^2
 &\le 3\,z_t^2+3\,[R_0+z_t]^2+12\,d_4\,z_t
 \\
 &\le d_5\,[\,z_t^2+R_0^2\,]
 &\text{for some constant $d_5>0$},
 \\
 &\le d_6(R_0)\,[\,z_t+\tfrac 12\,]^2
 &\text{for some $d_6(R_0)>0$},
 \\
 &\le d_6(R_0)\,\big[\,|\rho_n(t)|+e^{-T_n^2}+\tfrac 12\,\big]^2
 &\text{using Remark~\ref{RApproxSol} and~\eqref{bdsol3}},
 \\
 &\le d_6(R_0)\,\big[\,|\rho_n(t)|+ 1\,\big]^2
 &\text{if $n$ is large enough}.
 \end{align*}
 Now let $K(R_0):=\sqrt{d_6(R_0)}$.
 
 \endcproof

 \noindent{\bf Claim 3.} {\it There is $n_2>0$ such that
 if $n>n_2$ and  $t_1,t_2\in A_n$  then
 $$
 \lv H(y_n(t_1),p_n(t_1))-H(y_n(t_2),p_n(t_2))\rv <\tfrac 1n.
 $$
 }

 \begincproof
 By Claim~1 and the triangle inequality, it is enough to prove
 that if $t_0\in A_n$, $\lv p_n(t_0)\rv<R_0$ and $t_1\in A_n$
 then 
  $$
 \lv H(y_n(t_1),p_n(t_1))-H(y_n(t_0),p_n(t_0))\rv <\tfrac 1{2n}.
 $$ 

  Let $(x(t),q(t))$ be the solution of the hamiltonian equations
 with initial conditions $x(t_0)=y_n(t_0)$, $q(t_0)=p_n(t_0)$
 and let $z_t$ be as in Claim~2.
 
  If $t_1\in A_n$, $n>n_0$, then $\lv\rho_n(t_1)\rv<e^{-\frac 32 T_n^2}$
 and  by Remark~\ref{RApproxSol} and~\eqref{bdsol3},
 $$
 z_{t_1}\le\lv\rho_n(t_1)\rv+e^{-T_n^2}<2\,e^{-T_n^2}.
 $$
 From~\eqref{eqH0} and Claim~2, we get that
 \begin{align*}
  \big| H(y_n(t_1),p_n(t_1))-H&(y_n(t_0),p_n(t_0))\big|
 \le \big[K(R_0)\,(|\rho_n(t_1)|+1) +R_0+2\, d_4\big]\;z_{t_1}
 \\
 &\le \big[K(R_0)\,(e^{-\frac 32 T_n^2}+1) +R_0+2\, d_4\big]\; 2\,e^{-T_n^2}
 \\
 &<\tfrac 1{2n}
 \end{align*}
 if $n$ is large enough.
 \endcproof
 
 We now finish the proof of Lemma~\ref{Hkn}. Let $d_7>0$ be such that
 $$
 \lv H(x,p)-\tfrac 12\,\lv p\rv_x^2\rv\le d_7
 \qquad\text{for all }(x,p)\in TM.
 $$ 
 Then, using Claim~2, and the inequality $(a+b)^2\le 3(a^2+b^2)$,
 \begin{align*}
 \lv\int_{A_n^c}H\rv
 &\le\int_{A_n^c} \big(\tfrac 12\,\lv p_n(t)\rv_{y_n(t)}^2+d_7\big)\; dt
  \\
 &\le\tfrac 12\,K(R_0)^2 \int_{A_n^c}\big( \lv\rho_n(t)\rv + 1\big)^2\, dt
  + d_7\,m(A_n^c)
  \\
  &\le \tfrac 32\,K(R_0)^2\,\lV\rho_n\rV_{\cL^2}^2 
    + \big[\tfrac 32\,K(R_0)^2+d_7\big]\, m(A_n^c)
  \\
   &\le  \tfrac 32\,K(R_0)^2\,\tfrac 1{n^2}\,e^{-4T_n^2}
    + \big[\tfrac 32\,K(R_0)^2 + d_7\big]\,\tfrac 1{n^2}\,e^{-T_n^2},
 &\text{using~\eqref{rhon} and~\eqref{mAn},}
 \\
 &\le \tfrac 1n, &\text{if $n$ is large enough}.
  \end{align*} 

 Since $\lV d_{(x_n,T_n)}\cA_k\rV<\tfrac 1n$, 
 \begin{gather*}
 \lv\left.\frac{\partial\cA_k}{\partial T}\rv_{(x_n,T_n)}\rv
 =\frac 1{T_n}\lv\int_0^{T_n}\big[ H(y_n,p_n)-k\big]\,dt\rv <\frac 1n,
 \\
 \lv\int_0^{T_n}\big[ H(y_n,p_n)-k\big]\,dt\rv
 \le \frac{T_n}n.
 \end{gather*}
 Therefore
 \begin{align}
 \lv\int_{A_n}H\;dt-k \,m(A_n)\rv
 &\le \lv\int_{A_n^c}H\rv +\lv k\rv \,m(A_n^c)
  +\frac{T_n}n
  \notag\\
 &\le \frac 1n +\frac{|k|}{n^2}\,e^{-T_n^2}+\frac{T_n}n.
 \label{iH-km}
 \end{align}
 
 Let $t\in A_n$. By Claim~3, we have that
 \begin{equation}\label{iH-mH}
 \lv\int_{A_n}H\, dt- m(A_n)\,H(y_n(t),p_n(t))\rv
 \le \tfrac 1n\, m(A_n).
 \end{equation}
 Adding~\eqref{iH-km} and~\eqref{iH-mH} we get that
 \begin{align*}
 \lv H(y_n(t),p_n(t))-k\rv\;m(A_n)
 \le
  \frac 1n\, m(A_n)
 +\frac 1n+\frac{|k|\,e^{-T_n^2}}{n^2}+\frac{T_n}n
 \end{align*}  
 \begin{align*}
 \lv H(y_n(t),p_n(t))-k\rv
 &\le \frac 1n
 +\frac 1{n\,m(A_n)}
 +\frac{|k|\,e^{-T_n^2}}{n^2\,m(A_n)}
 +\frac{T_n}{n\,(T_n-\tfrac 1{n^2}\,e^{-T_n^2})},
 \\
 &\le\frac3n \qquad\text{if $n$ is large enough.}
 \end{align*}
 
 Since by~\eqref{R0}, 
 $\lv H(y_n(t),p_n(t))\rv\le |k|+1<\tfrac 12\,R_0^2$,
 by lemma~\ref{R1}, $|p_n(t)|<R_0$ if $t\in A_n$
 and $n$ is large enough.
\end{proof}

 \bigskip

 Let $\nu_n$ be the Borel probability measure defined by
 $$
 \int f\;d\nu_n=\frac 1{m(A_n)}\,\int_{A_n}f\big(y_n(t),p_n(t)\big)\;dt
 $$
 for any continuous function $f:TM\to\re$. By Lemma~\ref{Hkn} we
 have that
 $$
 \supp(\nu_n)\subseteq H^{-1}
 \big(\big[k-\tfrac 3{n},k+\tfrac 3{n}\big]\big)
 \subseteq H^{-1}\big([k-1,k+1]\big).
 $$
 Since $H^{-1}\big([k-1,k+1]\big)$ is compact, there exists a convergent
 subsequence $\nu_{n_i}$ in the weak* topology. Let
 $$
 \nu:=\lim_i\nu_{n_i}.
 $$
 Then
 $$
 \supp(\nu)\subseteq H^{-1}\{k\}.
 $$

 \bigskip
 
 \begin{Lemma}\label{incon}
 We can assume that $\nu$ is supported on a connected component
 of $H^{-1}\{k\}$.
 \end{Lemma}
 
 \begin{proof}
 If $k$ is a singular value of $H$ then $H^{-1}\{k\}$ contains a 
 singularity of the Hamiltonian flow. In that case a Dirac measure
 supported on the corresponding singularity of the Lagrangian flow
 satisfies the thesis of proposition~\ref{Tinf}.
 
 If $k$ is a regular value of $H$ then there is $\e>0$ such that 
 each of the finitely many connected components 
 of $H^{-1}(]k-\e,k+\e[)$ contains
 exactly one connected component of $H^{-1}\{k\}$. 
 Since the measures $\nu_n$ are supported on 
 the images of the connected curves $y_n$
 and $\supp(\nu)\subseteq H^{-1}(]k-\e,k+\e[)$ for
 $n$ large, we can take the convergent subsequence
 $\nu_{n_i}$ in a single connected component 
 of $H^{-1}(]k-\e,k+\e[)$.
  
 \end{proof}

 \bigskip
 
 \begin{Lemma}\label{invar}
 The probability $\nu$ is invariant under the hamiltonian flow.
 \end{Lemma}

 \begin{proof}
  Given $0<s<1$ let 
 \begin{align*}
 D_n(s):&=\big\{\,t\in[0,T_n]\;\big|\;|\rho_n(t)|<
 e^{-\frac 32 T_n^2},\;|\rho_n(t+s)|<e^{-\tfrac 32 T_n^2}\,\big\},
 \\
 &=A_n\cap (A_n-s).
 \end{align*}
 Then $m(D_n^c\setminus[T_n-s,T_n])
  \le m\big([0,T_n-s]\cap\{A_n^c\cup(A_n-s)^c\}\big)\le
  \tfrac 2{n^2} \, e^{-T_n^2}$, and
 $$
 m(D_n)\ge T_n-s-\tfrac 2{n^2} \, e^{-T_n^2}\ge T_n-2.
 $$
  Let $\psi_t:TM\hookleftarrow$ be the hamiltonian flow.
  Let $F:TM\to \re$ be a continuous function with compact support.
  
  Given $t\in D_n\subset A_n$, by Lemma~\ref{Hkn}, $\lv p_n(t)\rv\le R_0$.
 By Lemma~\ref{ApproxSol} we have that
 \begin{align*}
 d_{TM}\big[\psi_s\big(y_n(t),p_n(t)\big)&,\big(y_n(t+s),p_n(t+s)\big)\big]
 \le
 \lv \rho_n(t+s)\rv+e^{-T_n^2}
 \\
 &\le e^{-\frac 32 T_n^2}+e^{-T_n^2}
 \le \frac 1n,
 \qquad\text{if $n$ is large.}
 \end{align*}
   Since $F$ is uniformly continuous then
 $$
 {\mathcal O}(F,\tfrac 1{n}):=\sup_{d(z,w)<\frac 1{n}}
 \lv F(z)-F(w)\rv
 \xrightarrow{\;n\to\infty\;} 0.
 $$
 Since by~\eqref{mAn} $m(A_n)\ge T_n-1$, $m(D_n^c)\le 2$
  and $\frac{m(D_n)}{m(A_n)}\le 1$,
  we have that
 \begin{align*}
 \lv\,\int F\;d(\psi_s^*\nu_n)-\int F\;d\nu_n\,\rv
  &\le
   {2\,\lV F\rV_\infty} \; \frac{m(D_n^c)}{m(A_n)}
  +\frac 1{m(A_n)}\int_{D_n}\lv F\circ\psi_s-F\rv_{(y_n,p_n)} 
   \\
  &\le
  \frac{4\,\lV F\rV_\infty}{T_n-1}
 +  {\mathcal O}\big(F,\tfrac 1{n}\big)\xrightarrow{\;n\;} 0.
  \end{align*}
  Hence, for all $0<s<1$, 
  $$
  \int F\;d(\psi_s^*\nu)=\int F\;d\nu.
  $$
 \end{proof}

 Let $\mu=\fL_*(\nu)$ be the push forward of $\nu$ under the Legendre transform
 $\fL:TM\to TM$, $\fL(x,p)=\nabla_p H(x,p)$.
 
  \begin{Lemma}
  The homology class $\rho(\mu)\in H_1(M,\re)$ of $\mu$ is zero.
  \end{Lemma}

  \begin{proof}
  If we are working on $\Om_M(q_0,q_1)$,
  let $\ga=\ga_n$ be a minimizing joining the two common endpoints
  $\ga(1)=q_0=y_n(0)$, $\ga(0)=q_1=y_n(T_n)$ of all $y_n$.
  If we are working on $\La_M$, let $\ga_n$ be the constant curve 
  $\ga_n(t)\equiv y_n(0)=y_n(T_n)$, $t\in[0,1]$. 
  
  Let
  $\mu_{n_i}$ be the probability measure defined by
  $$
  \int_{TM}f\;d\mu_{n_i}:=
  \frac 1{T_{n_i}+1}\,
  \left[\,\int_0^{T_{n_i}}f\big(y_{n_i}(s),\dy_{n_i}(s)\big)\;ds
  +\int_0^1f\big(\ga_{n_i}(s),\dga_{n_i}(s)\big)\;ds\;\right]
  $$
  for any continuous function $f\in\cF$ with quadratic growth:
  $$
  \cF:=\left\{f\in C^0(TM,\re)\;\Big|\;\sup_{v\in TM}
  \tfrac{|f(v)|}{1+|v|^2}<+\infty\right\}.
  $$
  We show that for any $f\in\cF$,
  \begin{equation}\label{e6}
  \lim_i\int f\,d\mu_{n_i}=\int f\,d\mu=\lim_i\int f\;d(\fL_*\nu_{n_i})
  =\int f\;d(\cL_*\nu).
  \end{equation}
  We have that 
  $$
  \int f\,d\mu_{n_i}=\frac{m(A_{n_i})}{T_{n_i}+1}\int f\,d(\fL_*\nu_{n_i})
  +\frac 1{T_{n_i}+1}\oint_{\ga_{n_i}}f
  +\frac 1{T_{n_i}+1}\int_{A_{n_i}^c}\!\!\!f\circ \fL.
  $$
  Observe that either $\ga_n$ is a constant curve
  or $\ga_n$ does not depend on $n$. Then
  \begin{align*}
  \lim_i\int f\,d\mu_{n_i}
  &=\lim_i\int f\,d(\fL_*\nu_{n_i})+0+\lim_i\frac 1{T_{n_i}+1}
  \int_{A^c_{n_i}}\!\!\!f\circ \fL
  \\
  &=\int f\,d\mu+\lim_i
  \frac 1{T_{n_i}+1}\int_{A_{n_i}^c}\!\!\!f\circ \fL.
  \end{align*}
  Let 
  $\lV f\rV_{\cF}:=\sup_{v\in TM}\frac{|f(v)|}{1+|v|^2}$.
  Then
  $$
  \int_{A_n^c}f\circ \fL\le\lV f\rV_\cF\,
  \int_{A_n^c}\big[\,1+|v|^2\,\big]\circ \fL.
  $$
   
  Let $d_3>0$ be such that
  $$
  \lv \fL(x,p)\rv_x = \lv \nabla_pH(x,p)\rv_x\le \lv p\rv_x + d_3.
  $$
  Then, using Claim~2 in Lemma~\ref{Hkn}
  and the identity $(a+b)^2\le 3\,(a^2+b^2)$,
  \begin{align*}
  \int_{A_{n_i}^c}\big[\,1+|v|^2\,\big]\circ \fL
  &\le\int_{A_{n_i}^c}1+\big[\lv p_{n_i}(t)\rv_{y_{n_i}(t)}  +d_3\big]^2  
  \\
  &\le (1+3\,d_3^2)\,m(A_{n_i}^c)
  +3\,K(R_0)^2\int_{A_{n_i}^c}\big[\lv\rho_{n_i}(t)\rv+1\big]^2\;dt
    \\
  &\le (1+3\,d_3^2)\,m(A_{n_i}^c)+
  3\,K(R_0)^2\,\big[3\,\lV\rho_{n_i}\rV_{\cL^2}^2+3\,m(A_{n_i}^c)\big]
   \\
  &\le (1+3\,d_3^2)\cdot 2  
      +3\,K(R_0)^2\cdot [3\cdot \tfrac {e^{-4T_{n_i}^2}}{{n_i}^2} + 3\cdot 2]
     \qquad\text{ by \eqref{mAn}  and \eqref{rhon}.} 
    \end{align*} 
  So that
  $$
  \lim_i\frac 1{T_{n_i}+1}\int_{A^c_{n_i}}f\circ \fL=0,
  $$
  and hence~\eqref{e6} holds.
  
  Let $\eta_n$ be the closed curve
  $\eta_n=y_n*\ga_n$ and let $[\eta_n]\in H_1(M,\Z)$ be its homology
  class. Since all the $y_n$'s are in the same free homotopy class,  
  $\a=[\eta_n]$ is constant in $n$.
  Let $\om$ be a closed (bounded) 1-form and $[\om]\in H^1(M,\re)$ its
  cohomology class. Observe that $\om$ has linear growth, in particular
  $\om\in \cF$.
  Then
  \begin{align*}
  \int_{TM} \hskip -.15cm \om\;d\mu
  &=\lim_i\int_{TM}\hskip -.15cm\om\;d\mu_{n_i}
  =\lim_i\frac1{T_{n_i}}\oint_{\eta_{n_i}}\hskip -.2cm\om
  \\
  &=\lim_i\frac 1{T_{n_i}}\langle[\om],\a\rangle 
  =0.
  \end{align*}
  
  \end{proof}
  
    Finally, we prove that the $L+k$ action of $\mu$ is zero.
    Since $L$ is quadratic at infinity then $L\in \cF$.
    By~\eqref{e6}, we have that
    \begin{align*}
    \int_{TM}[L+k]\;d\mu
    &=\lim_i  \int_{TM}[L+k]\;d\mu_{n_i}
    \\
    &= \lim_i \frac {A_{L+k}(y_{n_i})+A_{L+k}(\ga_{n_i})}{T_{n_i}+1}
    \\
   &= \lim_i \frac {\cA_k(x_{n_i},T_{n_i})+A_{L+k}(\ga_{n_i})}{T_{n_i}+1}
   =0.
    \end{align*}
   This finishes the proof of Proposition~\ref{Tinf}.
    \end{proof}  

  \bigskip
 \section{energy levels satisfying the Palais-Smale condition.}
 \label{sSPS}

 In this section we prove corollary~\ref{SPS}. 
 Let $\tL=L:T\tM\to\re$ be the lift
 of $L$ to the universal cover $\tM$ of $M$.
 
 \begin{Lemma}\label{tbdd}
 Given $T_0,A_1>0$ there is $R=R(k,T_0,A_1)$
 such that if $0<T<T_0$ $y\in C^{ac}([0,T],\tM)$ and
 $A_{\tL+k}(y)<A_1$, then\;
 $d(y(0),y(t))\le R$ for all $t\in[0,T]$.
 \end{Lemma}
 
 \begin{proof}
 By the superlinearity there is $b>0$ such that
 $\tL(x,v)\ge|v|_x-b$ for all $(x,v)\in T\tM$.
 We have that
 $$
 d(y(0),y(t))\le\int_0^t|\dy|\;dt
 \le\int_0^T \big[\,\tL(y,\dy)+b \,\big]\;dt
 \le A_1 + (b-k)\,T_0.
 $$
 \end{proof}

 \begin{Lemma}\label{l27}\quad Identify $S^1=[0,1]/_{0\,\equiv\, 1}$.
 
  Let $\si\in[S^1,M]$ be a free homotopy class of closed curves
 in $M$. If $k\ge c_u(L)$ then
 $$
 \inf\{\,\cA_k(x,T)\;|\; x\in\si,\; T>0\,\}>-\infty.
 $$
 \end{Lemma}
 \begin{proof}
 Fix $y\in\si$. Using the homotopy between $x$ and $y$, there are points
 $p\in y([0,1])$, $q\in x([0,1])$ and a curve $z:[0,1]\to M$,
 with $z(0)=p$, $z(1)=q$ such that  
 $z*x*z^{-1}*y^{-1}$ is homotopic
 to a point. We can assume that $x(0)=q$ and $y(0)=p$.
 Then there are lifts $\tx$, $\ty$, $\tz_0$, $\tz_1$ 
 of $x$, $y$, $z$ such that  $\tz_0*\tx*\tz_1^{-1}*\ty^{-1}$
 is a closed curve in $\tM$. Let $\vr$ be the deck transformation
 of the covering $\tM\to M$ such that $\vr(z_0)=z_1$.
 Let $\txn:=\vr^n(\tx)$, $\tyn:=\vr^n(\ty)$ and $\tzn:=\vr^n(\tz)$.
 Since  the curves $\tzn*\txn*\tz_{n+1}^{-1}*\tyn^{-1}=
 \vr^n(\tz_0*\tx*\tz_1^{-1}*\ty^{-1})$ are closed in $\tM$,
 the curves $\tzn*\txn*\tz_{n+1}^{-1}$ and $\tyn$ have the same endpoints.
 Hence the curve $\tw:=\tz_0*(\tx*\tx_1\cdots*\txn)*\tz_{n+1}^{-1}*
 (\tyn^{-1}*\cdots*\ty_1^{-1}*\ty^{-1})$ is closed in $\tM$.
 Given $T>0$ let $S:=1+(n+1)\,T+1+(n+1)$ and $\teta(t):=\tw(t/S)$. 
 Since $k\ge c_u(L)$,
 $$
  A_{\tL+k}(\teta)=\cAk(z,1)+n\,\cAk(x,T)+\cAk(z^{-1},1)
 +n\,\cAk(y^{-1},1)\ge 0.
 $$
 Dividing by $n$,
 $$
 \tfrac 1n\, \cAk(z,1)+\cAk(x,T)+\tfrac 1n \cAk(x^{-1},1)
 \ge -\cAk(y^{-1},1).
 $$
 Letting $n\to+\infty$ we get that for all $T>0$ and $x\in\si$
 $$
 \cAk(x,T)\ge -\cAk(y^{-1},1).
 $$
  \end{proof}
 
 \noindent{\bf Proof of Corollary~\ref{SPS}:}\quad
 Let $\La_1$ be a connected component of
 $\Om_M(q_0,q_1)$ or $\La_M$.
 Let $(x_n,T_n)$ be a sequence in $\La_1$ such that 
 $$
 \lv\cA_k(x_n,T_n)\rv<A_1 \qquad\text{and}\qquad
 \lV d_{(x_n,T_n)}\cA_k\rV<\tfrac 1n.
 $$
 
  \begin{Claim} If the energy level $E^{-1}\{k\}$ does not contain
  singularities of the Euler-Lagrange flow
  and $T_n$ is bounded from above then there is a convergent
  subsequence of $(x_n,T_n)$.
  \end{Claim} 
  
  \begincproof
  By Lemma~\ref{tbdd},
  the curves $\tx_n$ stay in a compact ball $\ov{B(\tqo,R)}$.
  Hence we can apply propositions~\ref{La0} and~\ref{Tfinito}.  
  By Proposition~\ref{La0}.\eqref{La0ii}, $\liminf_n T_n>0$.
  Since $T_n$ is bounded, by proposition~\ref{Tfinito}, 
  $(\tx_n,T_n)$ has a convergent
  subsequence, and so does $(x_n,T_n)$. 

  \endcproof

  Suppose that $k>c_u$. If $(x_n,T_n)\in\La_1\subset\Om_M(q_0,q_1)$
  let $(z_n,S_n)=(x_n*x_0^{-1},T_n+1)\in\La_M$. 
  If $(x_n,T_n)\in\La_1\subset\La_M$ let $(z_n,S_n)=(x_n,T_n)$.
  By Lemma~\ref{l27}, 
  $
  B:=\inf_n
  \cA_{c_u}(z_n,S_n)
  $
  is finite and then
  $$
  A_1\ge \cA_{k}(z_n,S_n)\ge\cA_{c_u}(z_n,S_n)+(k-c_u)\, S_n
  \ge B + (k-c_u)\, S_n.
  $$  
       Hence $S_n$ is bounded and then $T_n$ is bounded. 
       Since $k>c_u\ge e_0$, the energy
       level $E^{-1}\{k\}$ does not contain singularities 
       of the Euler-Lagrange flow.    
       By the claim, $(x_n,T_n)$ has a convergent
       subsequence.


  Now assume that $\La_1\subset\Om_M(q_0,q_1)$.
 Since all the curves $y_n$
 have the same homotopy class with fixed endpoints
 there are lifts $\tqo,\,\tq1\in\tM$
 of $q_0$, $q_1$ and lifts $\tx_n$ of $x_n$ such that
 for all $n$, $(\tx_n,T_n)\in\Om_\tM(\tqo,\tq1)$. 

   Suppose  that $h_{c_u}\equiv +\infty$.
   Since
  $$
  A_1\ge\cA_{c_u}(x_n,T_n)\ge\Phi_{c_u}(\tqo,\tq1;T_n)
  $$
  and 
  $$
  h_{c_u}(\tqo,\tq1)=\liminf_{T\to+\infty}
  \Phi_{c_u}(\tqo,\tq1;T)=+\infty,
  $$
  we have that the sequence $T_n$ is bounded.
  If $E^{-1}\{c_u\}$ contains a singularity $(q_2,0)$ of the
  Euler-Lagrange flow then $L(q_2,0)+c_u=0$ and
  $$
  h_{c_u}(q_2,q_2)\le 
  \liminf_n\int_0^n\big[L(q_2,0)+c_u\big]\;dt =0.
  $$
  This contradicts $h_{c_u}\equiv +\infty$.
  By the claim, $(x_n,T_n)$ has a convergent
  subsequence.

 If $h_{c_u}\not\equiv+\infty$,  
 the same proof of Theorem~C in~\cite{CIPP2}
  applies to our action functional $\cA_{c_u}$
  with our riemannian metric, showing that
  $\cA_{c_u}$ does not satisfy the Palais-Smale
  condition. Indeed, the Palais-Smale sequence
  $(\tx_n,T_n)$ obtained there has $\lim_n T_n=+\infty$
  and is made with solutions of the
  Euler-Lagrange equation joining any two
  given points $\tqo$, $\tq1$  in the universal cover $\tM$.
  Take their projections $x_n:=\pi\circ\txn$.
  Then the curves $x_n$ are in the same homotopy class. 
  Also $\frac{\partial}{\partial x}\cA_{c_u}(x_n,T_n)=0$,
  and the theorem proves that 
  $$
  \lV d_{(x_n,T_n)}\cA_{c_u}\rV
  = \lv\tfrac{\partial\cA_{c_u}}{\partial T}(x_n,T_n)\rv
  =\lv k-E(x_n,\tfrac{\dx_n}{T_n})\rv
  \overset{n}\longrightarrow 0.
  \qquad\qquad\qed
  $$


 \section{The mountain pass geometry.}\label{sMPG}

 In this section we show that a small closed curve of a given length
 $\ell$ inside the projection of the energy level $E^{-1}\{k\}$
 has positive $(L+k)$ action bounded away from zero.
 This gives a mountain pass geometry when we consider families
 of curves going from a constant curve, with arbitrarily
 small action, to a curve with negative 
 action\footnote{This resembles the phrase by Ta\u{\i}manov in~\cite{tai1}:
 ``constant curves are local minimizers of the action'' for magnetic
 flows. But in our case the constant curves are not
 critical points for $\cA_k$, because they don't have energy $k$,
 and also the gradient flow of $-\cA_k$ is not complete on
 a constant curve $(x_0,T)$ because $T$
 reaches zero at a finite gradient flow time.}
 .

 In the case of closed curves in $\La_M$ without a basepoint, we
 need that $k>e_0(L)$, because otherwise, a family of curves could
 move by constant curves until it leaves the projection
 $\pi(E^{-1}\{k\})$, where it already becomes negative,
 without passing through a curve of length $\ell$ inside
 $\pi(E^{-1}\{k\})$.
 
  \begin{Lemma}\label{LU}
 Let $\theta_x$ be a 1-form in $M$. Let $x_0\in M$ and $x_0\in
 V\subset M$ be a neighbourhood of $x_0$.
 Then there exists an open ball  $U\subseteq V$ centered at
 $x_0$ in $M$ and $b>0$ such that if 
 $\ga$ is a closed curve in $U$ then
 $$
 \lv\int_\ga\theta_x\;\rv\le b\cdot \length(\ga)^2.
 $$
  \end{Lemma}
 
 \begin{proof}\quad

 Shrinking $V$ if necessary and
 using a local chart, we can assume that $V$ is a closed ball in
 $\re^m$ with the euclidean metric. 
 Moreover, we can assume that $\ga(0)=0\in V\subset\re^m$.
 Let $b>0$ be such that 
 $|d_x\theta(u,v)|\le b\;|u|\,|v|$ for all $u,\,v\in\re^m$, $x\in V$.
 Let $\ga:[0,T]\to V$ be a closed curve.
  Let $F:[0,1]\times[0,T]\to V$ be defined by
 $F(s,t):=s\,\ga(t)$. Then
 \begin{align*}
 \left|\int_\ga\theta_x\;\right|
 &=\left|\int_{F}d_x\theta\;\right|
 \le\int_0^1\int_0^T b\;\left|\frac{\partial F}{\partial s}\right|\,
 \left|\frac{\partial F}{\partial t}\right|\;dt\,ds
 \\
 &= b\;\int_0^1\int_0^T\lv\ga(t)\rv\cdot s\,\lv\dga(t)\rv\;dt\,ds
 \\
 &\le b\int_0^1\int_0^T\ell(\ga)\cdot s\,\lv\dga(t)\rv\;dt\,ds
 \\
 &\le b\int_0^1\ell(\ga)\cdot s\;\ell(\ga)\;ds
 \le b\cdot\ell(\ga)^2,
 \end{align*}
 where $\ell(\ga)$ is the length of $\ga$.
 Now let $U$ be an open ball for the Riemannian metric
 centered at $x_0$ and contained in $V$.
 
 \end{proof}

   \renewcommand{\theNThm}{\ref{CritVal}}
   
   \begin{NProp}\quad
    \begin{enumerate}
    \item[\eqref{CVloop}]
 Let $x_0\in M$ and $k>E(x_0,0)$.
 Then there exists $c>0$ such that
 if $\Ga:[0,1]\to \Om_M(x_0,x_0)$ is a path joining a constant loop
 $\Ga(0)=x_0:[0,T]\to \{x_0\}\subset M$ (with any $T>0$)
  to any closed loop $\Ga(1)\in\Om_M(x_0,x_0)$ 
  with negative $(L+k)$-action, $A_{L+k}(\Ga(1))<0$, then
 $$
 \sup_{s\in[0,1]} A_{L+k}(\Ga(s)) > c >0.
 $$
 
   \item[\eqref{CVper}]
 Let $k>e_0(L)$.  Then there exists $c>0$ such that
 if $\Ga:[0,1]\to \La_M$ is a path joining any constant curve
 $\Ga(0)=x_0:[0,T]\to \{x_0\}\subset M$ to any closed curve
 $\Ga(1)$ with negative $(L+k)$-action, $A_{L+k}(\Ga(1))<0$, then
 $$
 \sup_{s\in[0,1]} A_{L+k}(\Ga(s)) > c >0.
 $$
 \end{enumerate}
 \end{NProp}

 \begin{proof}\quad
 
  \eqref{CVloop}. Let $d_1\in\re$ be such that
 $$
  E(x_0,0)=-\psi(x_0)<d_1<k.
  $$
  Let $V$ be a neighbourhood of $x_0$ such that 
  $$
  \inf_{x\in V} \psi(x) \ge -d_1.
  $$
  Write, as in Lemma~\ref{Lo},
  $$
   L(x,v)\ge \tfrac 12\,a\,|v|_x^2+\theta_x(v)+\psi(x),
  $$
   where
  $\theta_x(v):=L_v(x,0)\cdot v$,\;
  $\psi(x) :=L(x,0)$
 and
 $a:=\inf_v \{v\cdot L_{vv}(x,v)\cdot v\}/{|v|_x^2}>0$.

  Let $U\subseteq V$ be an open ball centered at $x_0$ 
  given by Lemma~\ref{LU} for $(x_0,V)$.  
  Let
 \begin{equation}\label{lo}
 0<\ell_0<
 \min\left\{\tfrac 12\,\diam(U),\;\sqrt{\frac{a\,(k-d_1)}{2\,b^2}}\;\right\}.
 \end{equation}

 \begin{Claim}
 There exists $0<s_0<1$ such that
 $\length(\Ga(s_0))=\ell_0$. 
 \end{Claim}
 
 \begincproof
 Suppose that $\Ga(s_1)\not\subset U$ for some $s_1\in[0,1]$.
 Since  $s\mapsto\length(\Ga(s))$ is continuous,
 $\length(\Ga(0))=0$ and 
 $\length(\Ga(s_1))\ge d(x_0,U^c)\ge \tfrac 12\,\diam(U)>\ell_0$; 
 then there exists $0<s_0<s_1$ such that 
$\length(\Ga(s_0))=\ell_0< \tfrac
 12\diam(U)$.

 Now assume that $\Ga(s)\subset U$ for all $s\in [0,1]$.
 Writing $\ga_1:=\Ga(1):[0,T_1]\to M$ and 
 $\ell_1=\length(\ga_1)$. By lemma~\ref{LU} we have that
 \begin{align*}
 0> A_{L+k}(\ga_1)
 &\ge \frac 12\int_0^{T_1} a\,|\dga_1|^2\;dt
 -\lv\int_{\ga_1} \theta_x\;\rv +\int_0^{T_1}\psi(\ga_1(t))\;dt
 +k\,T_1
  \notag
 \\
 &\ge \frac a2 \int_0^{T_1}\lv\dot{\ga_1}\rv^2\;dt
      -b\,\ell_1^2 + (k-d_1)\,T_1.
      \label{cotaA}
 \end{align*}
 By the Cauchy-Schwartz inequality,
 $$
 T_1\;\int_0^{T_1}\lv\dga_1\rv^2\;dt
 \ge\left(\int_0^{T_1}\lv\dga_1\rv\;dt\right)^2
 = \ell_1^2.
 $$
 Hence
 \begin{equation}\label{ell1}
 0>A_{L+k}(\ga_1)\ge\left(\frac a{2 \,T_1} -b\right)\,\ell_1^2 
 + (k-d_1)\,T_1.
 \end{equation}
 Since $(k-d_1)>0$ and $T_1 >0$,  then $\frac a{2T_1} -b < 0$, i.e.
 $$
 T_1 > \frac{a}{2\,b}.
 $$
 From~\eqref{ell1}, we have that
 $$
 \ell_1^2 >\frac{(k-d_1)\, T_1}{b-\tfrac{a}{2T_1}}
 >\frac{(k-d_1)\,T_1}{b}
 >\frac {a\,(k-d_1)}{2\, b^2}
 >\ell_0^2.
 $$
 Since $\length(\Ga(0))=0$, there is $s_0\in [0,1]$
 with $\length(\Ga(s_0))=\ell_0$.
 \endcproof
 
  Since $\length(\Ga(s_0))=\ell_0<\tfrac 12\diam U$, $U$ is an open
  ball centered at $x_0$ and $\Ga(s_0)\in\Om_M(x_0,x_0)$,
  we have that $\Ga(s_0)\subset U$. In particular, 
  the right estimate in~\eqref{ell1} holds for $\Ga(s_0)$.
  Let
 $$
 f(t):=\left(\frac a{2\,t}-b\right)\,\ell_0^2 + (k-d_1)\,t.
 $$
 If $\Ga(s_0):[0,T_0]\to M$, then
 $$
 A_{L+k}(\Ga(s_0))\ge f(T_0)
 \ge \min_{t\in\re^+}f(t)
 =\ell_0\left[\sqrt{2a\,(k-d_1)}-b\,\ell_0\right]=:c >0,
 $$
 because
 $$
 \ell_0<\sqrt{\frac{a\,(k-d_1)}{2\,b^2}}<\frac{\sqrt{2\,a\,(k-d_1)}}b.
 $$
 
 \eqref{CVper}. Since $k>e_0(L)$, 
 \begin{equation}\label{ke0}
  k>\sup_{x\in M}E(x,0)=e_0(L).
 \end{equation}
  Let $U_1,\ldots,U_N$ be a finite cover of $M$
  by open balls given by Lemma~\ref{LU} with corresponding constants
  $b_i=b_i(U_i)>0$. Let $r_0>0$ be such that
  any ball of radius $r_0$ in $M$ is contained 
  in one $U_i$. Write $b=\max_{1\le i\le N}b_i$ and
  let
   \begin{equation}\label{lo1}
 0<\ell_0<
 \min\left\{r_0 \;,\;\sqrt{\frac{a\,(k-e_0)}{2\,b^2}}
 \;\right\}.
 \end{equation}
 Let $\Ga:[0,1]\to \La_M$ be a path joining a constant curve
 $\Ga(0)=x_0:[0,T]\to \{x_0\}\subset M$ to a closed curve 
 $\Ga(1)$ with negative $(L+k)$-action, $A_{L+k}(\Ga(1))<0$.

 \begin{Claim}
 There is $s_0\in[0,1]$ such
 that $\length(\Ga(s_0))=\ell_0$.
 \end{Claim}

   \begincproof

 Suppose that $\length(\Ga(s))<\ell_0$ for all $s\in[0,1]$.
 Since $\ell_0<r_0$,  for each $s\in[0,1]$,
 $\Ga(s)$ is contained in some $U_i$.
   From~\eqref{ke0} we have that
 $\psi(x)=-E(x,0)\ge -e_0$ for all $x\in M$.
  Let $\ell_1:=\length(\Ga(1))$.
 Then the same argument as in item~\eqref{CVloop} proves that
 $$
 \ell_1^2>\frac {a\,(k-e_0)}{2\, b^2}
 >\ell_0^2.
 $$
 Therefore there is  $s_0\in[0,1]$ such that 
 $\length(\Ga(s_0))=\ell_0$. 
 \endcproof

 Let
 $$
 g(t):=\left(\frac a{2\,t}-b\right)\,\ell_0^2 + (k-e_0)\,t.
 $$
 Since $\length(\Ga(s_0))=\ell_0<r_0$,  $\Ga(s_0)$ 
 is contained in some $U_i$ and we can apply Lemma~\ref{LU}.
 Therefore, if  $\Ga(s_0):[0,T_0]\to M$, then
 $$
 A_{L+k}(\Ga(s_0))\ge g(T_0)
 \ge \min_{t\in\re^+}g(t)
 =\ell_0\left[\sqrt{2a\,(k-e_0)}-b\,\ell_0\right]=:c >0,
 $$
 because
 $$
 \ell_0<\sqrt{\frac{a\,(k-e_0)}{2\,b^2}}<\frac{\sqrt{2\,a\,(k-e_0)}}b.
 $$

 \end{proof}

 \section{Some results on Morse theory}
 \label{sSRMT}
 \newcounter{cmorse}

 Let $X$ be an open set in a
 Riemannian manifold and 
 $f:X \to \re$ be a $C^2$ map. 
 Observe that if the vector field $-\nabla f$ 
 is not globally Lipschitz, the gradient flow $\psi_t$ 
 of $-f$ is a priori only a local flow. Given 
 $p\in X$, $t>0$, define
 \begin{align*}
 \a(p)&:=\sup\{\, a>0\,\vert\, s\mapsto\psi_s(p) 
        \text{ is defined on }s\in [0,a]\,\}.
 \end{align*}
  We say that the flow $\psi_t$ of $-\nabla f$ is
  {\it relatively complete on $[a\le f\le b]$}  
  if for $a\le f(p)\le b$,  either
 $\a(p)=+\infty$ or $f(\psi_\be(p))\le a$ for
 some $0\le \be < \a(p)$.

  We say that a function $\tau:X\to[0,+\infty[$
 is an {\it admissible time\/} if $\tau$ is differentiable and
 $0\le \tau(x)<\a(x)$ for all $x\in X$. 
 Given and admissible time $\tau$ and a subset $F \subset X$,
 define
 $$       
 F_\tau:=\{\,\psi_{\tau(x)}(x)\,\vert\, x\in F\,\}.         
 $$

 Given a closed subset $B\subseteq X$, 
 we say that the function $f$ satisfies the
 {\it Palais-Smale condition restricted to $B$ at level $c$}, 
 $(PS)_{c,B}$, if any  sequence 
 $\{x_n\} \subset B$ with $\lim_n\lV df(x_n)\rV=0$
 and $\lim_nf(x_n)=c$
 has a convergent subsequence.

 Given $c\in\re$, $\de>0$ and a closed subset $B\subseteq X$, 
 define
   \begin{align*}
   K_{c,B}:&=\{\,x\in B\,|\,f(x)=c,\;df(x)=0\,\},
    \\
   W_{c,\de,B}:&=\{\,x\in B\,|\, d(x,K_{c,B})<\de\,\},
    \\
   V_{c,\de,B}:&=\{\,x\in B\,|\, \lV df(x)\rV<\de,\;|f(x)-c|<\de\,\}. 
   \end{align*}
 
  \begin{Lemma}\label{nbhds}
  Let $X$ be a Riemannian manifold and $f:X\to\re$ a $C^1$ function.
  
  If $f$ satisfies the Palais-Smale condition $(PS)_{c,B}$ at level $c$
  restricted to $B$,  then 
   \renewcommand{\theenumi}{\roman{enumi}}
  \begin{enumerate}
  \item\label{nbhd0}
    $K_{c,B}$ is compact.
  \item\label{nbhd1}
     The family $\{W_{c,\de,B}\}_{\de>0}$  
     is a fundamental system of neighbourhoods
     of $K_{c,B}$ relative to $B$. 
  \item\label{nbhd2}
     The family $\{V_{c,\de,B}\}_{\de>0}$
     is a fundamental system of neighbourhoods 
     of $K_{c,B}$ relative to $B$.
  \end{enumerate}
  \end{Lemma}
 
  \begin{proof}
  \quad
  
  \eqref{nbhd0}. By $(PS)_{c,B}$ any sequence in $K_{c,B}$ 
                has a convergent
                subsequence. Since $df$ is continuous and $B$
		is closed, 
		the limit is also in $K_{c,B}$.
		
   \eqref{nbhd1}. Suppose item~\eqref{nbhd1} is false.
                  Then there is a relative neighbourhood $U$ of
		  $K_{c,B}$ with $K_{c,B}\subset U \subset B$, 
		  and a sequence $x_n\in  W_{c,1/n,B}\cap U^c$.		
		  Then there is a sequence $y_n\in K_{c,B}$ such that
		  $d(x_n,y_n)\le \tfrac 2n$. Since $K_{c,B}$ is compact
		  there is a convergent subsequence $z=\lim_k
		  y_{n_k}\in K_{c,B}$. Also, $\lim_k x_{n_k}=z\in
		  K_{c,B}$.
		  This contradicts  $x_n\notin U$ for all $n$.
		
   \eqref{nbhd2}. Suppose item~\eqref{nbhd2} is false.
                  Then there is a relative neighbourhood 
		  $U$ of $K_{c,B}$ with $K_{c,B}\subset U\subset
		  B$
		  and a sequence $x_n\in  V_{c, 1/n,B}\cap U^c$.
		  By $(PS)_{c,B}$ there is a convergent subsequence
		  $z=\lim_k x_{n_k}$. Since $f\in C^1$,
		  $z\in K_{c,B}$. This  contradicts
		   $x_n\notin U$ for all $n$.
		
  \end{proof}

  \begin{Lemma}\label{Deform}
  Let $X$ be a Riemannian manifold
  and $B\subset A\subset X$  closed subsets
  such that $A$ contains the $\e_1$-neighbourhood of B:
  $$
  \{\;x\in X\;|\;\exists\, y\in B,\;d(x,y)<\e_1\;\}
  \subset A.
  $$
  Let $f:X\to\re$ be a $C^2$ function. Let $c\in\re$
  and suppose that 
  \begin{enumerate}
  \renewcommand{\theenumi}{\roman{enumi}}
  \item $f$ satisfies the Palais-Smale condition 
        $(PS)_{c,A}$ at level  $c$ restricted to $A$.
  \item\label{relcoii} The flow of $-\nabla f$ is relatively complete
        on $[|f-c|\le \e_2]$ for some $\e_2>0$.	
  \end{enumerate}  
  Given any neighbourhood $N$ of $K_{c,A}$ relative to $A$ and $\e_3 >0$,
  there are  $0<\e<\de<\e_3$, 
  such that for all $0<\la\le \e$ there is 
  an admissible time  $\tau(x)$ such that
  $$
  F_\tau \subseteq N \cup [\,f\le c-\la\,]
  \qquad \text{ and } \qquad
   \tau(x)=0 \text{ on } [\,|f-c|\ge\de\,],
  $$
  where $F=\big([\,f\le c+\e\,]\cap B\big)\cup [f\le c-\la]$.
  \end{Lemma}
  
   \begin{proof}
   By Lemma~\ref{nbhds}, there are
    $0<\de,\,\rho,\,\eta<\min\{\,\e_2,\e_3,1\,\}$ such that
   $$
     V_{c,\de,A}\subset W_{c,\rho,A}\subset W_{c,2\rho,A}
     \subset V_{c,\eta,A}
     \subset N\subset A.
   $$
   Let 
   $$
   \e:=\tfrac 12\,\min\left\{\,\tfrac{\e_1^2}2,\,\e_2,\, \e_3,\,
   \tfrac{\de^2}2,\,\tfrac{\rho \,\de}{2}
   \,\right\}.
   $$
   Let $h:X\to[0,1]$ be a smooth function such that
   $$
   h(x)=\begin{cases}
      0 &\text{ if}\quad x\in \;\,V_{c,\frac \de 2,A}\;\,\cup [|f(x)-c|>\de], \\
      1 &\text{ if}\quad x\in(V_{c,\de,A})^c\cap [|f(x)-c|<\frac \de 2].
         \end{cases}
   $$

  Since $f$ is $C^2$, the vector field $Y(x):=-h(x)\,\nabla f(x)$
  is locally Lipschitz. Since $\de<\e_2$, by the hypothesis~\eqref{relcoii},
  the flow $\eta_s$ 
  of $Y$ is  complete.
  
  Let $\psi_t$ be the flow of $-\nabla f$.
  Define $\tau(x)$ by $\psi_{\tau(x)}(x)=\eta_1(x)$. We show that
 $\tau(x)$ is an admissible time. 
 Write $\eta_s(x)=\psi_{t(s)}(x)$,
 then
 $$
 Y(\eta_s(x))=
 -h(\eta_s(x)) \; \nabla f(\eta_s(x))
  = -\tfrac{dt}{ds}\;\nabla f(\eta_s(x)).
 $$
 So that
 $\tfrac{dt}{ds}=h(\eta_s(x))$ and
 $$
 \tau(x)=\int_0^1 h(\eta_s(x))\; ds\le 1.
 $$
 Therefore $\tau(x)$ is finite and differentiable.
 
 Let $x\in F$. If $x\notin B$ then $f(x)\le c-\la$ 
  and hence
 $f(\eta_1(x))\le f(x)\le c-\la$.
 
 Let $x\in F\cap B$.
 We can assume $|f(\eta_s(x))-c|<\frac{\de}2$ for all
 $s\in[0,1]$. For, if not, since $s\mapsto f(\eta_s(x))$ is
 non-increasing and at $s=0$, $f(x)\le c+\e<c+\frac\de 2$, then 
 we already have that
 $$
 f(\psi_{\tau(x)}(x))=f(\eta_1(x))
 =\inf_{s\in[0,1]}f(\eta_s(x))\le c-\tfrac \de 2\le c-\e.
 $$

  We can also assume that 
  \begin{equation}\label{inA}
  \eta_s(x)\in A\quad\text{ for all }
  s\in [0,1].
  \end{equation}
  For, if not, since $x=\eta_0(x)\in F\cap B$ and $A$ contains
  the $\e_1$-neighbourhood of $B$,
  \begin{align*}
  \e_1&\le \length\big[\eta_{[0,1]}(x)\big]
  =\int_0^1\lV h(\eta_s(x))\;\nabla f(\eta_s(x))\rV\; ds
  \\
  &\le \left[\int_0^1 \lv h\rv^2\,\lV\nabla f\rV^2\,ds\right]^{\frac 12}
  \le \left[\int_0^1 h\;\lV\nabla f\rV^2\;ds\right]^{\frac 12}.
  \end{align*}
  And then
  \begin{align*}
  f(\eta_1(x))&=f(x)-\int_0^1\langle \nabla f,\, Y\rangle\;ds
  \le c+\e-\int_0^1 h\,\lV\nabla f\rV^2\;ds
  \\
  &\le c+\e-\e_1^2<c-\e.
  \end{align*}

  Suppose that $x\in F\cap B$ and 
  $\eta_s(x)\notin V_{c,\de,A}$ for all $s\in[0,1]$.
  Then $h(\eta_s(x))=1$ for all $s\in[0,1]$.
  By~\eqref{inA}, 
  $\eta_s(x)\in A\cap (V_{c,\de,A})^c$ for all $s\in[0,1]$,
  then $\lV\nabla f(\eta_s(x))\rV\ge \de$
  for all $s\in[0,1]$.  Then, 
  \begin{align*}
  f(\eta_{1}(x))
  &\le c+\e-\int_0^1h(f(\eta_t(x)))\;\lV\nabla
  f(\eta_t(x))\rV^2\;dt
  \\
  &\le c+\e -\de^2 \le c-\e.
  \end{align*}
 
 Now suppose that $x\in F\cap B$, $\eta_{s_0}(x)\in V_{c,\de,A}$
 for some $s_0\in[0,1]$ and $\eta_1(x)\notin N$.
 Let
 \begin{align*}
 s_1:&=\inf\{ \, s> s_0\,|\, \eta_s(x)\notin V_{c,\de,A}\,\},
 \\
 s_2:&=\inf\{ \, s> s_1\,|\, \eta_s(x)\notin V_{c,\eta,A}\,\}
 \le 1.
 \end{align*}
 By~\eqref{inA}, the image of $[0,1]\ni s\mapsto \eta_s(x)$ is in $A$.
 Since the segment $[s_1,s_2]\mapsto \eta_s(x)$  crosses
 the annulus of width $\rho$:
 $W_{c,2\rho,A}\setminus W_{c,\rho,A}$, inside of $A$;
 we have that
 \begin{align*}
 \rho\le\length\left[\eta_{[s_1,s_2]}\right]
 &\le\int_{s_1}^{s_2} h(\eta_s(x))\,\lV\nabla f(\eta_s(x))\rV\;ds
 \\
 &\le \frac 1\de \int_{s_1}^{s_2} h\,\lV\nabla f\rV^2\;ds,
 \qquad\text{because }\; \eta_{]s_1,s_2[}(x)\subset A\cap
 (V_{c,\de,A})^c.
 \end{align*}
 Since $s_2\le 1$, 
 \begin{align*}
 f(\eta_1(x)) &\le c+\e-\int_0^1 h\,\lV\nabla f\rV^2\;ds
 \\
 &\le c+\e- \rho\,\de\le c-\e.
 \end{align*}
 Therefore
 $$
 F_\tau=\eta_1(F)  
 \subseteq N\cup [f\le c-\e].
 $$
 \end{proof}

  Given a function $f:X\to\re$ on a topological space $X$, 
 we say that $x\in X$ is a {\it strict local minimizer} of $f$ if
 there is a neighbourhood $V$ of $x$ in $X$ such that
 $f(y)> f(x)$ for all  $y\in V\setminus\{x\}$.

 Let $\cF$ be a family of subsets $F \subset X$.
 We say that $\cF$ is {\it forward invariant} if
 $F_\tau\in\cF$ for all $F\in\cF$ and any admissible time 
 $\tau$. Define
 $$
 c(f,{\cF})= 
 \inf_{F \in {\cF} }\;
  \sup_{x\in F}\; f(x).
  $$
  
 \begin{Proposition}\label{minimax}
 Let $X$ be a Riemannian manifold and $f:X\to\re$ a $C^2$ function.
 Let $\cF$ be a family of subsets of $X$. Suppose that
   \renewcommand{\theenumi}{\roman{enumi}}
   \begin{enumerate}
  \item\label{hypcon} The subsets $F\in\cF$ are connected.
  \item\label{hypfow} $\cF$ is a forward invariant family. 
 \item $c:=c(f,\cF)\in\re$.
 
 \item The flow of $-\nabla f$ is relatively complete on 
 $[c-\e_2\le f\le c+\e_2]$ for some $\e_2>0$.
 \item \label{hypB}
        There is a closed subset $B\subset X$ such that:
       \newline
       $\forall \e>0,\;\exists \la\in]0,\e[,\;\exists F\in\cF$ such that
       $F\subset B\cup[f\le c-\la]$ and $F\subset[f\le c+\e]$.
 \item There is $\e_1>0$ and  a closed subset $A\subset X$
        which contains the $\e_1$-neighbourhood of $B$
	such that
	$f$ satisfies the Palais-Smale condition $(PS)_{c,A}$,
                restricted to $A$, at level $c$.
 \end{enumerate}
 Then $K_{c,A}\ne\emptyset$, i.e.
 $f$ has a critical point $\ov{x}$ in $A$ with
 $f(\ov{x})=c$.
 
 Moreover, if
 \begin{equation}\label{extremes}
 \sup_{F\in\cF}\inf_{x\in F} f(x) <c,
 \end{equation}
 then there is a point in $K_{c,A}$ which is not a strict local
 minimizer.
 \end{Proposition}

\begin{Remark}\label{Rminimax}
 It is enough to consider admissible times $\tau(x)$ such that
 $\tau(x)=0$ if $f(x)\le c-\de$ for some $\de>0$.
 For, the value of $c(f,\cF)$ does not change and
 the proof of Proposition~\ref{minimax} only uses that
 kind of admissible times.
\end{Remark}

 \begin{proof}
 Suppose that $K_{c,A}$ is empty.
 Let $\e,\,\de>0$ be given by Lemma~\ref{Deform} for $N=\emptyset$
 and $\e_1$, $\e_2$, and $\e_3=1$.
 By the hypothesis~\eqref{hypB},
  there is  $F\in\cF$ such that
 $F\subset B\cup [f\le c-\la]$ and $F\subset [f\le c+\e]$. 
 By Lemma~\ref{Deform}
 there is an admissible time $\tau$, with $\tau=0$ on 
 $[f\le c-\de]$, such that $F_\tau\subset [f\le c-\la]$.
 This contradicts the definition of $c(f,\cF)$.

 Now suppose that $K_{c,A}$ consists entirely 
 of strict local minimizers of $f$ and that 
 inequality~\eqref{extremes} holds. 
 Let $\e_0>0$ be such that
 \begin{equation}\label{extremo2}
 \sup_{F\in\cF}\inf_{x\in F} f(x)< c-2\,\e_0.
 \end{equation}
 
  For each $x\in K_{c,A}$ let $N(x)$ be a
  neighbourhood of $x$ such that $f(y)> f(x)$ for all $y\in
  N(x)\setminus\{x\}$.
  Let $$
  N_0:={\textstyle\bigcup\limits_{x\in K_{c,A}}} N(x).
  $$ 
  Let $N:=A\cap N_0$.
  
  Let $0<\e<\de<\e_0$ be given by Lemma~\ref{Deform} for
  $c$, $B$, $A$ and $N_0$.
  By hypothesis~\eqref{hypB} there are $\la\in]0,\e[$ and $F\in\cF$
  such that $F\subset B\cup [f\le c-\la]$ and
   $F\subset [f\le c+\e]$.
    By Lemma~\ref{Deform},  there is an admissible time $\tau$
  such that  $\tau\vert_{[f\le c-\de]}\equiv 0$ and
  $$
  F_\tau
  \subseteq N\cup [f\le c-\la]
  \subseteq N_0\cup [f\le c-\la].
  $$
  
   By definition of $N_0$, the sets $N_0$ and $f\le c-\la$ are disjoint,
  in particular, disconnected. By hypothesis~\eqref{hypcon} and
  ~\eqref{hypfow}, $F_\tau\in\cF$ is connected.
  Then the set $F_\tau$ lies either
  in $F_\tau\subset N$ or $F_\tau\subset [f\le c-\la]$.
  Since $\la\le \e<\e_0$ and the value of $f$ decreases under the
  flow of $-\nabla f$, by~\eqref{extremo2},
  $F_\tau\cap[f\le c-\la]\ne \emptyset$.
  Hence $F_\tau\subset [f\le c-\la]$.
  This contradicts the definition of $c(f,\cF)$.
  
 \end{proof}

 We shall use   the following ``mountain pass'' theorem.
    
 \begin{Corollary}\label{MountainPass}

 Let $X$ be a $C^2$ Riemannian manifold and $f:X\to\re$ a $C^2$
 function.  
 Let $p,q\in X$ and
 $$
 c:=\inf_{\ga\in\Ga}\sup_{s\in[0,1]} f(\ga(s)),
 $$
 where $\Ga:=\{\,\ga:[0,1]\to X\,|\,\ga\in C^0,
 \;\ga(0)=p,\;\ga(1)=q\,\}$.
 
 Suppose that 
 \renewcommand{\theenumi}{\roman{enumi}}
 \begin{enumerate}
  \item $c\in\re$.
 \item The flow of $-\nabla f$ is relatively complete on 
 $[c-\e_2\le f\le c+\e_2]$ for some $\e_2>0$.
  \item $\max\{\,f(p),\;f(q)\,\}<c$.
 \item There are closed subsets $B\subset A\subset X$ such that
 \begin{enumerate}
 
 \item $f$ satisfies the Palais-Smale condition $(PS)_{c,A}$,
 restricted to $A$, at level $c$.
 \item For some $\e_1>0$,
         $A$ contains the $\e_1$-neighbourhood of $B$.
 \item For all $\e>0$, there are $\la_\e\in]0,\e[$ 
       and $\ga_\e\in\Ga$ such that
        $$
	\ga_\e([0,1])\subset \big(B\cup[f\le c-\la_\e]\big)
	\cap[f\le c+\e].
	$$	 
 \end{enumerate}
 \end{enumerate}
 Then $c$ is a critical value of $f$.
 Moreover the set
 $$
 K_{c,A}:=\{\; x\in X\;|\; x\in A,\; df(x)=0,\; f(x)=c\;\}
 $$
 contains a point which is not a strict local minimizer.
 \end{Corollary}

 \begin{proof}
 Let
 $\cF:=\{\,\ga([0,1])\,|\,\ga\in C^0([0,1],X),\;\ga(0)=p,\;\ga(1)=q\,\}$.
 Let $\de>0$ be such that $\max\{f(p), f(q)\}<c-2\,\de$. 
 Then $c=c(f,\cF,B)$ and the family $\cF$ is forward invariant
 when we consider only admissible times $\tau(x)$ such that
 $\tau(x)=0$ when $f(x)\le c-\de$.
 Then Proposition~\ref{minimax} and Remark~\ref{Rminimax}
 prove that $K_{c,A}\ne\emptyset$ and that $K_{c,A}$ contains a point
 which is not a strict local minimizer. 
 
 \end{proof}

 Applying Proposition~\ref{minimax} to the family
 of subsets $F=\{p\}$, $p\in X$, and $B=A=X$, we obtain

 \begin{Corollary}\label{minimum}\quad
 
 \noindent
  Let $X$ be a $C^2$ Riemannian manifold and $f:X\to\re$ a $C^2$
 function. Suppose that 
 \begin{enumerate}
 \renewcommand{\theenumi}{\roman{enumi}}
 \item $c:=\inf_{x\in X}f(x)>-\infty$.
 \item The flow of $-\nabla f$ is relatively complete on 
 $[c-\e\le f\le c+\e]$ for some $\e>0$.
 \item $f$ satisfies the Palais-Smale condition $(PS)_c$ at level $c$.
 \end{enumerate}
 Then $c$ is a critical value of $f$.
 
 \end{Corollary}

 \medskip

  \bigskip

  Now we concentrate on the relative completeness condition
  for the action functional $\cA_k$.
  For completeness we present the following statement:

 \medskip
 \begin{Lemma}\label{RelComp0}
 Suppose that  $f:X\to\re$ is $C^2$, 
 $\psi_t$ is the gradient flow of $-f$
 and the subset $[a\le f\le b]\subset X$ is complete.
 Then the flow  $\psi_t$ is relatively complete 
 on $[a\le f\le b]$.
 \end{Lemma}

 This lemma can be found in~\cite[Lemma~22]{CIPP2}, but its proof
 is similar to the first part of the proof of Lemma~\ref{RelComp}     
 below.
 Recall from remark~\ref{NCt}
 that the spaces $\La_M$ and $\Om_M(q_0,q_1)$ are not complete
 with our riemannian metric.
 Also the gradient flow of $-\cA_k$ is not complete
 as the following example shows.
 
 \begin{eexample}\label{cmetr}
 The gradient flow of $-\cA_k$ is not complete.
 \end{eexample}
  Let $(q_0,0)$ be a fixed point of the Euler-Lagrange
 flow, let $x_0:[0,1]\to\{q_0\}$ be the constant curve, let $T>0$ and 
 let $y(t)=x(t/T)$. Since $x_0$ is a solution of the 
 Euler-Lagrange equation, from~\eqref{deriv}, we have that
 the partial derivative  $(\partial_x\,\cA_k)(x,T)=0$. But
 $$
 \left.\frac{\partial \cA_k}{\partial T}\right|_{(x,T)}
 =\int_0^1\big[k-E(q_0,0)\big]\;dt=k-E(q_0,0).
 $$
 Suppose that  $k>E(q_0,0)$. Let $a:=k-E(q_0,0)>0$. Since the
 metric on the $\re^+$ factor of $\cH^1(M)\times\re^+$ is the
 euclidean metric, the gradient flow $\Psi_s$ of $-\cA_k$
 on $(x,T)$ has differential equations 
 $\frac{dx}{ds}=0$ and  $\frac{dT}{ds}=-a$. Then
  $$
 \Psi_s(x,T) = (x, T-as),
 $$
 which leaves the space $\cH^1(M)\times\re^+$
 through $(x,0)$ in finite time.
  
 One could change the riemannian metric~\eqref{1metric}
 in such a way that $\cH^1(M)\times\{0\}$ lies ``at infinity'',
 for example replacing the first term in~\eqref{1metric},
 \eqref{metric} by $\a \be/T$. In that case the Hilbert manifolds
 $\La_M$ and $\Om_M(q_0,q_1)$ become complete and the
 gradient flow of $-\cA_k$ becomes relatively
 complete by Lemma~\ref{RelComp0}.
 Indeed, in the example above we would have
 $$
 \left.\frac{\partial \cA_k}{\partial T}\right|_{(x,T)}\cdot\a
 =a\,\a=\frac 1T\;(a \,\a\,T)=\langle a\,T,\a\rangle_T.
  $$ 
  Thus the projection $\nabla_T\cA_k$ to the $\re^+$ component of
  $\cH^1(M)\times\re^+$ of the gradient of $\cA_k$ is $a\,T$.
  The differential equations for $\Psi_s$ become
  $\frac{dx}{ds}=0$,   $\frac{dT}{ds}=-a\,T$ 
  and $\psi_s(x,T)=(x,T e^{-as})$.
  
  Nevertheless, such a change in the riemannian metric
  would give a weaker hypothesis in Proposition~\ref{La0}:
  inequality~\eqref{mEn} would have a factor $T_n$ on the
  right hand side. 
  Indeed, Proposition~\ref{La0} would be false.
  Items~\eqref{La0i},~\eqref{La0ii} and~\eqref{La0iv} 
  would still hold but not item~\eqref{La0iii}.
  In the example above, the sequence 
  $(x_n,T_n)=\Psi_n(x,T)=(x,T e^{-an})$ would satisfy 
  \begin{gather*}
  \cA_k(x_n,T_n)=T e^{-an}
  \big[L(x,0)+k\big]\overset{n}\longrightarrow 0,
  \\
  \lV d_{(x_n,T_n)}\cA_k\rV
  =\frac{\lv\,\partial_T\cA_k(x_n,T_n)\cdot\a\,\rv}
  {\lV (0,\a)\rV_{(x_n,T_n)}}
  =\frac {\lv\a\cdot \int_0^1\big[k-E(x,0)]\;ds\rv}
   {\lv\a \rv/{T_n}}
   = T e^{-an}\cdot a \overset{n}\longrightarrow 0,
    \end{gather*}
    so that it is a non-convergent Palais-Smale sequence, but the
    energy $E(x,0)\equiv a$ is not converging to $k$. 
    In fact the energy level $E^{-1}\{k\}$ could be
    a regular energy level, above the critical value,
    which has no invariant measure as that of theorem~\ref{NPS}.
   
  \bigskip
 \begin{Lemma}\label{RelComp}
 
 For all $k\in\re$,
 if $0\notin[a,b]\subset\re$, then the gradient flow of $-\cA_k$ on
 $\La_M$ or $\Om_M(q_0,q_1)$ is relatively complete on $[a\le\cA_k\le b]$.
 \end{Lemma}

 \begin{proof}
 Write $f=\cA_k:\La_M\to\re$ [resp.  $f=\cA_k:\Om_M(q_0,q_1)\to\re$]
 and let $\psi_t$ be the flow of $Y=-\nabla f$. Then
 \begin{equation}\label{DifF}
 f(\psi_{t_1}(p))-f(\psi_{t_2}(p))
 =-\int_{t_1}^{t_2} \nabla f(\psi_s(p))\cdot Y(\psi_s(p))\;ds
 =\int_{t_1}^{t_2} \lV Y(\psi_s(p))\rV^2\;ds.
 \end{equation}
 Moreover, using the Cauchy-Schwartz inequality, we have that
 $$
 d(\psi_{t_1}(p),\psi_{t_2}(p))^2
 \le\Big[\int_{t_1}^{t_2}\lV Y(\psi_s(p))\rV\;ds\Big]^2
 \le|t_2-t_1|\;\int_{t_1}^{t_2}\lV Y(\psi_s(p))\rV^2\;ds.
 $$
 Thus
 \begin{equation}\label{e.d2}
 d(\psi_{t_1}(p),\psi_{t_2}(p))^2 
 \le |t_2-t_1|\; \;|\,f(\psi_{t_1}(p))-f(\psi_{t_2}(p))\,|.
 \end{equation}
 Let $I=[0,\a[$ a maximal interval of definition of
 $t\mapsto \psi_t(p)$. Suppose that
 $a\le f(\psi_t(p))\le b$ for $0\le t< \a<\infty$.
 By inequality~\eqref{e.d2}, for any sequence $s_n\uparrow\a$ we have
 that $n\mapsto \psi_{s_n}(p)=\big(x_{s_n},T(s_n)\big)$ 
 is a Cauchy sequence in $\La_M\cap[a\le f\le b]$ 
 (resp. in $\Om_M(q_0,q_1)\cap[a\le f\le b]$).
 Then
 $$
 T_0=\lim_{s\uparrow\a}T(s)\in[0+\infty[ \quad\text{ exists.}
 $$
 
 If $0<T_0<+\infty$, since all such $\{\psi_{s_n}(p)\}_n$ are Cauchy
 sequences, then $q=\lim_{s\uparrow\a}\psi_s(p)=\psi_\a(p)$ exists.
  Since $f$ is $C^1$, we can extend the solution 
 $t\mapsto\psi_t(p)$ at $t=\a$. This contradicts 
 the definition of $\a$.
 
 If $T_0=0$, then there is a sequence $s_n\uparrow\a$ such that
 $$
 \tfrac{d}{ds}T(s_n)\le 0.
 $$
 Since $L$ is quadratic at infinity, there exist constants
 $a_0$,~$a_1$,~$b_0$,~$b_1$,~$c_0$,~$c_1>0$ such that
 \begin{align}
 c_0\,|v|_x^2-c_1\le L(x,v)+k&\le b_0\,|v|_x^2+b_1, 
 \label{c0c1b0b1}\\
 E(x,v)&\ge a_0\,|v|_x^2-a_1.
 \notag
 \end{align}
 Write $T_n:=T(s_n)$, $y_n(t):=x_{s_n}(t/T_n)$. Then
 \begin{align*}
 0\ge \tfrac{d\;}{ds}T(s_n)
 =-\frac{\partial \cA_k}{\partial T}
 =-k+\frac 1{T_n}\int_0^{T_n}E(\yn,\dyn)\;dt
 \ge -k-a_1+\frac{a_0}{T_n}\int_0^{T_n}\lv\dyn\rv^2\;dt.
 \end{align*}
 Since $\lim_nT_n=0$, this implies that
 $$
 \lim_n\int_0^{T_n}\lv\dyn\rv^2\;dt=0.
 $$
 Also, from~\eqref{c0c1b0b1}
 $$
 c_0\int_0^{T_n}\lv\dyn\rv^2\;dt-c_1\,T_n
 \le\cA_k(x_{s_n},T_n)
 \le b_0\int_0^{T_n}|\dyn|^2\;dt+b_1\,T_n.
 $$
 Hence $\lim_n\cA_k(x_{s_n},T_n)=0$. This contradicts the hypothesis
 $\cA_k(x_{s_n},T_n)\in[a,b]\not\ni 0$.
 
 \end{proof}

 \medskip

 \begin{Corollary}\label{RelCompQ}

 For all $k\in \re$,
 if $q_0, q_1\in M$, $q_0\ne q_1$ and $b\in\re$, 
 then the gradient flow of $-\cA_k$ on $\Om_M(q_0,q_1)$ 
 is relatively complete  on $[\cA_k\le b\,]$.
 \end{Corollary}

 \begin{proof}
 By Lemma~\ref{Om0}, if $q_0\ne q_1$ then
 $$
 \inf\big\{\,T>0\,
  \big\vert\,(x,T)\in\Om_M(q_0,q_1)\cap[\cA_k\le b] \,
   \big\}>0.
 $$
 Therefore the case $T_0=0$ in the proof of Lemma~\ref{RelComp}
 does not happen.
 
 \end{proof}

 \begin{Corollary}\label{RelCompL}\quad

 Let $\La_1$ be a connected component of $\La_M$ or $\Om_M(q_0,q_0)$
  with non-trivial homotopy class.
 For all $k,b\in\re$,
 the gradient flow of $-\cA_k$ on $\La_1$ 
 is relatively complete  on $[\cA_k\le b\,]$.
 \end{Corollary}

  \begin{proof}
  By Lemma~\ref{La1}, the case $T_0=0$ in 
  the proof of Lemma~\ref{RelComp}
  does not happen.
  
  \end{proof}

 \bigskip

 \section{Generic Palais-Smale condition for the Mountain-Pass
 geometry.}\label{sPSMPG}

\bigskip

  \begin{Proposition}\label{GPS}

 Given $\ga_0,\, \ga_1\in \La_M$ ~[resp. $\Om_M(q_0,q_1)$] let
 $$
 \cC(\ga_0,\ga_1):=\big\{\,\Ga:[0,1]\to\La_M \;[\text{resp.
 }\Om_M(q_0,q_1)\,]\;\big\vert\; \Ga \text{ is continuous},
 \Ga(0)=\ga_0,\; \Ga(1)=\ga_1\;\big\}.
 $$
 For $k\in\re$, let
 $$
 c(k):=\inf_{\Ga\in\cC(\ga_0,\ga_1)}\;\;\max_{s\in[0,1]}\;\cA_k(\Ga(s)).
 $$
 Suppose that for some $k_0\in\re$ we have that
 $c(k_0)\ne 0$ and
 $$
 c(k_0)> \max \big\{\,\cA_{k_0}(\ga_0),\; \cA_{k_0}(\ga_1)\,\big\}
 $$
 Then there exists $\e>0$ such that for Lebesgue almost every
 $k\in]k_0,k_0+\e[$, $c(k)$ is a critical value for $\cA_k$ on 
 $\La_M$ [resp. $\Om_M(q_0,q_1)$], with a critical point which is not
 a strict local minimizer.

 \end{Proposition}

 \begin{proof}
 
 Observe that for all $k\in\re$, the number $c(k)$ is finite and
 $c(k)\ge \max\{\,\cA_k(\ga_0),\,\cA_k(\ga_1)\,\}$.
 Since for all $\ga\in \La_M$ [resp. $\Om_M(q_0,q_1)$], the function
 $k\mapsto \cA_k(\ga)$ is non-decreasing, then $k\mapsto c(k)$ is
 non-decreasing. By the continuity of $\cA_k$ on $k$, the
 functions $k\mapsto \cA_k(\ga_0)$ and $k\mapsto \cA_k(\ga_1)$ are
 continuous. Let $\e>0$ be such that\footnote{Observe that
 the function $k\mapsto c(k)$ may be discontinuous, in particular,
 we allow $k_0$ to be a discontinuity point for $k\mapsto c(k)$.}
  \begin{equation}\label{MPe}
 \max\big\{\, \cA_k(\ga_0),\,\cA_k(\ga_1)\,\big\}
 <c(k_0)\le c(k) \ne 0
 \qquad
 \text{ for all } k_0<k<k_0+\e.
 \end{equation}

 Write $I_\e:=]k_0,k_0+\e[$. Since the function $c:I_\e\to\re$ is
 non-decreasing,  by Lebesgue's Theorem there is a total measure
 subset of $I_\e$ where $c(\cdot)$ is locally Lipschitz, i.e. the subset
 $$
 \cK:=\big\{\;k\in I_\e\;\big\vert\;\exists M>0,\; \exists \de_0>0,\;
 \;\forall |\de|<\de_0:\;\;|c(k+\de)-c(k)|<M\;|\de|\;\big\}
 $$
 has total Lebesgue measure in $I_\e$.
 
 Now fix $\uk\in \cK$ and a sequence $k_n\ge \uk$ with $\lim_n
 k_n=\uk$. By the definition of $\e$ in~\eqref{MPe} the functionals
 $\cA_{k_n}$, $\cA_{\uk}$ both show a mountain pass geometry with
 the same set of paths $\cC(\ga_0,\ga_1)$. 
 
 Let $\Ga_n\in\cC(\ga_0,\ga_1)$ be a path such that
 \begin{equation}\label{Gan}
 \max_{s\in[0,1]} \cA_{k_n}(\Ga_n(s))\le c(k_n)+(k_n-\uk).
 \end{equation}

 Let $M=M(\uk)>0$ be given by the property $\uk\in\cK$.
 Let $B\subset A\subset \La_M$ [resp. $\Om_M(q_0,q_1)$] be the
 closed subsets defined by
 \begin{align*}
 B:&=\{\, (x,T)\in\La_M \;[\text{resp. }\Om_M(q_0,q_1)]\;|\;
 T\le M+2\,\},
 \\
 A:&=\{\, (x,T)\in\La_M \;[\text{resp. }\Om_M(q_0,q_1)]\;|\;
 T\le M+3\,\}.
 \end{align*}
 Then $A$ contains the $\tfrac 12$-neighbourhood of $A$ in $\La_M$
 [resp. $\Om_M(q_0,q_1)$].
 Since from~\eqref{MPe}, $c(\uk)\ne 0$,
 by Propositions~\ref{Tfinito} and~\ref{La0}.\eqref{La0i},
 the functional $\cA_\uk$ satisfies the Palais-Smale condition
 restricted to $A$, at level $c(\uk)$.

By the choice $\uk\in\cK$, the function $k\mapsto c(k)$ is continuous
 at $\uk$.
 Since $k\mapsto \cAk(\ga)$ is increasing, 
 \begin{equation*}\label{Gan2}
 \max_{s\in[0,1]} \cA_{\uk}(\Ga_n(s))\le
 \max_{s\in[0,1]} \cA_{k_n}(\Ga_n(s))\le
  c(k_n)+(k_n-\uk) \overset{n}{\longrightarrow} c(\uk).
 \end{equation*}
 
  If $s\in[0,1]$ is such that
  $$
  \cA_{\uk}(\Ga_n(s))>c(\uk)-(k_n-\uk),
  $$
  then $\Ga_n(s)=(x,T)$ with
  $$
  T=\frac{\cA_{k_n}(\Ga_n(s))-\cA_\uk(\Ga_n(s))}{k_n-\uk}
  \le\frac{c(k_n)-c(\uk)}{k_n-\uk}+2
  \le M(\uk)+2,
  $$
  if $n$ is large enough.
 
  Given $\de>0$, let $n$ be large enough so that
  \begin{gather*}
  c(k_n)-c(\uk)+(k_n-\uk)<\de,
  \\
  0<\la_n:=(k_n-k)<\de.
  \end{gather*}
  Then
  $$
  \Ga_n([0,1])\subset \big( B\cup [\cA_{\uk}\le c(\uk)-\la_n]\big)
  \cap [\cA_{\uk}\le c(\uk)+\de].
  $$
  Since $c(\uk)\ne 0$, by Lemma~\ref{RelComp} the gradient flow of
  $-\cAk$ is relatively complete on $[c(\uk)-\e,c(\uk)+\e]$
  for some $\e>0$.
  Now Corollary~\ref{MountainPass} implies that
  $\cA_\uk$ has a critical point in $A$ which is not
  a strict local minimizer.
  \end{proof}

  \bigskip
  
  \section{The displacement energy.}
  
   Write $I:=[0,1]$. Given a subset $A\subseteq T^*M$ 
   let $\cH_c(I\times A)$ be the set of smooth functions 
   $H:I\times T^*M\to\re$ whose support is compact 
   and contained in $I\times A$. To such $H\in \cH_c(I\times A)$ 
   we associate its Hamiltonian vector field $X_{H_t}$,
   defined by $\om(X_{H_t},\cdot\,)=-dH_t(\,\cdot\,)$, where 
   $\om=dp\wedge dx$ and its corresponding Hamiltonian
   flow $h_t$. The set of functions in $\cH_c(I\times A)$ 
   which do not depend on $t\in I$ is denoted by $\cH_c(A)$.
  
   We say that $F\in\cH_c(A)$ is {\it slow} if all non-constant
   contractible (in $T^*M$) periodic orbits of its hamiltonian
   flow $f_t$ have period $>1$.
   Define the $\pi_1$-sensitive Hofer-Zehnder capacity of $A$
   by
   $$
   c_{HZ}^\circ(A,T^*M,\om):=\sup\,\{\,\max F\;|\;F\in\cH_c(\interior A)
   \;\text{ is slow}\;\}.
   $$
   The equivalence of this definition of the Hofer-Zehnder capacity
   with the original definition in~\cite{HZ} is proven in theorem~2.9
   of~\cite{GiGu}.
  
   Given $H\in\cH_c(I\times A)$ define its norm $\lV  H \rV$  
   as
   $$
   \lV H\rV:= \int _0^1\left(\sup_{z\in A}H(t,z)-\inf_{z\in A} H(t,z)
   \right)\;dt.
   $$
   The displacement energy $e(A,T^*M,\om)$ of a {\it compact} subset 
   $A\subseteq T^*M$ is defined as
   $$
   e(A,T^*M,\om):=\inf\{\,\lV H\rV\;|\; H\in\cH_c(I\times T^*M),\;
   h_1(A)\cap A=\emptyset\;\},
   $$
   where $h_1$ is the time 1 map of the hamiltonian flow of $H$.
   
   \begin{Lemma}
   Given an open subset $U\subset M$ there is a smooth
   function $\phi:M\to\re$  whose critical 
   points are all in $U$.
   \end{Lemma}
   
   \begin{proof}
   Let $f:M\to\re$ be a Morse function. Its set of critical points
   $C(f)=\{ \,x_1,\ldots, x_N\,\}$ is finite.
   Let $\{\ga_i\}_{i=1}^N$ be a collection of disjoint
   smooth curves $\ga_i:[0,1]\to M$ such that 
   $\ga_i(0)=x_i$ and $\ga_i(1)\in U$. Let $\{B_i\}_{i=1}^N$
   be a collection of disjoint tubular neighbourhoods of 
   the curves $\ga_i$. For each $i$, let $h_i$ be a 
   smooth diffeomorphism of $M$ with support in $B_i$
   such that $h_i(\ga_i(1))=x_i$. Now let
   $\phi= f\circ h_1\circ \cdots\circ h_N$.
   \end{proof}
   
   \begin{Proposition}\label{displace}
   If $k<e_0(L)$ then $e([E\le k],T^*M,\om)<+\infty$.
   \end{Proposition}
   
   \begin{proof}
    Since $k<e_0(L)$ then $U:=M\setminus \pi([E\le k])$
    is a non-empty open subset of $M$. Let $\phi:M\to\re$ be a smooth
    function such that all its critical points are in $U$.
    Let $G:=\nabla \phi$ be the gradient vector field of $\phi$
    and $g_t$ its gradient flow. The $\om$-limit of 
    every orbit of $g_t$ is a critical point of $\phi$
    which is inside $U$. Since $M\setminus U$ is compact,
    there is a finite time $T>0$ such that 
    $g_T(M\setminus U)\subset U$. Let $F= T\cdot G$, then
    its time-1-flow $f_1=g_T$ satisfies $f_1(M\setminus U)\subset U$.
    
    Let $R:T^*M\to \re$ be $R(x,p):=\langle p, F(x)\rangle_x$.
    The Hamiltonian equations for $R$ are 
    \begin{align*}
    \dx &= \phantom{-} \nabla_pR(x,p) = F(x), \\
    \dot{p} &= -\nabla_xR(x,p) = - p \cdot D_xF.
    \end{align*}
    In particular, the Hamiltonian flow $r_t$ of $R$
    lifts the flow $f_t$.    
     Hence
     $$
     r_1([E\le k])\subset r_1(\pi^{-1}(M\setminus U))
     \subset \pi^{-1}(U).
     $$
     Let 
     $$
     A:= 1+ \sup\{\;\lv r_t(x,p)\rv\;|\; 
     H(x,p)\le k,\;t\in[0,1]\;\}.
     $$
     Let $\la:\re\to[0,+\infty[$ be a smooth function such that
     $\la(r)\equiv 1$ if $|r|\le A$ and $\la(r)=0$ if $|r|\ge A+1$.
     Now let $S:T^*M\to\re$ be $S(x,p):=\la(|p|_x)\;\langle
     p,F(x)\rangle_x$. We have that $S$ has compact support and 
     its flow $s_t$ satisfies $s_t(x,p)=r_t(x,p)$ when 
     $t\in[0,1]$ and $H(x,p)\le k$. In particular
     $s_1([H\le k])\cap [H\le k]=\emptyset$ because
     $s_1([H\le k])\subset \pi^{-1}(U)$.
     Therefore
     $e([H\le k],T^*M,\om)\le \lV S\rV<+\infty$.
    
   \end{proof}

    \begin{Corollary}\label{k<e0}
    
    For Lebesgue almost every $k<e_0(L)$ the
    energy level $[E=k]$ has a periodic orbit 
    which is contractible in $M$.
    \end{Corollary}
   
    \begin{proof}
    By theorem~1.3 in~\cite{Schlenk} (also~\cite{FraSch}),
    $c^\circ_{HZ}(A,T^*M,\om)\le 4\;e(A,T^*M,\om)$.
    From Proposition~\ref{displace}, we get that
    $c^\circ_{HZ}([E\le k],T^*M,\om)<+\infty$
    for all $k<e_0(L)$. A standard argument using the
    Hofer-Zehnder  capacity ~\cite[p. 118\,--119]{HZ}
    shows that almost all energy levels $[E=k]$, $k<e_0(L)$
    have a periodic orbit which is contractible in $T^*M$
    but possibly non-contractible in $[E\le k]$. Since $T^*M$ retracts 
    to the zero section $M\times \{0\}$, the projection
    of the closed orbit to $M$ is contractible in $M$.
    \end{proof} 
   
  \bigskip
  
  \section{Loops, closed orbits and conjugate points.}
  \label{sLC}

  {\bf Proof of Theorem \ref{LC}:}
  
  \eqref{LC.c}. We first prove that for all $k>c_u(L)$, 
       $E^{-1}\{k\}$ contains a periodic orbit.
       By Corollary~\ref{SPS}, $\cA_k$
       satisfies the Palais-Smale condition. 

       If $\pi_1(M)\ne 0$, by Lemma~\ref{l27}, $\cA_k$
       is bounded below on each non-trivial free homotopy class
       $\si\in[S^1,M]$. Let $\La_\si$ be the connected component of $\La_M$
       corresponding to $\si$. 
       By Corollary~\ref{RelCompL} the gradient flow
       of $-\cA_k$ is relatively complete on $[\cA_k\le b]\cap\La_\si$
       for any $b\in\re$. By corollary~\ref{minimum}
       there is a minimizer of $\cA_k$ on $\La_\si$.

       If $\pi_1(M)=0$, then $c_u(L)=c_0(L)=c(L)$ and $k>c(L)$. 
       Since $M$ is closed, there is some non-trivial homotopy group 
       $\pi_\ell(M)\ne 0$. Choose a non-trivial free homotopy class
       $0\ne \sigma\in[S^\ell,M]$. A map $f:S^{\ell}\to M$ with homotopy
       class $\sigma$ can be seen as a family $F$ of closed curves
       in $M$ (see e.g.~\cite[page 37]{Kli}). Let $\cF$ the set of all
       such families corresponding to the homotopy class $\sigma$. 
       Clearly $\cF$ is a forward invariant family.
       Since the homotopy class $\sigma$ is non-trivial,
       (c.f.~\cite[Th. 2.1.8, page 37]{Kli}):
       $$
       \inf_{F\in\cF}\sup_{(x,T)\in F}\length(x) =:a>0.
       $$
       By the superlinearity of $L$ there is $b>0$ such that
       $L(x,v)> \lv v\rv^2_x-2b$ for all $(x,v)\in TM$.
       We can assume that $b\gg k$.
       If $(x,T)\in \La_M$ is a closed curve  with length 
       $\ell\ge a$, bounded action $\cA_{k}(x,T)\le
       \a$ and speed $|v|$, then
       \begin{gather*}
       \ell^2=\left[\int_0^T |v|\right]^2
             \le T\;\int_0^T|v|^2,
       \\
       \a\ge \cA_{k}(x,T)\ge\int_0^T|v|^2-2b\,T+k\,T
       \ge\frac{\ell^2}T-(2b-k)\,T.
       \end{gather*}
     Hence $(2b-k)\, T^2+\a\,T-\ell^2\ge 0$. Since $T\ge 0$
     and $\ell^2\ge a^2$, we have that
     $$
     T\ge \frac{ -\a+\sqrt{\a^2+4\,(2b-k)\,a^2}}{2\,(2b-k)}=:d>0.
     $$
     And then
     \begin{align*}
     \cA_k(x,T)&\ge \cA_{c(L)}(x,T)+[k-c(L)]\, T
     \\
     &\ge 0+[k-c(L)]\,d >0.
     \end{align*}
     Thus
     $$
     c(\cF) := \inf_{F\in\cF}\sup_{(x,T)\in F}\cA_k(x,T)\ge
     [k-c(L)]\,d >0.
     $$
     Since $c(\cF)\ne 0$, Corollary~\ref{SPS}, Lemma~\ref{RelComp} 
     and Proposition~\ref{minimax} with
     $B=A=X=\La_M$ imply
     that there is a critical point on $\La_M$ for $\cA_k$.

       By definition of $c_u(L)$, if $e_0(L)<k<c_u(L)$ 
       then there is a closed curve
       $(x_1,T_1)\in\La_M$ homotopic to a point,
       such that $\cA_k(x_1,T_1)<0$.
       Then Proposition~\ref{CritVal}.\eqref{CVper}
       and  Proposition~\ref{GPS} imply that 
       for almost every $k\in]e_0,c_u[$ there is
       a critical point for $\cA_k$ in $\La_M$ with $c(k)>0$:
       i.e. a periodic orbit with trivial homotopy class and positive
       $(L+k)$-action which is not a strict local minimizer.

    The case $k<e_0(L)$ is proven in Corollary~\ref{k<e0}.
    The closed orbit obtained in Corollary~\ref{k<e0} could be
    a singularity of the Euler-Lagrange flow. But in that
    case $k$ is a critical value of the energy function.
    By Sard's theorem that can only happen on a set of measure
    zero of values of $k$.

    \eqref{LC.l}. For item~\eqref{LC.l} and $k>c_u(L)$ 
    the proofs are similar to those of item~\eqref{LC.c} 
    working on $\Om_M(q_0,q_0)$. Namely, if $\pi_1(M,q_0)\ne 0$
    one finds a minimizing loop in a non-trivial homotopy 
    class.
        If $\pi_1(M,q_0)= 0$ we decompose a map
    $(S^\ell,\text{N.Pole})\to (M,q_0)$ in a non-trivial
    homotopy class of $\pi_\ell(M,q_0)$ into a family of closed loops in
    $\Om_M(q_0,q_0)$.    
    
     For $E(q_0,0)<k<c_u(L)$, item~\eqref{LC.l} follows from
    Proposition~\ref{CritVal}.\eqref{CVloop}, 
    and Proposition~\ref{GPS} similarly to item~\eqref{LC.c}.

    \eqref{LC.ps}.  If $k<c_u$ and the Palais-Smale condition
	  holds, the proof is similar to items~\eqref{LC.c}
	  and~\eqref{LC.l}, but now using 
	  Corollary~\ref{MountainPass}, with $B=A=X$ instead of
	  Proposition~\ref{GPS}.
 
 \qed

 \bigskip

  Now we will prove Theorem~\ref{CP}.
  Let $H:T^*M\to\re$ be the hamiltonian associated to $L$ 
 and $\psi_t$ its hamiltonian flow.
 Recall that two points 
 $\theta_1, \theta_2$ are said {\it conjugate }
 if there is $T\in\re\setminus\{0\}$ such that
 $$
 \theta_2=\psi_T(\theta_1) \quad\text{ and }\quad
 d_{\theta_1}\psi_T(V(\theta_1))\cap V(\theta_2)\ne\{0\},
 $$
 where $V(\theta):= \ker d_\theta\pi\subset T_\theta(TM)$
 is the {\it vertical subspace } and $\pi:TM\to M$ 
 is the projection. This definition coincides with
 the one given in page~\pageref{conjugate} because
 the Legendre transform $\cL(x,v)=(x,L_v)$
 maps the vertical subspace
 of $T_{v_x}TM$ to the vertical subspace of $T_{\cL(v_x)}T^*M$.
  
 \bigskip

 \begin{Proposition}\label{SLM}\quad
 
 Suppose that the forward orbit of $(x_0,v_0)$ has no conjugate points.
 Let $\ga:[\e,T]\to M$ be the solution
 $\ga(t)=\pi(\vr_t(x_0,v_0))$. Let $x_{\e,T}(s):=\ga\big(\e+s(T-\e)\big)$,
 $s\in[0,1]$ and $k:=E(x_0,v_0)$. Then for all $T>\e>0$ the solution
 $(x_{\e,T},T-\e)$, is a strict local minimizer 
  of the free time $(L+k)$-action $\cA_k$ 
  on $\Om_M(\ga(\e),\ga(T))$.

 \end{Proposition}

 \bigskip

 In~\cite[page 663]{CIPP2} we gave an example of
 an orbit segment $(x,T=\pi)$ which has no conjugate points
 and which is not a local minimizer of $\cA_k$ on 
 $\Om_M(x(0),x(1))$. But in that example the forward
 orbit of $(x(0),\tfrac{\dx(0)}T)$ has a conjugate point at 
 time $t=2\pi$ 
 (cf. it is the same lagrangian as in example~A.3, p.~949 in~\cite{CI}). 

 \begin{proof}
 Let $H:T^*M\to\re$ be the hamiltonian associated to $L$,
 $\psi_t$ its hamiltonian flow and $X$ its hamiltonian vector field.
 Let $\om=dp\wedge dx$ be the canonical symplectic form on $T^*M$.
 Given $\theta\in H^{-1}\{k\}$, let
 \begin{gather*}
 \Si:=H^{-1}\{k\},\quad 
 \Si_{\pi(\theta)}:=T_{\pi(\theta)}^*M\cap \Si,\quad
 V(\theta):=\ker d_\theta\pi=T_\theta (T^*_{\pi(\theta)}M)
 \quad \text{and}\\
  \La_{\pi(\theta)}:=\cup_{t>0}\;\psi_t(\Si_{\pi(\theta)}).
 \end{gather*}
 Then  
 $T_\theta\Si_{\pi(\theta)}=V(\theta)\cap T_\theta\Si$
  and 
 \begin{gather*}
 T_{\theta}\La_{\pi(\theta)}
 =(V(\theta)\cap T_\theta\Si)\oplus \langle X(\theta)\rangle
 =:W(\theta),
 \\
 T_{\psi_t(\theta)}\La_{\pi(\theta)}
 =d\psi_t\big(W(\theta)\big).
 \end{gather*}
 By definition $i_X\om=-dH$, hence $i_X\om\vert_{T\Si}\equiv 0$.
 Since the vertical subspace $V(\theta)$ is lagrangian, we get
 that $\La_{\pi(\theta)}$ is a invariant lagrangian submanifold
 of $T^*M$ inside the energy level $\Si$.
 Since $V=\ker d\pi$, the kernel of the projection 
 $d\pi\vert_{\La_{\pi(\theta)}}$ restricted to $\La_{\pi(\theta)}$
 is
 $$
 \ker d_{\psi_t(\theta)}\pi\vert_{\La_{\pi(\theta)}}
 = V(\psi_t(\theta)) \cap d\psi_t\big(W(\theta)\big). 
 $$
 By Proposition~1.16 and Remark~1.17 in~\cite{CI} if the {\it whole}
 forward orbit of $\theta$ has no conjugate points
 then $V(\psi_t(\theta)) \cap d\psi_t\big(W(\theta)\big)=\{0\}$
 for $t>0$, and hence the derivative of the projection 
  $d_{\psi_t(\theta)}\pi\vert_{\La_{\pi(\theta)}}$
 is injective along the forward orbit $\psi_t(\theta)$, $t>0$.
 If $T>\e>0$, then the projection $\pi\vert_{\La_{\pi(\theta)}}$ 
 is an immersion in  a small tubular 
 neighbourhood $N\subset \La_{\pi(\theta)}\subset H^{-1}\{k\}$  
 of the compact orbit segment 
 $\psi_{[\e,T]}(\theta)$.

 Now fix $\theta_0:=L_v(x_0,v_0)$.
 Observe that if the tubular neighbourhood $N$ is small enough, 
 then $N$ is either contractible 
 or $N$ is homeomorphic to a solid torus and the orbit of $\theta_0$
 is periodic with period smaller than or equal to $T-\e$.

 If $(x,p)\in N$ we have that
 $$
 k=H(x,p)=\sup_{v\in T_xM}p\cdot v-L(x,v).
 $$
 Since $N\subset H^{-1}\{k\}$, for any curve  $(z(s),q(s))$ inside $N$,
 \begin{equation*}\label{strictineq}
 q(s)\cdot \dz(s) \le L(z,\dz)+k,
 \end{equation*}
 with strict inequality if $L_v(z,\dz)\ne q\in N$.

  Now let $(y,S)\in\Om_M(\ga(\e),\ga(T))$  be a curve near
  $(x_{\e,T},T-\e)$ in the metric of $\cH^1(M)\times\re^+$.
  Since the  time parameters $S$, $T-\e$, are bounded,
  by Lemma~\ref{nearby}, if $(y,S)$ is sufficiently
  near $(x_{\e,T},T-\e)$ then the Hausdorff distance 
  $d_H\big(y([0,1]),x_{\e,T}([0,1])\big)$ is small.
  In particular, $y$ is homotopic to $x_{\e,T}$ with fixed endpoints
  and $y([0,1])\subset\pi(N)$.
  Let $z(t):=y(t S)$ and let $(z(s),q(s))$ be 
  the lift of $z$ to $N$ with $q(0)=\psi_\e(\theta_0)$.
  Then $(z,q)$ is homotopic in $N$ to the orbit segment 
  $\psi_{[\e,T]}(\theta_0)$ with fixed endpoints.
 Since $N$ is a lagrangian submanifold of $T^*M$, the Liouville
 1-form $p\,dx$ is closed on $N$. Then
 \begin{equation}\label{slm}
 \cA_k(x,T-\e)
 = \oint_{\ga}\big( L+k \big)
 =\oint_{\psi_t(\theta_0)}p\,dx
 =\oint_{(z,q)} p\,dx  \le \oint_{z}\big( L+k \big)
 =\cA_k(y,S),
 \end{equation}
 with strict inequality if $q(s)\ne L_v(z,\dz)$
 on a set of positive measure. Thus $\ga\vert_{[\e,T]}$ is 
 a local minimum of the $(L+k)$-action.

 We now see that $\ga\vert_{[\e,T]}$ is a strict local minimum
 of the $(L+k)$-action. 
 Let $\cL:TM\to T^*M$ be the Legendre transform 
 $\cL(x,v)=L_v(x,v)$. Observe that the hamiltonian vector field $X$
 satisfies
 $$
 d\pi\circ X(x,p) =\cL^{-1}(x,p) \quad \text{ for all } (x,p)\in T^*M.
 $$ 
 Suppose that~\eqref{slm} is an equality.
 Then $q(s)=L_v(z,\dz)\in N$ for almost every $s\in[0,S]$.
 Therefore 
 \begin{equation}\label{hamEL}
 \dz=\cL^{-1}(z,q)=d\pi(X(z,q)) \quad \text{ for almost every }s\in[0,S].
 \end{equation}
 Since $z(s)$ is continuous,  its lift $q(s)$ is continuous.
 Hence, from~\eqref{hamEL}, its derivative $\dz$ is 
 continuous.
 Then equation~\eqref{hamEL} says that the curve $z$ is an orbit of the
 projection of the hamiltonian vector field on $N$. Since $X$ 
 is tangent to $N$, $N$ is $\psi_t$-invariant and 
 the lift $(z,q)$ is unique, we have that $(z,q)$
 must be an orbit of $X$.
  Since $z(0)=\ga(\e)$ and $(z(0),q(0))=\psi_\e(\theta_0)$, 
 we have that $z(t)=\ga(t)$ for  all $t\in[\e,T]$.
 Since $z(S)=x_{\e,T}(1)=\ga(T)$,  either $S=T-\e$ 
 or $\theta_0$ is a periodic point and $\vert S-(T-\e)\vert$ is a 
 multiple of its period. Since we are assuming that $\vert S-(T-\e)\vert$ 
 is small,  $S=T-\e$.
 Therefore~\eqref{strictineq} is a strict inequality unless 
 $(z,S)\equiv (x_{\e,T},T-\e)$. 

 \end{proof}

 \medskip

\noindent{\bf Proof of Theorem~\ref{CP}:}

 Observe that the convexity of $L$  
 implies\footnote{e.g. Lemma~\ref{Lo}.} 
 that  $\min_{v\in T_xM}E(x,v)=E(x,0)$.
 If $k>e_m(L)$ then there is $x_0\in M$ such that
$k>E(x_0,0)$. Then Theorem~\ref{LC}.\eqref{LC.l} says 
that for almost every $k\in]e_m(L),c_u(L)[$ there is 
an orbit segment with energy $k$ which is not a strict
local minimizer of the action functional $\cA_k$.
Then Proposition~\ref{SLM} implies that the forward orbit
of the initial point of such orbit segment must have 
a conjugate point. From the definition of conjugate point
and the continuity of the derivative of the hamiltonian flow, 
it is easy to see that having a conjugate point is an
open condition. 

If for a specific $k\in]e_m(L),c_u(L)[$, the energy level
$E^{-1}\{k\}$ satisfies the Palais-Smale condition
then the same argument, now using Theorem~\ref{LC}.\eqref{LC.ps}
and Proposition~\ref{SLM}, implies that the energy level
$k$ has conjugate points.

 \qed

\bigskip

  \section{Proof of Proposition~\ref{kPS}.}\label{CTO}
  
   Proposition~\ref{kPS} follows from Corollary~\ref{SPS}
   and Lemmas~\ref{injPS} and~\ref{pinj} below.

   \bigskip

   Fix $k\in\re$ and let $\Si:=\Hk$ be the energy level. 
   Let $X$ be the hamiltonian vector field for $H$ and
   $\psi_t$ be its flow. Let $\pi:T^*M\to M$ be the projection,
   $\om=dp\wedge dx$ the canonical symplectic form on $T^*M$ and
   $\Theta=p\,dx$ the Liouville 1-form on $T^*M$.

   Given $\psi_t$-invariant Borel probability measure $\nu$
   supported on $\Si$,  the {\it Schwartzman asymptotic cycle}
   \label{schw}
   $\cS(\nu)\in H_1(\Si,\re)\approx H^1(\Si,\re)^*$ of $\nu$
   is defined by
   $$
   \langle \cS(\nu),[\eta]\rangle = \int_\Si \eta(X)\;d\nu
   $$
   for every closed 1-form $\eta$ on $\Si$.
    The map $(\pi|_\Si)_*:H_1(\Si,\re)\to H_1(M,\re)$ maps
   $\cS(\nu)$ to the homology class $\rho(\nu)$
   of $\nu$.

   \begin{Lemma}\label{injPS}\quad
      If $\Si$ is of contact type and 
   $\pi_*:H_1(\Si,\re)\to H_1(M,\re)$ is injective,
   then $\cA_k$ satisfies the Palais-Smale condition.
   \end{Lemma}
   
   \begin{proof}
   Let $\Theta=p\,dx$ be Liouville 1-form.
   Observe that on the energy level $\Si$:
   $$
   \Theta(X)=p\cdot H_p=v\cdot L_v=L+k.
   $$
   Suppose that $\cA_k$ does not satisfy the Palais-Smale 
   condition. Let $\mu$ be the measure given by Theorem~\ref{NPS}
   and let $\nu=\cL_*(\mu)$ be its push-forward under the
   Legendre transform $\cL(x,v)=L_v(x,v)$.
   Let $\la$ be a contact-type form on $\Si$.
   Since $\la(X)\ne 0$, $\la(X)$ has a
   single sign on each connected component 
   of $H^{-1}\{k\}$, in particular in the support
   of $\nu$.
   Since $d\la=\om=d\Theta$, the form $\eta:=\la-\Theta$
   is closed on $\Si$.   
   Since
   $\pi_*(\cS(\nu))=\rho(\mu)=0$ and $\pi_*$ is injective,
    $\cS(\nu)=0$. Then
    \begin{align}
    A_{L+k}(\mu)
    &=\int_\Si \Theta(X)\;d\nu + 0
    =\int_\Si \Theta(X)\;d\nu +\langle\cS(\nu),[\eta]\rangle
    \notag\\
    &=\int_\Si (\Theta+\eta)(X)\;d\nu
    =\int_\Si \la(X)\;d\nu\ne 0.
    \label{A+}
    \end{align}
    This contradicts Theorem~\ref{NPS}.
    
    \end{proof}
    
    \bigskip
    
    \begin{Lemma}\label{pinj}
       
    If $\dim M\ge 2$ and either 
    \begin{itemize}
    \item $M\ne\T^2$ or 
    \item $M=\T^2$ and $k<e_0$, 
    \end{itemize}then
    $\pi_*:H_1(\Si,\re)\to H_1(M,\re)$ is an isomorphism.
    \end{Lemma}

   In the following proof we shall use the lagrangian\footnote{The
   hamiltonian version, $\Si=\Hk$, may not contain the zero section in
   its interior} version
   $\Si=E^{-1}\{k\}$. Its intersections with the fibers of $TM$,
   $\Si\cap T_{x}M$ are convex subsets containing $(x,0)$ in its
   interior which are either homeomorphic to a sphere $S^{n-1}$ 
   or to a point, when $E(x,0)=k$.

    \begin{proof}
          
     Suppose first that $k>e_0(L)$ and $M\ne \T^2$.
     
 Since $k>e_0(L)$, the energy level $\Si:=\Ek$ is isomorphic to
 the unit tangent bundle of $M$ with the projection $\pi:\Si\to M$.
 If $M$ is orientable, the Lemma follows from an argument
 using the Gysin exact sequence, e.g.~\cite[Lemma~1.45]{Patb}.

 If $M$ is not orientable and $n=\dim M\ge 3$, from the exact homotopy
 sequence of the fiber bundle $\pi:\Si\to M$:
 $$
 0=\pi_1(S^{n-1})\overset{i_*}\longrightarrow \pi_1(\Si)
 \overset{\pi_*}\longrightarrow \pi_1(M)
 \longrightarrow \pi_0(S^{n-1})=0,
  $$ 
  we get  that $\pi_*:\pi_1(\Si)\to\pi_1(M)$ is an isomorphism.
  This implies that $\pi_*:H_1(\Si,\re)\to H_1(M,\re)$ is an isomorphism.
 
  If $M$ is not orientable and $\dim M=2$, from the homotopy sequence
  above we get an isomorphism $f:\pi_1(\Si)/_{\ima i_*}\to \pi_1(M)$.
  Let $h:\pi_1(M)\to H_1(M,\re)$ 
  and $k:\pi_1(\Si)\to H_1(\Si,\re)$ be the natural homomorphisms. 
  We show that $\ima i_*\subset \ker k$ and therefore,
  $k$ induces a homomorphism $k_1:\pi_1(\Si)/_{\ima i_*}\to H_1(\Si,\re)$.
  Indeed, 
  the fiber $F=S^1$ lies inside a Klein bottle $K$ 
  inside $\Si$, which is $\pi^{-1}(\ga)$, where $\ga$
  is a closed curve containing the base point
  $\pi(F)$ along which $M$ is non-orientable. 
  Then if $1_F$ is a generator of the fundamental group
  of the fiber $1_F\in\pi_1(F)=\Z$, its image $i_*(1_F)$
  has order at most 2 in $\pi_1(\Si)$. Hence $k(i_*(1_F))=0\in
  H_1(\Si,\re)$.  
  
  The following diagram commutes. There $k_1$ and $h$ are surjective
  and $f$ is an isomorphism.
  $$
  \begin{CD}
  \pi_1(\Si)/_{\ima i_*} @>f>\approx> \pi_1(M) \\
  @V k_1 VV @VV h V \\
  H_1(\Si,\re) @>\pi_*>> H_1(M,\re)
  \end{CD} 
   $$ 
   Then $\pi_*$ is surjective. Suppose that $\pi_*(a)=0$.
   Let $b\in \pi_1(\Si)/_{\ima i_*}$ be such that $k_1(b)=a$
   and let $c=f(b)$. Then $h(c)=0$. Hence $c$ is in the commutator
   subgroup of $\pi_1(M)$. Since $f$ is an isomorphism, 
   $b$ is in the commutator subgroup of $\pi_1(\Si)/_{\ima i_*}$.
   Therefore $a=k_1(b)=0$. Thus $\pi_*$ is injective.

    Now assume that $k<e_0(L)$.   
     Let  $B:=\pi(\Si)$ and let $E:=T^1_BM$ 
 be the restriction  of the unit tangent bundle to $B$.
 Let $\equiv$ be the equivalence relation on $E$
 defined by $(x,v)\sim (y,w)$ iff either
 $(x,v)=(y,w)$  or   $E(x,0)=k$ and $x=y$.
 Then the energy level $\Si$ is homeomorphic to
 $E/_\sim$, i.e. the one point compactification 
 of the fibers over the points $x$ with $E(x,0)=k$. 
 
 We can  assume that $B$ is connected, for the connected
 components of $\Si$ are in 1-1 correspondence with the connected
 components of $B$ under the projection $\pi$.
 
  We can also assume that there is $b_1\in B$ 
 such that $E(b_1,0)\ne k$. For, if not,
 then $E^{-1}\{k\}=\{\,(x,0)\,|\,x\in B\,\}$
 and the Lemma becomes trivial.

 Let $p:E\to B$ be the restriction of the projection
 of the unit tangent bundle and $f:E\to\Si= E/_\sim$
 the canonical projection. Then $p=\pi\circ f$.
 The homotopy exact sequence
 of the fibering $S^{n-1}\hookrightarrow E\to B$ gives
 $$
 \pi_1(S^{n-1})\overset{i_*}\longrightarrow \pi_1(E,e_1)
 \overset{p_*}\longrightarrow \pi_1(B,b_1)
 \longrightarrow \pi_0(S^{n-1})=0.
 $$
 Then $p_*$ is an epimorphism and it induces an
 isomorphism $g:\pi_1(E)/_{\ima i_*}\to \pi_1(B)$.
 We see that $\ima i_*\subset \ker f_*$, so that
 $f_*$ induces a homomorphism 
 $\hat{f}_*:\pi_1(E)/_{\ima i_*}\to\pi_1(\Si)$.
 Indeed, if $n\ge 3$ then $\pi_1(S^{n-1})=0$ and then
 $\ima i_*=0$. If $n=2$, let $1_F$
 be a generator of the fundamental group of the fiber
 $\pi_1(S^1)=\Z$. Since $k<e_0(L)$ there is a point $x_1\in B$
 such that $E(x_1,0)=k$.
 Let $\la$ be a curve in $B$ joining $b_1$ to $x_1$.
 The fiber bundle $E$ over the interval $\la$ is trivial
 $E\vert_\la \approx S^1\times [0,1]$. Observe that
 the inverse image $\pi^{-1}(\la)\subset\Sigma$
 has the topology of a cylinder with one of its
 boundary circles compactified to a point. 
 Hence it is homeomorphic to a 2-disc, and the class $1_F$
 is represented by its boundary circle.
 Hence $f_*( i_*(1_F))=0\in \pi_1(\Si)$.

 We prove that $\pi_*:\pi_1(\Si)\to\pi_1(B)$
 is an isomorphism. This implies the Lemma.
 Since $g=\pi_*\circ \hat{f}_*$, it follows 
 that $\hat{f}_*$ is injective
 and $\pi_*$ is surjective. In order to prove that
 $\pi_*$ is injective it is enough to prove that 
 $\hat{f}_*$ is onto.
 Since $f_*$ and $\hat{f}_*$ have the same image,
 it is enough to prove that $f_*$ is surjective.

 Let $\si_1=f(e_1)\in \Si$. 
 Let $\Ga:(S^1,1)\to (\Si,\si_1)$ be a loop in $\Si$ 
 based at $\si_1$. We want a preimage under $f_*$
 of the homotopy class of $\Ga$.
 If $\Ga(s)\ne 0$ for all $s\in S^1$,
 such preimage is the homotopy class
 of $\Psi_E(s)=\frac{\Ga(s)}{\lV \Ga(s)\rV}$.
 In general, the problem is that 
 such definition of $\Psi_E$ may have no
 continuous extension to the $s$ where $\Ga(s)=0$.
 Assume now that there is $s_0\in S^1$ such that
 $\Ga(s_0)=0$.

 Let $\ga=\pi\circ\Ga:(S^1,1)\to (B,b_1)$
 be the projection of $\Ga$.
  Let $C:=\{\,x\in M\,|\,E(x,0)=k\,\}$ and 
  $D:=\{\, t\in S^1\,|\,\ga(t)\in C\,\}$.  
  Then $D$ is a compact subset of 
  $S^1$ and 
   its complement is a union of open intervals
   $I_i$.
   Choose any continuous loop $\La_E:(S^1,1)\to(E,e_1)$
   such that $p\circ\La_E=\ga$.
  
 The pullback of the sphere bundle $E$ along each segment 
 $\ga(I_i)$ is trivial $(\ga\vert_{I_i})^*E\approx I_i\times S^{n-1}$.
 Then the inverse image $\pi^{-1}\{\ga(\ov{I_i})\}\subset\Si$
 of the {\it closed} segment $\ga(\ov{I_i})$ has the topology
 of a cylinder $S^{n-1}\times [0,1]$ with its boundary
 spheres $S^{n-1}\times\{0\}$ and $S^{n-1}\times\{1\}$
 compactified to two points $\{A,B\}$ or to a single point
 $A=B$.

 Both segments $f\circ\La_E\vert_{\ov{I_i}}$ and
 $\Ga\vert_{\ov{I_i}}$ must have the same endpoints $A$ and $B$.
 Hence they are homotopic with fixed endpoints inside
 $\pi^{-1}\{\ga(\ov{I_i})\}\subset\Si$.
 Joint all these homotopies for each interval $I_i$, to
 obtain a homotopy in $\Si$ between $f\circ\La_E$ and $\Ga$.
 Therefore $f_*([\La_E])=[\Ga]\in \pi_1(\Si)$.
 
 \end{proof}

 \newpage 
 
 \appendix
 

 \section{A non ergodic measure in Theorem~\ref{NPS}.}\label{Reeb}


  Consider the flat metric on the 2-torus $\T^2$.
  Let $X$ be a vector field with norm 1 on $\T^2$
  whose orbits form a Reeb foliation. Let $L:T\T^2\to\re$
  be the lagrangian
  $$
  L(x,v):= \tfrac 12 \,\lv v-X(x)\rv^2.
  $$
  Its Euler-Lagrange flow is the same as the exact
  magnetic flow with lagrangian
  $$
  L(x,v) -\tfrac 12 = \tfrac 12 \,|v|^2 -\eta_x(v),
  $$
  where $\eta_x(v)=\,\langle X(x),v\rangle$.
  It is easy to see from the definition of critical value
  that $c(L)=0$. The vector field $X$ has two closed orbits
  $\ga_1$ and $\ga_2$ with {\it opposite} homology classes.
  Since $L\ge 0$ the Euler-Lagrange flow has only two
  ergodic invariant measures $\mu_1$, $\mu_2$, 
  with zero $L$-action, corresponding
  to the periodic orbits $\Ga_i=(\ga_i,\dga_i)$.
  The unique invariant probability $\mu$ with 
  $A_L(\mu)=0$ and zero homology class is 
  $\mu=\tfrac 12\,\mu_1+\tfrac  12\, \mu_2$.
  It follows that $c_u(L)=c_0(L)=c(L)=0$.
  Let $\L:T\re^2\to\re$ be the lift of $L$
  to the universal cover $\re^2$ of $\T^2$.

\centerline{
\parbox{12cm}{
\includegraphics[scale=.6]{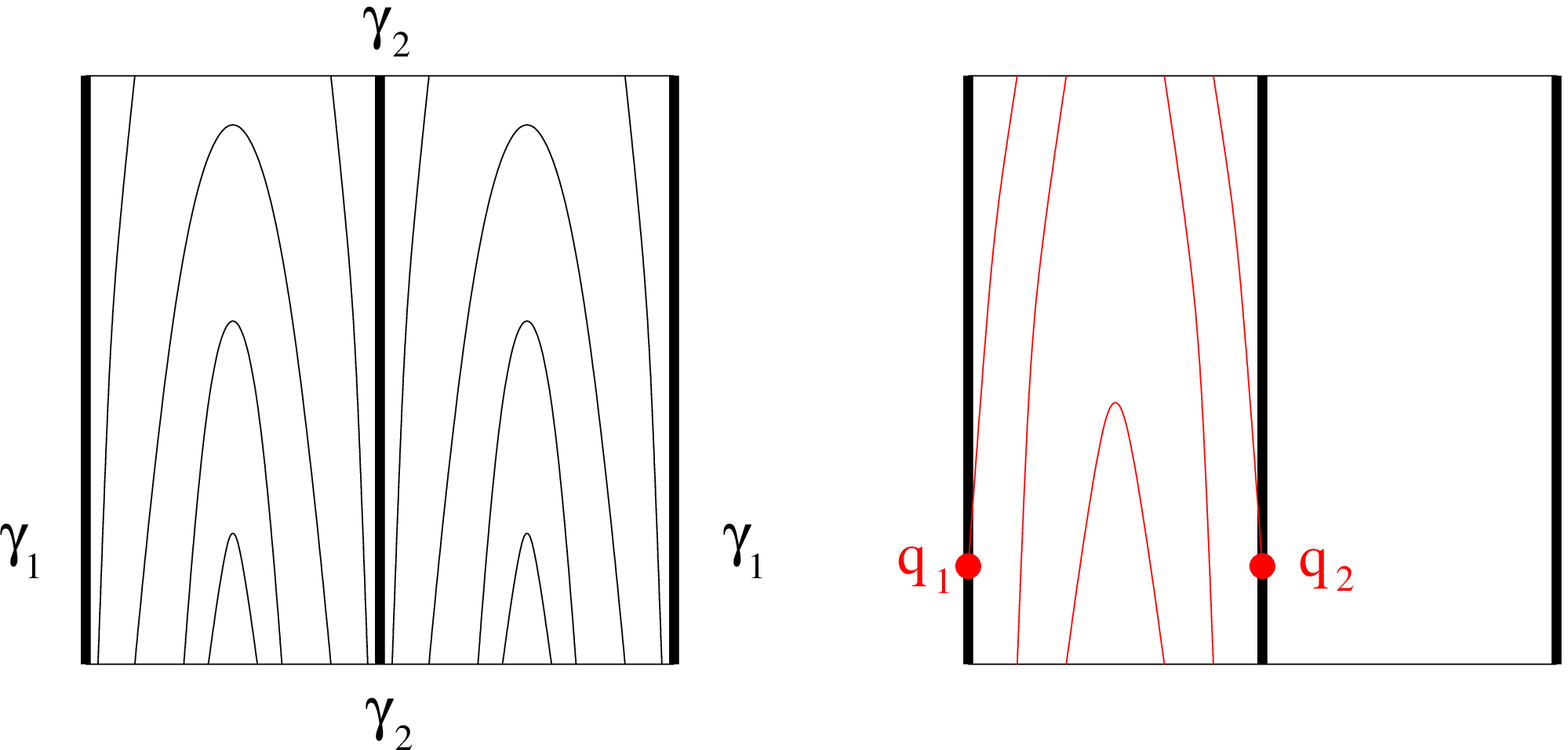}
\refstepcounter{figure}\label{reeb}
\\
{ \openup -2pt {\small {\sc Figure~\ref{reeb}.} 
The left figure shows the flow lines of the
vector field $X$ in the example.
The right figure shows a curve $(x_n,T_n)\in\Om_M(q_1,q_2)$
in an unbounded Palais-Smale sequence.
The probability measures associated to $(x_n,T_n)$
converge to the non-ergodic measure
$\tfrac 12\,\mu_1+\tfrac 12\,\mu_2$,
where $\mu_1$ and $\mu_2$ are the invariant
probabilities for the Euler-Lagrange flow corresponding
to the periodic solutions $\ga_1$ and $\ga_2$. 
}}
}
}
\hskip 1cm

\bigskip
  
  It is easy to see that the Peierls barrier
  for $\L$ is finite, because one can join two 
  points $q_0,\,q_1\in\re^2$ by curves with 
  bounded action which spend long time on a
  lift of $\ga_1$ and come back on a lift of
  $\ga_2$. Since the action of $\ga_1$ and $\ga_2$
  is zero, the total action spent on them is bounded.
  By  Corollary~\ref{SPS}, there is an unbounded
  Palais-Smale sequence $(x_n,T_n)$ with
  $\lim_n T_n=+\infty$. 
  Nevertheless, there is no {\it ergodic}
  invariant probability  in $TM$ with zero action and
  zero homology class.

 \newpage

 \section{Energy levels of non-contact type.}\label{nctp}

 The following theorem is not explicitly stated in~\cite{P9}.
 
 \begin{Theorem}[G. Paternain]\label{tacto}\quad
 Suppose that $\dim M\ge 2$.
 
 If $k>c_0(L)$ then $H^{-1}\{k\}$ is of contact type.
 
 If $M\ne\T^2$ and 
 $c_u(L)<k\le c_0(L)$ then $H^{-1}\{k\}$ is not of contact type.
 \end{Theorem} 
 
 There is an example in~\cite{clor} of a lagrangian in $\T^2$
 for which the energy level $E=c_0(L)$ is of contact type.
 
 As an application (c.f.~\cite[Th.~1.1]{P9}), if $M$ admits a metric with negative
 curvature and if the Lagrangian flow on
  an energy level $H^{-1}\{k\}$ with 
 $c_u(L)<k<c_0(L)$ is Anosov, then the strong stable and 
 unstable subbundles $E^{ss}$, $E^{uu}$ can not be $C^1$.
 For if they were $C^1$, the form $\la$ defined by
 $\la(X)\equiv 1$ and $\la\vert_{E^{ss}\oplus E^{uu}}\equiv 0$
 is a contact form for $H^{-1}\{k\}$ (c.f. U. Hamendst\"adt~\cite{Hams},
 G. Paternain~\cite[Th.~5.5]{P9}). Examples of such Anosov energy levels
 appear in G. Paternain \& M. Paternain~\cite{PP4}.

   If $\Si$ is a regular energy level,
   the {\it Liouville measure} $m$ on $\Si$ is the smooth measure induced
   by the volume form $i_Y \om^n$, where $Y$ is a vector field on
   $T^*M$ such that $\om(Y,X)\equiv 1$ on $\Si$. 
   It is invariant under the hamiltonian flow because
   $L_X (i_Y\om^n)=0$. We choose the orientation
   on $\Si$ that makes $m$ a positive measure.

  \begin{Lemma}\label{la+}\quad
  
  If $\la$ is a 1-form on $\Si$ such that
  $d\la=\om$ and $\la(X)\ne 0$ then $\la(X)>0$ on $\Si$.
  \end{Lemma}

  \begin{proof}

  Let $\xi$ be a 1-form on $\Si$. We have that
  \begin{align}
  \xi(X)\, i_Y \om^n &= i_X\left(\xi\wedge i_Y\om^n\right)
   +\xi\wedge\left(i_X\,i_Y\,\om^n\right)
   \notag\\
    &=0+\xi\wedge i_X\big[n\,(i_Y\,\om)\wedge\om^{n-1}\big]
   \notag
   &\text{because }\dim\Si=2n-1\\
   &=\xi\wedge \big(n\,\om(Y,X)\,\om^{n-1}-0\,\big)
   \notag
   &\text{because  on } \Si,\;i_X\om = -dH\equiv 0\\
   &=n\;\xi\wedge\om^{n-1}.
   \label{etaX}
  \end{align}
  
   We show first that the asymptotic cycle (c.f. page~\pageref{schw})
   of the Liouville measure
   $m$ is zero. Indeed, let $\Theta:= p\; dx$ be Liouville 1-form
   on $T^*M$ and $\tau:=\Theta\wedge \om^{n-2}$.
   If $\eta$ is a closed 1-form on $\Sigma$ then
   $\eta\wedge \om^{n-1}=\eta\wedge d\tau= d(\eta\wedge\tau)$.
   Hence by~\eqref{etaX},
   $$
   \int_\Si\eta(X)\;dm=\int_\Si\eta(X)\;i_Y\om^n
   =n\int_\Si\eta\wedge\om^{n-1} = n\int_{\Si} d(\eta\wedge\tau)
   =0.
   $$

  Since $\la(X)\ne 0$ on $\Si$, it is enough to prove that
  for any connected component $\cN$ of $\Si$ we have
  $\int_\cN\la(X)\;dm>0$.
  Since $d\la=\om=d\Theta$, the form $\eta=\la-\Theta$ is closed on
  $\Si$. Then
  $$
  \int_\cN\la(X)\;dm=\langle\cS(m),[\eta]\rangle+\int_\cN\Theta(X)\;dm
  =\int_\cN\Theta(X)\;dm.
  $$
   From~\eqref{etaX} we have that
  \begin{equation*} 
  \Theta(X)\;i_Y\om^n\big\vert_\Si=
  n\; \Theta\wedge\om^{n-1}\big\vert_\Si.
  \end{equation*}
  Let $\cW$ be the fiberwise convex hull of $\cN$ in $T^*M$. 
  Then $\partial\cW=\cN$.  By Stokes Theorem,
  $$
  \int_\cN\Theta(X)\;dm=
  \int_{\partial \cW}   \Theta(X)\, i_Y\om^n
  =  \int_{\partial \cW} n\; \Theta\wedge\om^{n-1}
  =n\int_\cW \om^n.
  $$
  
  We prove that the last integral is 
  positive.\footnote{Alternatively, 
  the integral is equal to the
  $(L+k)$ action of the Liouville measure $m_k$ on $E^{-1}\{k\}$.
  If $\ell>c_0(L)$ by~\eqref{c0}, $A_{L+\ell}(m_\ell)>0$,
  and one can show that  the
  orientation on $E^{-1}\{\ell\}$ defined by the Liouville 
  measure varies continuously  with $\ell$.}
  Since $dH(Y)=i_Y(-i_X\om)=1$,
  the convexity of $H$ implies that $Y$ is 
  an {\it outwards} pointing vector in
  $\Si=\partial \cW$. 
  A basis $(v_1,\ldots,v_{2n-1})$ of 
  $T\Si$ is positively oriented iff 
  $i_Y\om^n(v_1,\ldots,v_{2n-1})=\om^n(Y,v_1,\ldots,v_{2n-1})>0$
  and Stokes theorem uses $(Y,v_1,\ldots,v_{2n-1})$ 
  as a positively oriented basis for $T\cW$.
 
  \end{proof}
  
  \begin{Remark}\label{leb+}\quad
  
  Lemma~\ref{la+} also says that the Liouville measure has
  always positive $(L+k)$-action, for in the energy level
  $E^{-1}\{k\}$ we have that $\Theta(X)=v\cdot L_v=L+k$.
  \end{Remark}
   
  {\bf Proof of Theorem~\ref{tacto}:}
  
  From~\eqref{c0}, there is a closed 1-form $\eta$ on $M$
  such that $c_0(L)=c(L-\eta)$. The Hamiltonian of
  $L-\om$ is $\H(x,p)=H(x,p+\eta)$. If $k>c_0(L)$, 
  by Theorem~A in~\cite{CIPP} there is a smooth function
  $u:M\to\re$ such that $\H(x,d_xu)<k$ for all $x\in M$.
  From the definition~\eqref{H} of the Hamiltonian 
  $\H$ we have that
  $$
  L(x,v)-\eta_x(v)-d_xu(v)+k>0\quad\text{ for all }(x,v)\in  TM.
  $$
    Let $\Theta=p\, dx$ be the Liouville 1-form on $T^*M$.  
    Let $\la:=\Theta-\eta\circ d\pi-du\circ d\pi$, where
    $\pi:T^*M\to M$ is the projection.  Since $\eta$ is closed,
    $d\la=d\Theta=\om$. On $H^{-1}\{k\}$ we have that
    $$
    \Theta(X) = p\cdot H_p = L(x,v)+ k,
    $$
    where $v=H_p(x,p)$. Since $X=(H_p,*)$,
    on $H^{-1}\{k\}$ we have that
    $$
    \la(X)=L(x,v)+k-\eta_x(v)-d_xu(v)>0,
    $$
    where $v=H_p(x,p)$.
    
    Now assume that $c_u(L)<k<c_0(L)$.
    Let $\hL$ be the lift of $L$ to the abelian cover $\hM$.
    Since $k<c_0(L)=c(\hL)$, there exists a closed curve 
    $\hga$ in $\hM$ with negative $(L+k)$-action.
    Observe that the projection $\ga$ of $\hga$ to $M$
    has trivial homology class.
    The homotopy class of $\ga$ can not be trivial
    because if it where, its lift to the universal 
    cover would be closed and since $k>c_u(L)$ its
    $(L+k)$-action would be non-negative.
    Let $\si$ be the free homotopy class of  $\ga$ 
     and let
    $$
    \La_\si:=\{\,(x,T)\in\La_M\;|\; x\in\si,\; T>0\,\}. 
    $$
    Since $k>c_u(L)$, by lemma~\ref{l27},
    $$
    -\infty<\inf_{(x,T)\in\La_\si}\cA_k(x,T)
    <A_{\hL+k}(\hga)<0.
    $$
    Since $k>c_u(L)$, by Corollary~\ref{SPS}, 
    $\cA_k$ satisfies the Palais-Smale
    condition. By Corollary~\ref{RelCompL}
    the gradient flow of $-\cAk$ is relatively complete
    on $[\cAk\le 0]$. Then by Corollary~\ref{minimum}
    there is a  minimizer $(x,T)$ of $\cAk$ on $\La_\si$.
    The curve $y(t):=x(t/T)$ is a periodic orbit of
    the Euler-Lagrange flow with negative $(L+k)$-action
    and energy $k$. Let $\mu$ be the invariant probability
    measure supported on the periodic orbit $(y,\dy)$ and
    let $\nu=\cL_*(\mu)$, $\cL(x,v)=(x,L_v(x,v))$. 
    Since the homology class corresponding to $\si$
    is trivial, $\rho(\mu)=0$. By lemma~\ref{pinj},
    since $\pi_*(\cS(\nu))=\rho(\mu)=0$ and $M\ne\T^2$,
    $\cS(\nu)=0$.
    If $\la$ is a contact-type form, since $\supp(\mu)\subset
    E^{-1}\{k\}$, the same calculation as in~\eqref{A+} gives
    \begin{equation}\label{nctct}
    \int \la(X)\; d\nu = A_{L+k}(\mu) < 0.
    \end{equation}
    This contradicts lemma~\ref{la+}.

    When $k=c_0(L)$, by~\eqref{c0} there is an invariant
    probability $\mu$ such that $\rho(\mu)=0$ and $A_{L+c_0}(\mu)=0$.
    The same argument as in~\eqref{nctct} shows that $H^{-1}\{c_0\}$
    is not of contact type.

        \qed 
    
    
 \section{Non-magnetic Lagrangians.}\label{non-magnetic}
 
 The following result was suggested by R. Ma\~n\'e in~\cite{Ma4}.
 
   \begin{Theorem}
  If $L$ is a convex superlinear lagrangian
  and the  1-form $\theta_x:=L_v(x,0)$ is closed,
  then every energy level contains a closed orbit.
  \end{Theorem}

  Since we are looking for closed orbits on specified energy
  levels, by~\cite[prop. 18]{CIPP2} we can assume that
  $L$ is Riemannian at infinity. Since the 1-form 
  $\theta_x$ is closed the lagrangian $\L=L-\theta_x$
  has the same Euler-Lagrange flow as $L$. Replacing
  $L$ by $\L$ we can also assume that $\theta_x\equiv 0$.  
  If $k\ge e_0(L)$, by lemma~\ref{Lo},
  $$
  L(x,v)+k\ge\tfrac 12\,a_0\,|v|_x^2+\big[k-E(x,0)\big]
  \ge 0.
  $$
  Therefore $c(L)=e_0(L)$. By theorem~\ref{LC}, when 
  $k>e_0=c(L)$ the energy level $E^{-1}\{k\}$ has
  (non-trivial) closed orbits and at $k=e_0$ it has a singularity.

  By Proposition~\ref{displace} and Frauenfelder-Schlenk
  theorem~\cite{FraSch} for $k<e_0$ the set $[E\le k]$
  has finite Hofer-Zehnder capacity. Then~\cite[Th. 5, p.123]{HZ} 
  it has a closed orbit on every contact type energy level
  $E^{-1}\{k\}$ with $k<e_0$.
  
  Singular energy levels have a singularity of the Euler-Lagrange
  flow.  Now we see that the regular energy levels of $L$
  are of contact type. We use the following
  
  \begin{Proposition}\label{cct}
  [McDuff~\cite{mcduff}, also {\cite[sec. 2 \& app.  B]{clor}}]
  Suppose that the flow of a vector field $X$ on a compact orientable
  manifold $\Si$ does not admit a global cross section.
  Let $\Theta$ be a smooth 1-form on $\Si$. Then the following
  are equivalent
  \begin{enumerate}
  \item $\int\Theta(X)\;d\mu\ne 0$ for every invariant
        Borel probability with zero asymptotic cycle.
  \item There exist a smooth closed 1-form such that
        $\Theta(X)+\vr(X)$ never vanishes.	
  \end{enumerate}
  \end{Proposition}
  
  Let $k$ be a regular value of the energy function.
  Let $H$ be the hamiltonian of $L$ and $X$ its hamiltonian
  vector field.
  Since the Liouville measure has zero asymptotic
  cycle (inside Lemma~\ref{la+}) the energy level
  $\Si=H^{-1}\{k\}$ has no global cross section.
  For, by Poincar\'e duality such cross section 
  would give a closed 1-form $\eta$ such that 
  $
  \langle \cS(\mu),\eta\rangle =\int_\Si\eta(X)\;d\mu>0
  $
  for every invariant probability $\mu$.
  Since $E^{-1}\{k\}$ is a regular energy level, 
  if $(x,v)\in E^{-1}\{k\}$ then
  $k=E(x,v)>E(x,0)$. By Lemma~\ref{Lo}, writing $\Theta=p\,dx$,
  $$
  \Theta(X)=v\cdot L_v(x,v)
  =L(x,v)+k\ge\tfrac 12\,a_0\,|v|_x^2+\big[k-E(x,0)\big]
  > 0 \qquad\text{ if }\quad E(x,v)=k.
  $$
  By proposition~\ref{cct} there is a closed form $\vr$ on
  $\Si$ such that $\Theta(X)+\vr(X)\ne 0$. Let $\la=\Theta+\vr$. 
  Then $\la(X)\ne 0$ and $d\la=d\Theta=\om$.
  
  \bigskip
  
  Indeed, as above, from Lemma~\ref{la+} and proposition~\ref{cct}
  we have
  
  \begin{Proposition}\quad
  
  A compact energy level $E^{-1}\{k\}$ of a convex lagrangian
  is of contact type
  if and only if $A_{L+k}(\mu)>0$ for every Borel 
  invariant probability in $E^{-1}\{k\}$ with
  zero asymptotic cycle.
  \end{Proposition}
  
   
   
 \def\cprime{$'$} \def\cprime{$'$} \def\cprime{$'$}
 \providecommand{\bysame}{\leavevmode\hbox to3em{\hrulefill}\thinspace}
 \providecommand{\MR}{\relax\ifhmode\unskip\space\fi MR }
 \providecommand{\MRhref}[2]{%
   \href{http://www.ams.org/mathscinet-getitem?mr=#1}{#2}
 }
 \providecommand{\href}[2]{#2}

 \end{document}